\documentclass[twoside]{article}
\usepackage{tocstyle}

\usepackage[sc]{mathpazo}
\usepackage{graphicx,amsmath,amssymb,alltt}
\usepackage[T1]{fontenc} 
\usepackage{microtype}
\usepackage{pgfplots,tikz,wrapfig}
\usetikzlibrary{cd}
\usepackage{lscape}

\usepackage{wasysym}
\usepackage{fancyvrb,color}
\def\nicefrac#1#2{{}^{#1}\!\!/_{\!\!#2}}
\def\twov{2_{{}_{V}}}
\makeatletter
\def\PY@reset{\let\PY@it=\relax \let\PY@bf=\relax%
    \let\PY@ul=\relax \let\PY@tc=\relax%
    \let\PY@bc=\relax \let\PY@ff=\relax}
\def\PY@tok#1{\csname PY@tok@#1\endcsname}
\def\PY@toks#1+{\ifx\relax#1\empty\else%
    \PY@tok{#1}\expandafter\PY@toks\fi}
\def\PY@do#1{\PY@bc{\PY@tc{\PY@ul{%
    \PY@it{\PY@bf{\PY@ff{#1}}}}}}}
\def\PY#1#2{\PY@reset\PY@toks#1+\relax+\PY@do{#2}}

\expandafter\def\csname PY@tok@gd\endcsname{\def\PY@tc##1{\textcolor[rgb]{0.63,0.00,0.00}{##1}}}
\expandafter\def\csname PY@tok@gu\endcsname{\let\PY@bf=\textbf\def\PY@tc##1{\textcolor[rgb]{0.50,0.00,0.50}{##1}}}
\expandafter\def\csname PY@tok@gt\endcsname{\def\PY@tc##1{\textcolor[rgb]{0.00,0.27,0.87}{##1}}}
\expandafter\def\csname PY@tok@gs\endcsname{\let\PY@bf=\textbf}
\expandafter\def\csname PY@tok@gr\endcsname{\def\PY@tc##1{\textcolor[rgb]{1.00,0.00,0.00}{##1}}}
\expandafter\def\csname PY@tok@cm\endcsname{\let\PY@it=\textit\def\PY@tc##1{\textcolor[rgb]{0.25,0.50,0.50}{##1}}}
\expandafter\def\csname PY@tok@vg\endcsname{\def\PY@tc##1{\textcolor[rgb]{0.10,0.09,0.49}{##1}}}
\expandafter\def\csname PY@tok@vi\endcsname{\def\PY@tc##1{\textcolor[rgb]{0.10,0.09,0.49}{##1}}}
\expandafter\def\csname PY@tok@vm\endcsname{\def\PY@tc##1{\textcolor[rgb]{0.10,0.09,0.49}{##1}}}
\expandafter\def\csname PY@tok@mh\endcsname{\def\PY@tc##1{\textcolor[rgb]{0.40,0.40,0.40}{##1}}}
\expandafter\def\csname PY@tok@cs\endcsname{\let\PY@it=\textit\def\PY@tc##1{\textcolor[rgb]{0.25,0.50,0.50}{##1}}}
\expandafter\def\csname PY@tok@ge\endcsname{\let\PY@it=\textit}
\expandafter\def\csname PY@tok@vc\endcsname{\def\PY@tc##1{\textcolor[rgb]{0.10,0.09,0.49}{##1}}}
\expandafter\def\csname PY@tok@il\endcsname{\def\PY@tc##1{\textcolor[rgb]{0.40,0.40,0.40}{##1}}}
\expandafter\def\csname PY@tok@go\endcsname{\def\PY@tc##1{\textcolor[rgb]{0.53,0.53,0.53}{##1}}}
\expandafter\def\csname PY@tok@cp\endcsname{\def\PY@tc##1{\textcolor[rgb]{0.74,0.48,0.00}{##1}}}
\expandafter\def\csname PY@tok@gi\endcsname{\def\PY@tc##1{\textcolor[rgb]{0.00,0.63,0.00}{##1}}}
\expandafter\def\csname PY@tok@gh\endcsname{\let\PY@bf=\textbf\def\PY@tc##1{\textcolor[rgb]{0.00,0.00,0.50}{##1}}}
\expandafter\def\csname PY@tok@ni\endcsname{\let\PY@bf=\textbf\def\PY@tc##1{\textcolor[rgb]{0.60,0.60,0.60}{##1}}}
\expandafter\def\csname PY@tok@nl\endcsname{\def\PY@tc##1{\textcolor[rgb]{0.63,0.63,0.00}{##1}}}
\expandafter\def\csname PY@tok@nn\endcsname{\let\PY@bf=\textbf\def\PY@tc##1{\textcolor[rgb]{0.00,0.00,1.00}{##1}}}
\expandafter\def\csname PY@tok@no\endcsname{\def\PY@tc##1{\textcolor[rgb]{0.53,0.00,0.00}{##1}}}
\expandafter\def\csname PY@tok@na\endcsname{\def\PY@tc##1{\textcolor[rgb]{0.49,0.56,0.16}{##1}}}
\expandafter\def\csname PY@tok@nb\endcsname{\def\PY@tc##1{\textcolor[rgb]{0.00,0.50,0.00}{##1}}}
\expandafter\def\csname PY@tok@nc\endcsname{\let\PY@bf=\textbf\def\PY@tc##1{\textcolor[rgb]{0.00,0.00,1.00}{##1}}}
\expandafter\def\csname PY@tok@nd\endcsname{\def\PY@tc##1{\textcolor[rgb]{0.67,0.13,1.00}{##1}}}
\expandafter\def\csname PY@tok@ne\endcsname{\let\PY@bf=\textbf\def\PY@tc##1{\textcolor[rgb]{0.82,0.25,0.23}{##1}}}
\expandafter\def\csname PY@tok@nf\endcsname{\def\PY@tc##1{\textcolor[rgb]{0.00,0.00,1.00}{##1}}}
\expandafter\def\csname PY@tok@si\endcsname{\let\PY@bf=\textbf\def\PY@tc##1{\textcolor[rgb]{0.73,0.40,0.53}{##1}}}
\expandafter\def\csname PY@tok@s2\endcsname{\def\PY@tc##1{\textcolor[rgb]{0.73,0.13,0.13}{##1}}}
\expandafter\def\csname PY@tok@nt\endcsname{\let\PY@bf=\textbf\def\PY@tc##1{\textcolor[rgb]{0.00,0.50,0.00}{##1}}}
\expandafter\def\csname PY@tok@nv\endcsname{\def\PY@tc##1{\textcolor[rgb]{0.10,0.09,0.49}{##1}}}
\expandafter\def\csname PY@tok@s1\endcsname{\def\PY@tc##1{\textcolor[rgb]{0.73,0.13,0.13}{##1}}}
\expandafter\def\csname PY@tok@dl\endcsname{\def\PY@tc##1{\textcolor[rgb]{0.73,0.13,0.13}{##1}}}
\expandafter\def\csname PY@tok@ch\endcsname{\let\PY@it=\textit\def\PY@tc##1{\textcolor[rgb]{0.25,0.50,0.50}{##1}}}
\expandafter\def\csname PY@tok@m\endcsname{\def\PY@tc##1{\textcolor[rgb]{0.40,0.40,0.40}{##1}}}
\expandafter\def\csname PY@tok@gp\endcsname{\let\PY@bf=\textbf\def\PY@tc##1{\textcolor[rgb]{0.00,0.00,0.50}{##1}}}
\expandafter\def\csname PY@tok@sh\endcsname{\def\PY@tc##1{\textcolor[rgb]{0.73,0.13,0.13}{##1}}}
\expandafter\def\csname PY@tok@ow\endcsname{\let\PY@bf=\textbf\def\PY@tc##1{\textcolor[rgb]{0.67,0.13,1.00}{##1}}}
\expandafter\def\csname PY@tok@sx\endcsname{\def\PY@tc##1{\textcolor[rgb]{0.00,0.50,0.00}{##1}}}
\expandafter\def\csname PY@tok@bp\endcsname{\def\PY@tc##1{\textcolor[rgb]{0.00,0.50,0.00}{##1}}}
\expandafter\def\csname PY@tok@c1\endcsname{\let\PY@it=\textit\def\PY@tc##1{\textcolor[rgb]{0.25,0.50,0.50}{##1}}}
\expandafter\def\csname PY@tok@fm\endcsname{\def\PY@tc##1{\textcolor[rgb]{0.00,0.00,1.00}{##1}}}
\expandafter\def\csname PY@tok@o\endcsname{\def\PY@tc##1{\textcolor[rgb]{0.40,0.40,0.40}{##1}}}
\expandafter\def\csname PY@tok@kc\endcsname{\let\PY@bf=\textbf\def\PY@tc##1{\textcolor[rgb]{0.00,0.50,0.00}{##1}}}
\expandafter\def\csname PY@tok@c\endcsname{\let\PY@it=\textit\def\PY@tc##1{\textcolor[rgb]{0.25,0.50,0.50}{##1}}}
\expandafter\def\csname PY@tok@mf\endcsname{\def\PY@tc##1{\textcolor[rgb]{0.40,0.40,0.40}{##1}}}
\expandafter\def\csname PY@tok@err\endcsname{\def\PY@bc##1{\setlength{\fboxsep}{0pt}\fcolorbox[rgb]{1.00,0.00,0.00}{1,1,1}{\strut ##1}}}
\expandafter\def\csname PY@tok@mb\endcsname{\def\PY@tc##1{\textcolor[rgb]{0.40,0.40,0.40}{##1}}}
\expandafter\def\csname PY@tok@ss\endcsname{\def\PY@tc##1{\textcolor[rgb]{0.10,0.09,0.49}{##1}}}
\expandafter\def\csname PY@tok@sr\endcsname{\def\PY@tc##1{\textcolor[rgb]{0.73,0.40,0.53}{##1}}}
\expandafter\def\csname PY@tok@mo\endcsname{\def\PY@tc##1{\textcolor[rgb]{0.40,0.40,0.40}{##1}}}
\expandafter\def\csname PY@tok@kd\endcsname{\let\PY@bf=\textbf\def\PY@tc##1{\textcolor[rgb]{0.00,0.50,0.00}{##1}}}
\expandafter\def\csname PY@tok@mi\endcsname{\def\PY@tc##1{\textcolor[rgb]{0.40,0.40,0.40}{##1}}}
\expandafter\def\csname PY@tok@kn\endcsname{\let\PY@bf=\textbf\def\PY@tc##1{\textcolor[rgb]{0.00,0.50,0.00}{##1}}}
\expandafter\def\csname PY@tok@cpf\endcsname{\let\PY@it=\textit\def\PY@tc##1{\textcolor[rgb]{0.25,0.50,0.50}{##1}}}
\expandafter\def\csname PY@tok@kr\endcsname{\let\PY@bf=\textbf\def\PY@tc##1{\textcolor[rgb]{0.00,0.50,0.00}{##1}}}
\expandafter\def\csname PY@tok@s\endcsname{\def\PY@tc##1{\textcolor[rgb]{0.73,0.13,0.13}{##1}}}
\expandafter\def\csname PY@tok@kp\endcsname{\def\PY@tc##1{\textcolor[rgb]{0.00,0.50,0.00}{##1}}}
\expandafter\def\csname PY@tok@w\endcsname{\def\PY@tc##1{\textcolor[rgb]{0.73,0.73,0.73}{##1}}}
\expandafter\def\csname PY@tok@kt\endcsname{\def\PY@tc##1{\textcolor[rgb]{0.69,0.00,0.25}{##1}}}
\expandafter\def\csname PY@tok@sc\endcsname{\def\PY@tc##1{\textcolor[rgb]{0.73,0.13,0.13}{##1}}}
\expandafter\def\csname PY@tok@sb\endcsname{\def\PY@tc##1{\textcolor[rgb]{0.73,0.13,0.13}{##1}}}
\expandafter\def\csname PY@tok@sa\endcsname{\def\PY@tc##1{\textcolor[rgb]{0.73,0.13,0.13}{##1}}}
\expandafter\def\csname PY@tok@k\endcsname{\let\PY@bf=\textbf\def\PY@tc##1{\textcolor[rgb]{0.00,0.50,0.00}{##1}}}
\expandafter\def\csname PY@tok@se\endcsname{\let\PY@bf=\textbf\def\PY@tc##1{\textcolor[rgb]{0.73,0.40,0.13}{##1}}}
\expandafter\def\csname PY@tok@sd\endcsname{\let\PY@it=\textit\def\PY@tc##1{\textcolor[rgb]{0.73,0.13,0.13}{##1}}}

\def\PYZbs{\char`\\}
\def\PYZus{\char`\_}
\def\PYZob{\char`\{}
\def\PYZcb{\char`\}}

\def\PYZam{\char`\&}
\def\PYZlt{\char`\<}
\def\PYZgt{\char`\>}
\def\PYZsh{\char`\#}
\def\PYZpc{\char`\%}

\def\PYZhy{\char`\-}
\def\PYZsq{\char`\'}
\def\PYZdq{\char`\"}

\makeatother

\usepackage[hmarginratio=1:1,top=32mm,columnsep=20pt]{geometry} 
\usepackage{multicol} 
\usepackage[hang, small,labelfont=bf,up,textfont=it,up]{caption} 
\usepackage{booktabs} 
\usepackage{float} 
\usepackage[hidelinks]{hyperref} 
\usepackage{lettrine} 
\usepackage{paralist} 
\usepackage{abstract} 
\usepackage{outlines} 
\usepackage{mathtools}

\usepackage{titlesec} 
\renewcommand\thesection{\Roman{section}} 
\renewcommand\thesubsection{\Roman{subsection}} 
\titleformat{\section}[block]{\large\scshape\centering}{\thesection.}{1em}{} 
\titleformat{\subsection}[block]{\large}{\thesubsection.}{1em}{} 

\usepackage{fancyhdr} 
\title{\vspace{-15mm}{\large{
\selectfont\textbf{The Constituents of Sets, Numbers, and Other Mathematical Objects\\
Part Two}}}}

\small
\author{
\textsc{\small Ruadhan O'Flanagan}\footnote{\href{mailto:rof@ruadhan.net}{rof@ruadhan.net}}
\vspace{-3mm}
}
\vspace{3mm}
{\small\date{\small May 2019}}
\usepackage[format=plain,labelsep=newline,singlelinecheck=false]{caption}
\begin{document}

\maketitle 

\thispagestyle{fancy} 

\vspace{-5mm}

\begin{abstract}
\hrule
\vspace{0.5mm}
\hrule

\vspace{1mm}

\noindent

\small 

The arithmetic of natural numbers has a natural and simple encoding within sets, and the simplest set whose structure is not that of any natural number extends this set-theoretic representation to positive and negative integers. 

The operation that implements addition when applied to sets that represent natural numbers yields both addition and subtraction when used with the sets that encode integers.

The encoding of the integers naturally extends beyond them and identifies sets that encode arithmetic expressions and rational numbers. The sets that encode arithmetic expressions naturally specify the set operations that should be performed to evaluate those expressions.

The natural encoding of rational numbers within sets expresses each rational number as a novel order-preserving form of continued fraction, which provides a new efficient algorithm for finding rational approximations to irrational numbers. It also arranges all rational numbers within a tree that shows which numbers are constituents of which others. 

The part of this tree containing the positive rationals coincides with the well-known Stern-Brocot tree, which it extends to all rationals, including zero and negative numbers, introducing new non-trivial symmetries. These symmetries overlap with each other and form a group isomorphic to the group of symmetries of an equilateral triangle.

\vspace{2mm}

\hrule
\vspace{0.5mm}
\hrule
\end{abstract}

\section*{Introduction}

In Part One of this series\cite{partone}, we introduced the concept of a set's constituents and its constituent structure, and considered how natural numbers and ordered pairs and tuples can be encoded within sets in a natural way that results in a clear and simple alignment of set operations with the operations natural to the objects that the sets encode.

In this part, we apply the concept of constituent structure, and the necessity of encoding mathematical objects within it, beyond the natural numbers, to integers, arithmetic expressions, and rational numbers.

To avoid confusion regarding what exactly is being done and why, it is worth explaining in advance how the current project differs from what has been done before, in previous instances in which sets were used to define products, sums and other constructions involving numbers and operations on them.

One primary difference is that in most of the previous cases, the purpose of defining numbers, and their operations, using sets, was to use set theory to provide an unambiguous construction of the numbers that mathematicians used in practice, so that no dispute or confusion can ever arise regarding what a number actually is, or whether the type of existence that a number has is sufficient to justify its use in a given circumstance without introducing possible contradictions.

In those cases, it didn't matter precisely what sets were chosen to represent specific numbers. As long as some collection of sets could be found that had all the needed properties, arithmetic, along with the rest of mathematics, was on solid ground.

When it came to defining the sum, product, or difference of two numbers, when each number had already been assigned a set to represent it, the task at hand was to specify a procedure that a mathematician could follow, starting from the two sets representing the two numbers, to identify a third set, and that third set would have to be the set that represents the number reached by applying the arithmetic operation to the two numbers.

When this has been achieved, the product, sum or difference of any two numbers can be referenced with confidence and without any ambiguity, because it refers to a specific set and can no longer be interpreted in a way that might cause doubt to arise. ${}^1\!\!/_{\!\!7}$ has a decimal expansion that never terminates, and so does ${}^1\!\!/_{\!\!13}$. Adding these numbers by using the digits in their decimal expansions is a process that never finishes, so doubt can arise if we don't have a clear definition of the sum of two numbers, and one of the great achievements of twentieth century mathematics was to banish such doubt and confusion forever by grounding all of mathematics on the solid foundation of set theory.

With the current project, we start with the observation that we have a new way to understand the structure within a set, which allows us to draw a diagram that reveals what type of mathematical objects the set can represent. It discards the details of the set that belong only to sets - elements, subsets and cardinality, and displays the structure that remains when we remember what objects went into its construction, and what objects were used to construct those objects, but forget which sets we chose to represent those objects.

We also have a principle that tells us whether some way of encoding other mathematical objects within sets will lead to simple formulae and conceptual clarity later on, or will lead to complexity and confusion. For a number or other mathematical object that we wish to construct or represent using sets, we can choose the simplest set whose constituent structure matches that of the object, and the resulting set will have properties that correspond exactly to those of the object.

Along with these new instruments, we have a new algebra of sets, which was made visible by the concept of a constituent, and which operates specifically on constituent structure. We can construct a new set by replacing specific constituents within a set with different sets. Successively replacing the empty set, which every set contains, with other sets, allows us to construct or specify a complicated structure using a string of symbols such as $xyz$, denoting simpler sets used in the construction. The operation is associative and invertible, and there is a set, the empty set, that has no effect when added to the left or the right of such a string of sets.

A promising sign of the usefulness of this new algebra is given by the simplest sets with the constituent structure of the natural numbers, for which this basic set operation coincides with addition: $n(m)=m(n)=n+m$. We regard this as an indication that those sets are not just good enough to serve as representatives of the natural numbers, but that the natural numbers are {\em naturally} encoded by those sets. The arithmetic algebra of the natural numbers is the same as the constituent algebra of those sets, not by happy coincidence, but because each of those sets has exactly the same structure as the corresponding natural number. This brings the set algebra and the arithmetic algebra into alignment.

As we saw in the case of ordered pairs and tuples, other structures, more complicated than simple counts, can be constructed using the constituent substitution operation, and result in simple conceptually clear formulae.

So we now direct our attention to numbers with more complicated structure than a simple count, namely the integers, and, later on, the rational numbers. Unlike previous projects in which sets were found to represent the integers, we are not trying to find any satisfactory collection of sets that can be put into perfect correspondence with integers and in terms of which the arithmetic of integers can be defined.

Instead, we are trying to find the sets that have exactly the same structure as the integers. We conjecture that, like the natural numbers, each integer has a constituent structure that specifies all of its properties and nothing else, and that when the simplest set with the same constituent structure is identified for each integer, the constituent algebra of those sets will coincide with the arithmetic of integers.

It is certainly true that addition of natural numbers expressed in unary is the only arithmetic operation for which the result of performing the operation is always identical with the expression for it: Adding 1 to 11 gives 111 when 1 is the only digit, and concatenation is then identical with addition. This is why $n(m)=n+m$. In the case of operations such as multiplication, division, and subtraction, there will necessarily be a distinction between an expression for an operation and the result of the operation.

So unlike previous cases where products of numbers were defined using sets, in which there was a set for each of the numbers to be multiplied and a third set for the number that resulted from the multiplication, we will have four sets, with the additional set representing the arithmetic expression that designates the multiplication, and a procedure for evaluating the expression to reach the result.

We also depart from previous practice by imposing on ourselves the requirement that no choice of ours may be introduced into the doctrine. Our intention is to find the structure of arithmetic within the structure of sets, not to put it there. We are not building a copy of the numbers, or setting a new convention, whereby this or that set is selected to represent something because it serves the purpose as well as any other. We intend to identify, and make visible for the first time, in constituent structure diagrams, the structure of numbers themselves, by studying the simplest sets and their constituent algebra.

In the next section, where we consider the integers, we will find that concepts which we might consider artificial and contrived, such as positive and negative numbers, along with subtraction and addition of those numbers, are already encoded perfectly within the simplest possible constituent structures that aren't structures of natural numbers.

We will later see that the sets that naturally represent the rational numbers, identified here for the first time, specify an internal structure for each rational number consisting of a finite sequence of integers with remarkable properties. 

The sequence of integers is similar to the representation of the rational number as a continued fraction, but the sequence of integers is significantly more compact than that given by the standard continued fraction representation, and the algorithm that converts between the numerator/denominator form of the rational number and the sequence of integers is significantly more efficient, particularly in computationally complex cases. It provides a new and faster way to find rational approximations to irrational numbers, by achieving the same precision as the existing continued fraction algorithm in fewer steps.

This new representation of rational numbers is related to a known way of arranging the positive rational numbers in a binary tree known as the Stern-Brocot tree, and it provides new information about the structure of the tree, including an annotation of the edges, specifying the sets that each edge, when traversed, adds to the expression that specifies the set representing the rational number reached.

It also uniquely extends the tree to include all rational numbers including zero and negative numbers, preserving existing symmetries and introducing new symmetries which overlap with one another and have the same algebra as the group of symmetries of an equilateral triangle.

The arrangement of the negative rationals in the tree is quite non-trivial, but there is no choice of convention or new definition required to introduce them or determine their locations. The same rules that specify where each positive rational number appears in the tree, based on the sequence of integers in its natural representation, apply to negative numbers automatically.

The rational numbers and their natural representations as sequences of integers share the same order. At any given level in the tree, the rational numbers appear in increasing order from left to right and the sequences of integers that specify the internal structure of each rational number also appear in increasing order from left to right, ordered lexicographically. The novel form of continued fraction is an order-preserving isomorphism between those sequences and the rational numbers.

No considerations of algorithmic efficiency or order isomorphisms are involved in the process of identifying the sets that naturally encode the rational numbers. It is the study of the structure of the simplest sets, specifically the sets $\{0, 1\}$ and $\{\{1\}, \{0, 1\}\}$, which leads to the resulting representation of the rationals as finite sequences of integers.

Sets are nowhere to be seen in the novel type of continued fraction or the extended Stern-Brocot tree without annotations, but the consideration of sets and constituent structure is what leads to the improvements in efficiency and in the understanding of the relationship between rational numbers and integers. This suggests that the natural encoding of integers, arithmetic, and rational numbers within finite sets has been correctly identified.

That encoding is shown in its entirety in table \ref{tablearithmetic}, which specifies the sets that encode numbers and arithmetic expressions as well as the set operations that evaluate those expressions to yield the resulting numbers. The meaning of the equations defining the encoding of the natural numbers should already be clear from what has been said in Part One. The encoding of the integers and the rational numbers will be derived and explained in this part.

\begin{table}
\small
\caption{The Encoding of Arithmetic Within Finite Sets\label{tablearithmetic}}
\bgroup
\def\arraystretch{1.25}
\begin{tabular}{|c|c|c|c|}
\hline
            & Natural Numbers & Integers & Rational Numbers \\
\hline
Numbers     & $0=\{\}$        & $\Diamond=\{2, 2_{{}_{V}}\}$                      & $2_{{}_{V}}=\{0, 1\}$  \\
            & $n=\{n-1\}$     &  $+n=\Diamond n$                         & $\frac{1}{2+n}=n2_{{}_{V}}$  \\
            &                 &  $-n=n\Diamond$                          &  \\
\hline
Expressions & $``n+\cdots"=n $ &  $``x-\cdots" = x$    & $``\frac{1}{2+n+\cdots}"=n2_{{}_{V}}$\\
        & $``n+m+\cdots"=mn $  &  $``x-y-\cdots"=\Diamond yx$ & $``\frac{1+\cdots}{2+n}"=(0, n2_{{}_{V}})$ \\
            &                 &  $``x+y-\cdots"=y\Diamond x$ & $``\frac{1}{2+n}-\cdots"=(n2_{{}_{V}}1, n2_{{}_{V}})$ \\
            &                 &  $``(1+\cdots)\times y"=(0, y) $  & \\
            &                 &  $``(e)-\cdots"=(e1, e) $   & \\
\hline
Evaluation &                 &  $1\Diamond1\rightarrow \Diamond$  & $\Diamond 2_{{}_{V}} \rightarrow 2_{{}_{V}}\Diamond1$ \\
 & & $\Diamond\Diamond \rightarrow \{\}$ & $m2(0, m2_{{}_{V}}) \rightarrow (0, m2_{{}_{V}})1$ \\
            &                 &  $1(0, e) \rightarrow (0, e)\Diamond(e1, e)$ & $n2(m2_{{}_{V}}\Diamond, 0) \rightarrow m\Diamond n(n2_{{}_{V}}\Diamond, 0)1\Diamond2_{{}_{V}}$\\
            &                 &  $( b, (c, d))\rightarrow (b, d)$ & \\
            &                 &  $( (a, b), d)\rightarrow (b, d)$ & \\
\hline
Notation &   $ab\equiv a(\{\}\rightarrow b) $   & $(a,b)\equiv \{\{a\}, \{a, b\}\}$     & \\
\hline
\end{tabular}
\egroup
\end{table}

\newpage
\tableofcontents
\newpage

\section{Integers and Arithmetic Expressions}

The standard construction of the integers uses an ordered pair of natural numbers. In Part One, the simplest encoding of a finite sequence of natural numbers as a set emerged naturally as the set that specifies the coordinates of a set positioned deep within tuples of tuples. The diamond set, $\Diamond$,
separated each coordinate from the next, acting as a delimiter.

\begin{wrapfigure}{R}{3.6cm}
\centering
\scriptsize
\begin{tikzcd}
{} & (2, 3) \ar[-]{d} & {} \\
{} & (1, 3) \ar[-]{d} & {} \\
{} & (0, 3) \ar[-]{ld} \ar[-]{rd} & {}\\
\bullet \ar[-]{rd} & {} & \bullet \ar[-]{ld} {} \\
{} & \bullet \ar[-]{d} & {}\\
{} & 3 \ar[-]{d} & {} \\
{} & 2 \ar[-]{d} & {} \\
{} & 1 \ar[-]{d} & {} \\
{} & 0 & {} \\
\end{tikzcd}
\caption{{\scriptsize{When we use the set $n\Diamond m$ to represent the ordered pair of natural numbers, $(n, m)$, the set representing the pair $(2, 3)$ has the structure shown above.}}}
\end{wrapfigure}

That experience suggests that an ordered pair of natural numbers, $(n, m)$, can be represented by the set $n\Diamond m$, where $n$ and $m$ refer to Zermelo's constructions of those natural numbers.

When the corresponding set is constructed for the ordered pair $(2, 3)$, its structure diagram consists of 4 vertices arranged vertically at the bottom, corresponding to the sets representing Zermelo's natural numbers 0, 1, 2, and 3, followed by the diagram for $\Diamond$, followed by sets that have the proposed structure of ordered pairs, reaching $(2, 3)$ at the top.

When constructing the integers using ordered pairs, we attend to the difference between the two natural numbers in the pair. In the standard construction, an integer is identified with the infinite set of ordered pairs of natural numbers with a given difference, and whether the first or the second natural number in the pair is the negative, or subtracted, number, is a matter of convention.

In our case, it's not a matter of convention. We have a correspondence between set operations and natural number addition, given by $n(m) = n + m$, and we need to understand how the set $\Diamond$ interacts with the set operations and the arithmetic operations, in order to say whether $n$ or $m$ plays the role of a negative number within the structure $n\Diamond m$.

So rather than considering $\Diamond$ to be useful merely as something that is not a number, it is worth considering what its structure suggests about numbers before and after it.

The set $\Diamond$ is $\{\{1\}, \{1, 0\}\}$, which is the Kuratowski ordered pair $(1, 0)$, encapsulating the concept of 1 followed by 0. If the empty sets within this are replaced by a natural number, $n$, resulting in the set $\Diamond n$, then the result will be the Kuratowski ordered pair $(n+1, n)$. 

This ordered pair has a constituent structure that, among ordered pairs of arbitrary sets, doesn't unambiguously specify $n$ as the second entry in the pair. When representing integers, we need to ensure that no two integers are encoded in the same constituent structure, requiring us to look to other information to distinguish them. Among sets constructed using only $\Diamond$ and the sets representing natural numbers, however, $(n+1, n)$ is the only set with that structure, making it unambiguous and an acceptable representation.

The question then becomes how to interpret $1\Diamond n$, if $\Diamond n$ is $(n+1, n)$.

\begin{wrapfigure}[24]{R}{3.8cm}
\centering
\scriptsize
\begin{tikzcd}
{} & (4, 3, 2, 1) \ar[-]{d} & {} \\
{} & (4, 3, 2) \ar[-]{d} & {} \\
{} & (4, 3) \ar[-]{ld} \ar[-]{rd} & {}\\
\{4\} \ar[-]{rd} & {} & \{4, 3\} \ar[-]{ld} {} \\
{} & 4 \ar[-]{d} & {}\\
{} & 3 \ar[-]{d} & {} \\
{} & 2 \ar[-]{d} & {} \\
{} & 1 \ar[-]{d} & {} \\
{} & 0 & {} 
\end{tikzcd}
\caption{\label{countdown}\scriptsize{The set $\Diamond 3$ is the Kuratowski pair $(4,3)_K$. If the sets that follow it, $1\Diamond3=\{\Diamond3\}$ and so on, are interpreted as continuations of a sequence started by the pair $(4, 3)$, that sequence counts down.}}
\end{wrapfigure}

The set $\ 1(x)=\{\varnothing\}(x)=\{x\}$, which adds a single vertex to the top of $x$'s constituent structure, clearly corresponds to proceeding to the next number when $x$ is a natural number, but when $x$ is not a natural number but rather a pair of consecutive numbers in reverse order, such as $(n+1, n)$, a different sequence must be conceived of in order to say what comes next.

When we consider what mathematical object could be represented by the set $1\Diamond n$, by asking what naturally comes after $(n+1, n)$, the most obvious candidates are $(n+2, n+1)$, $(n, n-1)$, and $(n+1, n, n-1)$. The first two, however, can be ruled out, because $(n+2, n+1)=\Diamond (n+1)$ and $(n, n-1)=\Diamond(n-1)$, which are distinct sets within the same system of sets. 

That leaves $(n+1, n, n-1)$ as the natural interpretation of $1\Diamond n$. That is, $(n+1, n, n-1)$ is
a way of expressing the information given by $1\Diamond n$ -- ``that which comes next after $(n+1, n)$'' --
as a mathematical object of the same type, without losing or adding information.

This reveals that, when successive addition of a vertex is interpreted as proceeding to the next step of a sequence, $\Diamond n$ has a natural interpretation as the beginning of a sequence that counts downwards from the natural number $n$.

Sets of the form $m\Diamond$, where $m$ is a natural number, have natural interpretations as continuations of a downward count which goes past zero, reaching the negative integers.

\subsection{Expressions and Evaluation}

$m\Diamond n$ is obviously not the same set as the natural number $n-m$, presuming $n\geq m$, but its structure designates that number as the current step of a downward count. A downward count starting from that number is also designated by $\Diamond(n-m)$, which is a different set.

The set $\Diamond(n-m)$ has a structure consisting of that of a single natural number, $n-m$, topped by $\Diamond$, and can be considered to be the result of evaluating $m\Diamond n$, which has the structure of an arithmetic expression consisting of unary representations of the natural numbers $n$ and $m$ connected by a non-numeric structure designating subtraction. 

This is illustrated in figure \ref{countdown}, which shows how the structure of the set $2\Diamond 3$ represents an upwards count from 0 to 3, then a turnaround operation leaving the count at 3 and changing the direction of the count to downwards, and finally two steps further along the new sequence which
counts down from 3 to reach 1.

The fact that the structure of a set, $x$, appears to specify a procedure for reaching another set, $y$, implies that there is a natural encoding of certain operations on sets within the structure of sets like $x$.

To make the recognition of this encoding of set operations within set structure explicit, and to clarify the requirement that the encoding should be natural and consistent, we define:

\noindent {\bf Definition:} For a given class of finite pure sets, called expressions, which includes the empty set, an evaluation rule is a rule that assigns to every expression, $x$, an expression, $[x]$, called the evaluation of $x$, satisfying:
\begin{align}
[ab] = [ [a]b ] = [a[b]]
\end{align}
for any two expressions, $a$ and $b$, where $ab\equiv a(b)$.

This condition requires evaluation to respect associativity, $a(b(c))=a(b)(c)=abc$, and also ensures that $[[a]]=[a]$, since $b$ can be the empty set.

\subsection{The Operation Encoded in $\Diamond$}

We can use the requirement that evaluation must respect associativity to specify the evaluation of every expression, starting from a number of set replacement operations that define the encoding of an operation.

In the case of $\Diamond$, the defining replacement is\footnote{A more complete notation could express this as $[1\Diamond 1]_{1\Diamond 1 \rightarrow \Diamond} = \Diamond$. We will omit the subscript for brevity.}:
\begin{align}
[1\Diamond 1] &= \Diamond.
\end{align}

This allows us to conclude that:
\begin{align}
[11\Diamond11]=[1\Diamond1]=\Diamond
\end{align}
which shows that $[(n+m)\Diamond n]=m\Diamond$ and $[n\Diamond (n+m)]=\Diamond m$ for natural numbers, $n$ and $m$.

\begin{wrapfigure}[20]{R}{3.6cm}
\centering
\scriptsize
\begin{tikzcd}
{} & (1, 0, -1, -2, -3) \ar[-]{d} & {} \\
{} & (1, 0, -1, -2) \ar[-]{d} & {} \\
{} & (1, 0, -1) \ar[-]{d} & {} \\
{} & (1, 0) \ar[-]{ld} \ar[-]{rd} & {}\\
\{1\} \ar[-]{rd} & {} & \{1, 0\} \ar[-]{ld} {} \\
{} & 1 \ar[-]{d} & {} \\
{} & 0 & {}
\end{tikzcd}
\caption{{\scriptsize{The set $3\Diamond$, interpreted as the continuation of the sequence specified by $\Diamond$.}}}
\end{wrapfigure}
With this understanding of the relation between natural numbers, subtraction, and the diamond set, we can specify the construction of the positive and negative integers, $+n$ and $-n$, corresponding to the natural number $n$ as:
\begin{align}
+n  &\ \equiv \ \Diamond n \\
-n  &\ \equiv \ n \Diamond.
\end{align}

These sets include every result of evaluating an expression of the form $n\Diamond m$. They are the simplest sets that can be constructed from the natural number $n$ that don't have the structure of a natural number, and the set operation that coincides with addition of natural numbers produces both addition and subtraction when used with sets of this form.

Given an integer, $x$, and a natural number, $m$, we can specify the construction as a set of the unary expression for the addition, $x+m$, as:
\begin{equation}
``x+m" \equiv x(m)
\end{equation}
and the unary expression for the subtraction, $x-m$, as:
\begin{equation}
``x-m" \equiv m(x)
\end{equation}
where the quotation marks indicate that we are referring to the expressions, not the results of evaluating them, which are\footnote{Gratuitously abusing notation, we use $x$ and $m$ to refer to both the numbers and the sets representing those numbers. When $x$ appears beside a sign denoting an arithmetic operation such as $+$ or $-$, it denotes the number, and when it appears in a set operation such as $x(y)=xy$, it denotes the set. Arithmetic multiplication will always be explicitly denoted using $\times$ to prevent ambiguity.}:
\begin{align}
[x(m)] &= x+m \\
[m(x)] &= x-m
\end{align}
although the unary expression will be identical with the result of evaluating it if $m$ is added to a positive integer, $\Diamond n(m)$, or subtracted from a negative integer, $m(n)\Diamond$.

To specify constructions for expressions involving addition to integers and subtraction from integers of other integers, it is necessary to consider structures containing multiple diamond sets, such as $a\Diamond b\Diamond c$, where $a$, $b$, and $c$ are natural numbers.

When a sequence element, $s_n$, is used to replace the empty sets within $\Diamond$, the result is an ordered pair containing that sequence element occurring after its successor in that sequence: $\Diamond(s_n) = (s_{n+1}, s_n)$. So $\Diamond$ reverses the direction of the current sequence. Adding $\Diamond$ to a descending sequence will therefore result in an ascending sequence.

The structure $a\Diamond b\Diamond c$ will evaluate to an increasing sequence starting from $c-b+a$. An expression, $e$, that evaluates to an integer such as $m\Diamond$ or $\Diamond m$ must contain an odd number of diamonds in its structure, so that $1e$ indicates $e-1$ and $e1$ indicates $e+1$.

To get an integer from an expression with three $\Diamond$ structures, it is necessary to include an additional rule of replacement stating that reversing direction twice has no effect:
\begin{align}
[\Diamond\Diamond] &= \{\}.
\end{align}

The expressions for the sum and difference of two integers will therefore need to include a $\Diamond$ alongside the integers themselves. If $x$ and $y$ are integers of the form $\Diamond n$ and $\Diamond m$, which have the same sign, then $xy=x(y)=\Diamond n \Diamond m$ will perform a subtraction using $n$ and $m$. So the set that represents the expression for the subtraction, $x-y$, of one integer, $y$, from another, $x$, can be specified as:
\begin{equation}
``x-y" \equiv \Diamond y(x) = \Diamond yx
\end{equation}
and the expression for the addition, $x+y$, can be specified as:
\begin{equation}
``x+y" \equiv y(\Diamond x) = y\Diamond x.
\end{equation}

These expressions yield the corresponding results when evaluated:
\begin{align}
[\Diamond yx] &= x-y \\
[y\Diamond x] &= x+y
\end{align}
which shows that the integers, arithmetic expressions involving addition and subtraction, and the evaluation of those expressions to yield integers, all have natural representations as finite sets, given by the defining equations:
\begin{align}
+n  &\ \equiv \ \Diamond n \nonumber\\
-n  &\ \equiv \ n \Diamond \nonumber\\
``x+m" &\ \equiv x(m) \nonumber\\
``x-m" &\ \equiv m(x) \\
``x-y" &\ \equiv \Diamond y(x) \nonumber\\
``x+y" &\ \equiv y(\Diamond x)  \nonumber\\
[1\Diamond 1] &= \Diamond  \nonumber\\
[\Diamond\Diamond] &= \{\} \nonumber
\end{align}
where $n$ and $m$ are natural numbers, $x$ and $y$ are integers, quotation marks denote an expression, square brackets denote evaluation, and $\Diamond=(1, 0)=\{\{1\}, \{1, 0\}\}$.

\section{Multiplication}

\subsection{Arithmetic Expressions and Rules of Evaluation for Products of Integers}

The natural representation of the integers in the form of sets reveals that the Kuratowski ordered pair, $\Diamond=(1, 0)$,
naturally starts a sequence that counts down, with the succeeding terms in the sequence, $1\Diamond = \{\Diamond\}, 2\Diamond=\{\{\Diamond\}\}$, and so on, continuing the sequence, by successively making a step of the same magnitude in the same direction.

In the same way, an ordered pair of distinct integers naturally encodes a step from the first integer in the pair to the second, and a natural interpretation of the successor of that pair is another step with the same magnitude and direction, continuing the sequence started by the pair. 

The Kuratowski pair $(x, y)=\{\{x\}, \{x, y\}\}$ can encode all pairs of distinct positive integers unambiguously within its constituent structure, because $\Diamond n$ is not a constituent of $\Diamond(n+m)$ unless $m=0$. 

In the case when both integers are negative, $(x, y)=(n\Diamond, m\Diamond)$, or one integer is negative and the other is zero, the set representing one integer will be a constituent of that of the other, and the constituent structure of the Kuratowski pair will not unambiguously specify the order and structure of the elements it contains.

This ambiguity also occurs with Kuratowski ordered pairs of natural numbers, $(n, m)$, since that also causes the set representing one entry in the pair to be a constituent of the set representing the other.

A small modification to the encoding of a pair permits pairs containing negative integers, zero, and also natural numbers, to be unambiguously encoded within the constituent structure of sets:
\begin{equation}
(x, y) = (\Diamond x, \Diamond y)_K = \{\{\Diamond x\}, \{\Diamond x, \Diamond y\}\}
\end{equation} 

The only case when this pair fails is when the entries in the pair are equal, in which case $\{\{a\}, \{a, b\}\}$ reduces to $\{\{a\}\}=2a$, which is the set two steps after $a$ in whatever sequence $a$ itself naturally specifies. This is conceivably a more natural way to  interpret $(a, a)$ than as the specification of the start of a sequence that goes nowhere. 

We will take the hint and abstain from multiplying by zero. When we consider $a$ and $b$ to be distinct integers, the result, after evaluation, of proceeding one step further along the sequence started by the pair $(a, b)$ is the pair $(b, b+(b-a))$, which is the pair representing the step from $b$ that proceeds the same distance in the same direction as the previous step.

The evaluation rule for stepping forward one term in a sequence started by the pair of integers $(a, b)$, is then given by:
\begin{equation}
[1(a, b)] = [(a, b)ab]
\end{equation}
where we have left the brackets around the result on the right-hand side to indicate the fact that it is not the final result of the evaluation, since, for example, $[ab]$ involves a subtraction.

$\Diamond ab$ is the expression for the integer $b-a$. As an integer, this subtracts numbers supplied to it on the left, while $ab$ has the same numeric value as that integer, but adds numbers supplied to it on the left, since it contains two instances of $\Diamond$, one inside $a$ and one inside $b$.  The set $(a, b)ab$ is equal to $(aab, bab)$. The entries in this pair are expressions for $a-a+b$ and $b-a+b$, which provide the result $(b, b+(b-a))$.

For a natural number, $n$, the expression $n(a, b)$ will step the sequence forward $n$ times from the starting point, $b$, to $b+n(b-a)$. 

Multiplication of integers by natural numbers greater than zero can then be achieved by setting $a$ equal to the integer 0, which is represented by the set $\Diamond$:
\begin{equation}
[1(\Diamond, b)]=[(b, b\Diamond b)] = [(\Diamond, b)\Diamond b].
\end{equation}

The 1 on the left of the equation above steps the sequence forward from $(0, b)$ to $(b, b+b)=(b, 2\times b)$. The $n^{\rm th}$ step after $(0, b)$ will reach the numeric value $n\times b + b=(n+1)\times b$, encoded by the expression $n(\Diamond, b)$.

We can use trailing dots, $\cdots$, when we want to make it clear which part of an expression will be updated by taking that expression's successor:
\begin{equation}
``(1+\cdots)\times b" \equiv (\Diamond, b).
\end{equation}
This allows us to indicate that it does not subtract numbers supplied to it on the left like integers do:
\begin{equation}
``x-\cdots" \equiv x.
\end{equation}

To evaluate such an expression and get an integer as a result, it is necessary to start a new sequence counting downwards in steps of 1 from latest number of a given sequence.

The following rule of evaluation allows sequence step sizes to be changed in the necessary way:
\begin{equation}
[((a, b), (c, d))] = [((a, b), d)] = [(b, (c, d))] = [(b, d)].
\end{equation}

That is, only the last number reached by a sequence matters when switching to a new sequence; the steps leading to it, if there are any, are forgotten during evaluation when a new sequence is started.

This matches the structure of integers and their multiplicative sequences, since each positive integer is a pair of natural numbers, $\Diamond n=(1, 0)n=(n+1, n)$, and a pair of integers starts a sequence, $( (n+1, n), (m+1, m))\rightarrow (n, m)$ with the step size and position determined by the second entries in the original ordered pairs.

Given an expression, $e$, with an integer numeric value, but which is part of a sequence with a step size other than $-1$, an integer with that value can be obtained evaluating the expression $(e1, e)$. This starts a downward sequence at that numeric value with a step size of 1. We can explicitly denote this by showing the trailing dots subtracted from it: $``(n\times m)-\cdots"$.

In the case of multiplication, this converts the sequence generator $(n, m)$ into $(m+1, m)=(1, 0)m = \Diamond m$, allowing an integer result to be obtained.

This sequence-switching rule also automatically handles multiplication by negative integers:
\begin{equation}
[\Diamond(a, b)] = [(1, 0)(a, b)] = [(1(a, b), (a, b))] = [(b, bab), (a, b)] = [(bab, b)]
\end{equation}
where $(bab, b)$ refers to $(b+(b-a), b)$. This shows that $\Diamond$ naturally reverses the direction of an existing sequence. It changes the step size from $(b-a)$ to $-(b-a)$, while keeping the current position of the sequence at $b$.

Using this rule, the still-multiplying expression for the product of integers $x+1$ and $y$:
\begin{equation}
``(x+1-\cdots)\times y" \equiv x(\Diamond , y)
\end{equation}
can be turned into an expression for the integer that results from the multiplication:
\begin{equation}
``((x+1)\times y)-\cdots" \equiv \left(x(\Diamond, y)1,\ x(\Diamond, y)\right).
\end{equation}

\subsection{Multiplication With Ordered Pairs of Natural Numbers}

When the ordered pair contains natural numbers, $(n, m)$, the next term in the sequence is the pair $(m, m+(m-n))$, whose elements exceed those of the original pair by the amount $m-n$.

The evaluation rule for the corresponding sets is:
\begin{equation}
[1(n, m)] = [(n, m)\Diamond n\Diamond m].
\end{equation}

In the 
special case when $n=1$ and $m=0$, the pair is $(1, 0)=\Diamond$ and this gives the result:
\begin{equation}
[1\Diamond] = [1(1, 0)] = [(1, 0)\Diamond 1\Diamond0] = [\Diamond\Diamond 1 \Diamond] = 1\Diamond
\end{equation}
which shows that negative integers, $n\Diamond$, do not evaluate to anything simpler.

The evaluation of the expression for multiplication of natural numbers, $1(n, m) = (n, m)\Diamond n\Diamond m$, appears more complex than the result for integers, $1(a, b)=(a, b)ab$, but it coincides with the evaluation of the expression obtained by replacing $n$ and $m$ with the equivalent positive integers, $\Diamond n$ and $\Diamond m$:
\begin{equation}
[1(\Diamond n, \Diamond m)] = [(\Diamond n, \Diamond m)\Diamond n\Diamond m].
\end{equation}

This is because:
\begin{equation}
[(\Diamond n, \Diamond m)] = [( (n+1, n), (m+1, m) )] = [(n, m)]
\end{equation}
which shows that encoding an ordered pair, $(x, y)$ as $(\Diamond x, \Diamond y)_K=\{\{\Diamond x\}, \{\Diamond x, \Diamond y\}\}$ leads to exactly the same final evaluation of expressions as that obtained by using the pair $(x, y)_K=\{\{x\}, \{x, y\}\}$.

We can therefore use ordered pairs of natural numbers such as $(2, 0)$ without worrying about the encoding of their information within constituent structure, since the evaluation rule ensures that the result of evaluating expressions involving them will be the same as the result for the corresponding positive integers, whose order and constituent structures are unambiguously encoded in their ordered pairs.

One simple expression involving natural numbers occurs when $n=0$ and the pair is $(0, m)$:
\begin{equation}
[1(0, m)] \equiv (0, m)\Diamond \Diamond m = (0, m)m
\end{equation}
leading to a multiplication formula for natural numbers without $\Diamond$:
\begin{equation}
[n(0, m)]=(0, m)(n\times m).
\end{equation}

\subsection{The Minimal General Evaluation Rule for Multiplication}

The evaluation rule $[1(a, b)]=(a, b)ab$ is valid for integers because they subtract numbers supplied to them on the left. 

The set $(e1, e)$, has this property for any expression, $e$, while having the same numeric value as $e$, so the most general form of the multiplication rule, which is valid for any two expressions, $e_1$ and $e_2$, is:
\begin{equation}
[1(e_1, e_2)]=[(e_1, e_2)(e_11, e_1)(e_21, e_2)].
\end{equation}

It suffices to specify the evaluation rule for the case when one of the entries in the ordered pair is zero, since associativity will extend this rule to all pairs.

This allows multiplication of expressions in general to be evaluated with the rule:
\begin{equation}
[1(0, e)]=[(0, e)\Diamond(e1, e)].
\end{equation}

\section{Division and the Rational Numbers}

\subsection{Reciprocals}

Starting from the set, $2_{{}_{V}}=\{0,1\}$, and proceeding $n$ steps forward, to get the set $n2_{{}_{V}}=n(2_{{}_{V}})$, results in a set that, on its own, has a constituent structure indistinguishable from $2+n$, since the set $2_{{}_{V}}$ is von Neumann's construction of the number 2. The distinction becomes apparent in sets that contain successors of $2_{{}_{V}}$ and also natural numbers in Zermelo's representation.

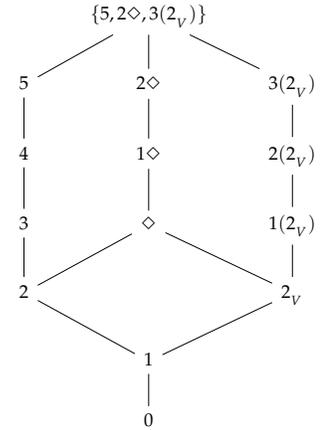
\begin{wrapfigure}{R}{3.6cm}
\centering
\scriptsize
\begin{tikzcd}
{} & \{5, 2\Diamond, 3(2_{{}_{V}})\} \ar[-]{ld}\ar[-]{d}\ar[-]{rd} & {}\\
5 \ar[-]{d}& 2\Diamond \ar[-]{d} & 3(2_{{}_{V}}) \ar[-]{d}\\
4 \ar[-]{d}& 1\Diamond \ar[-]{d} & 2(2_{{}_{V}}) \ar[-]{d}\\
3 \ar[-]{d}& \Diamond \ar[-]{ld} \ar[-]{rd} & 1(2_{{}_{V}})\ar[-]{d}\\
2 \ar[-]{rd} & {} & 2_{{}_{V}} \ar[-]{ld} {} \\
{} & 1 \ar[-]{d} & {} \\
{} & 0 & {} \\
\end{tikzcd}
\caption{\label{reciprocals}{\scriptsize{$2_{{}_{V}}=\{0, 1\}$ starts a sequence of sets that are distinguished from natural numbers and integers within structures that contain them.}}}
\end{wrapfigure}

When the natural numbers and $2_{{}_{V}}$ and its successors are seen within the same structure, it becomes clear that, alongside the sequence, $0, 1, 2, 3, \cdots$, there is another sequence, $0, 1, 2_{{}_{V}}, 1(2_{{}_{V}}), \cdots$, mirroring the sequence of natural numbers.

From the point of view of constituent structure alone, either sequence could be taken to be the natural numbers; all that the structure shows is that they are distinct from the number 2 onwards, with 2 and each succeeding natural number having a counterpart in the other sequence.

The number $1=\{\varnothing\}$ is the origin of both sequences. There are exactly two sets whose only constituents apart from themselves are $0=\{\}$ and $1=\{0\}$. One of these is $\{1\}$, the set that, in a sense, loses sight of 0 as it steps forward from 1, and which we, following Zermelo, identify with the number 2. 

The other set, $2_{{}_{V}}=\{0, 1\}$, ``straddles'' 1 and 0. In Conway's construction of surreal numbers\cite{conway}, the number $\{0|1\}$ is identified with ${}^1\!\!/_{\!\!2}$, and it appears to be natural to interpret $2_{{}_{V}}=\{0, 1\}$ as ${}^1\!\!/_{\!\!2}$ in our case, since the sequence of reciprocals, $1, {}^1\!\!/_{\!\!2}, {}^1\!\!/_{\!\!3}, \cdots$, does have the relationship to the sequence $1, 2, 3, \cdots$ that is visible in figure \ref{reciprocals}.

There are several other indications within the structure of sets that this is the right set to start the sequence of reciprocals. It is the simplest set available: $\{0, 1\}$ is the $4^{\rm th}$ simplest set, and is the simplest set that hasn't previously been assigned any numeric value.

The integer $0=\Diamond$, is the simplest set that cannot be a natural number because of its constituent structure, and it introduces subtraction, negative numbers, and the integers, extremely naturally, by introducing the sequence that starts with $(1, 0)$. This set, $\Diamond$, is also unique among all orderings of sets because it is the only structure of its kind that determines an unambiguous ordering of sets already ordered by constituency.

So $\Diamond$ is not an arbitrarily chosen set; it's a structure within sets that genuinely encodes properties of integers, and $\Diamond$ contains $2_{{}_{V}}=\{0, 1\}$ and $2=\{1\}$ as its only elements. So $2_{{}_{V}}$ is already within the system of numbers, but the number it represents has not previously been articulated. If $2_{{}_{V}}$ represents ${}^1\!\!/_{\!\!2}$ and its successors represent ${}^1\!\!/_{\!\!3}$, ${}^1\!\!/_{\!\!4}$ and so on, then every set within the system represents a number and has a successor that represents a number.

$2_{{}_{V}}$ also contains within its structure the assertion that only natural numbers that are greater than or equal to 2 have reciprocals, which is quite intelligible and natural. Finally, $2_{{}_{V}}$ is the only choice of set that produces the result that each number has a constituent structure isomorphic to that of its reciprocal. This is the property that specifies that a number and its reciprocal are a match for each other: It seems natural that only numbers with matching structures should be able to cancel each other arithmetically, and no other set can achieve this.

\begin{wrapfigure}{R}{3.6cm}
\centering
\scriptsize
\begin{tikzcd}
{} & \{5, -2, {}^1\!\!/_{\!\!5}\} \ar[-]{ld}\ar[-]{d}\ar[-]{rd} & {}\\
5 \ar[-]{d}& -2 \ar[-]{d} & {}^1\!\!/_{\!\!5}  \ar[-]{d}\\
4 \ar[-]{d}& -1 \ar[-]{d} & {}^1\!\!/_{\!\!4}  \ar[-]{d}\\
3 \ar[-]{d}& -0 \ar[-]{ld} \ar[-]{rd} & {}^1\!\!/_{\!\!3}  \ar[-]{d}\\
2 \ar[-]{rd} & {} & {}^1\!\!/_{\!\!2}  \ar[-]{ld} {} \\
{} & 1 \ar[-]{d} & {} \\
{} & 0 & {} \\
\end{tikzcd}
\caption{{\scriptsize{The constituent structure of a set containing natural numbers, negative numbers, and reciprocals.}}}
\end{wrapfigure}

With this interpretation, for any natural number $n$, we can specify the reciprocal number $1/(2+n)$ using the expression:
\begin{equation}
\frac{1}{2+n} \equiv n2_{{}_{V}}
\end{equation}
where $2_{{}_{V}}=\{0, 1\}$.

However, we need to remember that the next set in the sequence, $1n2_{{}_{V}}$, is the reciprocal of $n+3$. The expression $n2_{{}_{V}}$ is still in the process of accepting additions to the number whose reciprocal is to be constructed:
\begin{equation}
``\frac{1}{2+n+\cdots}" = n2_{{}_{V}}
\end{equation}
so the expression cannot be simply added and subtracted like an integer.

We can apply the operation $e \rightarrow (e1, e)$ to get a different expression that ends the process of constructing the denominator:
\begin{equation}
``\frac{1}{2+n}-\cdots" = (n2_{{}_{V}}1, n2_{{}_{V}})
\end{equation}
and this can then be added to and subtracted from integers using the existing rules:
\begin{align}
x+({}^1\!\!/_{\!\!2+n}) &= (n2_{{}_{V}}1, n2_{{}_{V}})\Diamond x  \\
x-({}^1\!\!/_{\!\!2+n}) &= \Diamond (n2_{{}_{V}}1, n2_{{}_{V}})x.
\end{align}

\subsection{Fractions}

The set $2_{{}_{V}}=\{0, 1\}$ and its successors, $n2_{{}_{V}}$, provide sets that represent fractions of the form $\frac{1}{2+n}$. The next task is to determine which sets provide the simplest and most natural representations of fractions with numerators other than 1.

A way to represent $\frac{m}{n}$ for any two natural numbers, $n$ and $m$, is actually automatically specified once reciprocals have been assigned a representation. Any rational number can be expressed in the form of a continued fraction:
\[
  a_0 + \cfrac{1}{a_1 + \cfrac{1}{a_2 + \cfrac{1}{ \ddots + \cfrac{1}{a_k} }}}
\]
where $a_0$ is an integer and $a_i$ is a positive number for $i>0$.

All rational numbers can be expressed as a continued fraction with a finite number of terms in the sequence, $[a_0; a_1, \cdots, a_k]$. The first term, $a_0$, specifies the integer part of a rational number, which can be positive, negative or zero, and the remaining terms, $a_1, \cdots, a_k$, specify a fraction between 0 and 1.

When $m_1<n_1$, the fraction $m_1/n_1$ has a reciprocal, $n_1/m_1$, which is greater than one, and therefore has an integer part, $a_1$. This provides the next term in the sequence $[a_0; a_1, \cdots, a_k]$. Subtracting $a_1$ from $n_1/m_1$ produces another fraction, $m_2/n_2$, which is less than one and has a reciprocal with another integer part. 

The numbers $n_i$ and $m_i$ get smaller each time because of the subtractions, and the sequence eventually terminates.

So we already have an expression for a fraction with a denominator other than 1, given by:
\begin{equation}
2_{{}_{V}}2_{{}_{V}} = \cfrac{1}{2 + \frac{1}{2}} = \cfrac{1}{\frac{5}{2}}=\frac{2}{5}
\end{equation}
where we have used the fact that:
$$
2_{{}_{V}} = \frac{1}{2+\cdots}
$$
and added another reciprocal rather than a natural number to the denominator.

Unfortunately, but interestingly, we cannot use the continued fraction representation to construct all rational numbers using $2_{{}_{V}}$ to represent reciprocals. The terms, $a_1, a_2, \cdots$, in the continued fraction representation of a rational number must be positive, but they can be equal to 1, while $2_{{}_{V}}$ can only construct reciprocals of natural numbers greater than or equal to 2.

We have no way to construct a fraction such as:
\[
  \cfrac{1}{1 + \cfrac{1}{1 + \frac{1}{2 }}}
\]
using only natural numbers and $2_{{}_{V}}$ in an expression of the form:
\begin{equation}
2_{{}_{V}}m2_{{}_{V}}n2_{{}_{V}} = \cfrac{1}{2 + n +\cfrac{1}{2+m+\frac{1}{2}}}.
\end{equation}

Our ability to represent reciprocals as sets starts with the number 2 for a good reason: 2 is the first natural number that has a reciprocal distinct from any natural number or integer. The reciprocal of 1, if it can be called that, is 1, and the next number after 1 is $\{1\}=2$, so the set that represents the number 1, namely $\{\varnothing\}$, is an expression for $1+\cdots$, not an expression for $\frac{1}{1+\cdots}$.

The reciprocal of $2$ is distinct from $2$, and is a completely different type of number, represented by a different type of set, $2_{{}_{V}} = \{0, 1\}$. This set can't be expressed as $x1=x(1)$ for any set, $x$, so it isn't an expression of the form $1+x$ for any number, $x$. It also can't be expressed as $1x$ for any set, $x$. This makes $\{0, 1\}$ suitable for starting the sequence of reciprocals and specifying the sets that encode them. 

So the structures of reciprocals and sets agree with each other that reciprocals start at 2, not 1. This raises the intriguing possibility that the natural representation of rational numbers within sets, when we find it, will tell us something we didn't already know about continued fractions.

\subsection{Representing Fractions with Positive Integers}

Since some fractions between 0 and 1 can't be constructed with any combination of natural numbers and reciprocals generated using $2_{{}_{V}}$, it will be necessary to use the $\Diamond$ set and introduce subtraction, which doesn't appear in standard continued fractions.

We can start by considering just positive, or, rather, non-negative integers, $\Diamond m$, which are expressions of the form $\Diamond m = m-\cdots$,
and in fact doing this immediately produces a way to construct every fraction between 0 and 1:
\begin{equation}
\Diamond m_3 2_{{}_{V}}\Diamond m_2 2_{{}_{V}}\Diamond m_1 2_{{}_{V}}
= \cfrac{1}{2 + m_1 -\cfrac{1}{2+m_2-\frac{1}{2+m_3-\cdots}}}.
\end{equation}

Since each reciprocal is between 0 and 1, subtracting it from the number $2+m$ results in a rational number between $1+m$ and $2+m$, where $m\geq 0$.
So each denominator in this variant of a continued fraction has an integer part greater than or equal to $1$, and a fractional part between 0 and 1. These are exactly the denominators that can be specified in standard continued fractions.

Positive integers, $\Diamond m$, and reciprocals generated by $2_{{}_{V}}$, are therefore able to represent all fractions between 0 and 1 in the form of a continued fraction that subtracts reciprocals instead of adding them.

We can specify the procedure for generating the sequence of numbers, $m_1, m_2, \cdots, m_k$, that are added to 2 in each denominator in the continued fraction expression for a fraction $\frac{n}{d}$ whose numerator, $n$, and denominator, $d$, are positive natural numbers with $n<d$:
\begin{enumerate}
\item 
Let $i=1$ and $n_1=n$ and $d_1=d$.
\item
If $n_i=1$, compute $m_i=d_i-2$ and stop.
\item
Compute $m_{i}=\rm{floor}(d_i/n_i)-1$.
\item
Compute $d_{i+1}=n_i$ and $n_{i+1}=(2+m_i)\times n_i - d_i$.
\item 
Increment $i$ and go to step 2.
\end{enumerate}

This procedure repeatedly finds the smallest multiple of the numerator that exceeds the denominator, and uses the amount by which it exceeds the denominator as the next numerator, terminating when the numerator reaches 1, indicating that a reciprocal of a natural number has been reached, providing the final reciprocal in the continued fraction.

The procedure for reconstructing the numerator and denominator of a rational number from the sequence, $m_1, \cdots, m_k$, is simply to define a function, $f(m_1, m_2, \cdots)$, as:
\begin{equation}
f(m_1, m_2, \cdots, m_k) = \frac{1}{2+m_1-f(m_2, \cdots, m_k)}
\end{equation}
with $f$ having a value of zero when it has zero arguments.

The sequence of numbers that this type of continued fraction assigns to any specific fraction will differ from the sequence that appears in its standard continued fraction, because in this case the numbers $n_1, n_2, \cdots$, can be zero.

\begin{table}
\centering
\caption{\label{positiveintegercontinuedfractions}Representations of Rational Numbers as Continued Fractions}
\begin{tabular}{|r|l|l|}
\hline
{} & {\scriptsize Standard Continued Fraction} & {\scriptsize Non-Negative Integers and $2_{{}_{V}}$ }\\
\hline
1/2 & 2 & 0 \\
1/3 & 3 & 1 \\
1/4 & 4 & 2 \\
1/5 & 5 & 3 \\
2/3 & 1, 2 & 0, 0\\
2/5 & 2, 2 & 1, 0\\
3/4 & 1, 3 & 0, 0, 0\\
3/5 & 1, 1, 2 & 0, 1\\
4/5 & 1, 4 & 0, 0, 0, 0 \\
\hline
\end{tabular}  
\end{table}

Table \ref{positiveintegercontinuedfractions} shows the sequences that this continued fraction expression and the standard one generate for fractions with numerators and denominators below 6.

Although this novel kind of continued fraction generates sequences with smaller numbers than the standard version, it is evidently very inefficient at encoding fractions of the form $m/(m+1)$, for which it uses a sequence of $m$ zeros. This will become awkward for large values of $m$. The standard version is able to encode those fractions using just two numbers rather than a long sequence.

We have yet to consider how the set $2_{{}_{V}}$ interacts with negative integers, of the form $n\Diamond$, which might provide a more efficient way to encode fractions that are awkward to represent using only non-negative integers.

\subsection{Reciprocals and Negative Integers}

When we consider the expression:
\begin{equation}
n\Diamond 2_{{}_{V}}
\end{equation}
where $n\Diamond$ is a negative integer, it seems natural to insert $n\Diamond$ into the denominator of the reciprocal in the place where $n$ appears in its corresponding expression:
\begin{equation}
n 2_{{}_{V}} = \frac{1}{2+n+\cdots}\ \ \longrightarrow\ \  n\Diamond 2_{{}_{V}} = \frac{1}{2-n-\cdots}.
\end{equation}

There are a few indications that this might not be right, though. Since 2 is positive, $2-n$ is not just a negative version of $2+n$. The negative integer whose reciprocal generated by $2_{{}_{V}}$ is equal in magnitude and opposite in sign that of $n$ would be $-(n+4)$, which is somewhat bizarre and breaks symmetries for no apparent reason.

It would also be the case that $4\Diamond 2_{{}_{V}}=2(2\Diamond 2_{{}_{V}})$ would be the number that is obtained by continuing two steps past the reciprocal of zero, which requires us to assume that the occurrence within a numerical sequence of a term with no definable numerical value does not disrupt the previous pattern.

The appearance and irremovability of the set $2\Diamond 2_{{}_{V}}$ representing the non-numeric expression $1/0$ would require special rules of evaluation for expressions containing it, and there would be no clear way for us to determine what the result of adding such a number to another, or multiplying by it, would be. 

For these reasons, we must consider the possibility that $n\Diamond 2_{{}_{V}}$ is not the reciprocal of $2-n-\cdots$, and that some considerations will lead us to another value for it.

We can list a few desiderata:

\begin{itemize}
\item
Associativity requires that $\Diamond 2_{{}_{V}}$ can be evaluated, and that it yields an expression, $[\Diamond 2_{{}_{V}}]$, with the same numeric value as $2_{{}_{V}}$, namely $\nicefrac{1}{2}$.
\item
Associativity also requires $[[\Diamond\Diamond] 2_{{}_{V}}]=2_{{}_{V}}$, so the expression for $[\Diamond 2_{{}_{V}}]$ would need to satisfy $[\Diamond[\Diamond 2_{{}_{V}}]]=2_{{}_{V}}$.
\item
$1[\Diamond 2_{{}_{V}}]$ should be larger than $\nicefrac{1}{2}$, and $n[\Diamond 2_{{}_{V}}]$ should increase as $n$ increases, because it goes in the opposite direction to $n 2_{{}_{V}}$, which decreases with $n$.
\item
The numeric value of $1[\Diamond 2_{{}_{V}}]$ should be determined by a rule that applies to sets in general and isn't just made up for this one case.
\end{itemize}

One candidate for a rule that could specify a value for $1[\Diamond 2_{{}_{V}}]$ is the rule that steps a sequence forward given an ordered pair of numeric values to start it.

The set $\Diamond 2_{{}_{V}}$ is an ordered pair, $(1, 0) 2_{{}_{V}}=(1(2_{{}_{V}}), 2_{{}_{V}})$, and the entries in the pair have numeric values, $(\frac{1}{2+1}, \frac{1}{2})=(\frac{1}{3}, \frac{1}{2})$.

This identifies the numeric value of the next term in the sequence started by that pair as:
\begin{equation}
1\Diamond 2_{{}_{V}} = \frac{1}{2} + (\frac{1}{2} - \frac{1}{3}) = \frac{1}{2} + \frac{1}{6} = \frac{2}{3}.
\end{equation}

This value can be obtained from $\frac{1}{2}$ by:
\begin{equation}
\frac{2}{3} = \frac{1+1}{2+1}
\end{equation}
which suggests that:
\begin{equation}
[\Diamond 2_{{}_{V}}] = \frac{1+\cdots}{2+\cdots} = 1 - \frac{1}{2+\cdots}=2_{{}_{V}}\Diamond 1.
\end{equation}

This evaluation rule satisfies all of the desiderata above: $1-\frac{1}{2}=\frac{1}{2}$ so the numerical value is correct, $1-(1-x)=x$ so $[\Diamond[\Diamond 2_{{}_{V}}]]=2_{{}_{V}}$, two thirds is larger than one half, and the rule that specified this value is simple sequence progression. It also creates a symmetry around the number $\frac{1}{2}$.

Equipped with this understanding of how negative integers interact with $2_{{}_{V}}$, we can consider the question of how they appear in the natural representations of rational numbers.

\subsection{An Efficient Continued Fraction Representation}

The evaluation rule for reciprocals involving negative integers:
\begin{equation}
[\Diamond 2_{{}_{V}}] = 2_{{}_{V}}\Diamond 1 \ \  \Longrightarrow\ \ [n\Diamond 2_{{}_{V}}] = n2_{{}_{V}}\Diamond 1
\end{equation}
immediately tells us how to handle negative integers that appear in a sequence of numbers representing a continued fraction constructed using $2_{{}_{V}}$.

Specifically, if the sequence of integers is $s_1, s_2, \cdots, s_k$, then the function:
\begin{equation}
\label{fractionfunction}
f(s_1, s_2, \cdots, s_k) = \begin{cases}
1-f(-s_1, -s_2, \cdots, -s_k) & s_1 < 0 \\
\cfrac{1}{2+s_1-f(s_2, \cdots, s_k)} & s_1 \geq 0 \\
0 & k=0 \\
\end{cases}
\end{equation}
computes the value of the fraction represented by that sequence.

The first case above implements the evaluation rule $\Diamond 2_{{}_{V}}=1-\frac{1}{2+\cdots}$. It reverses the sign of every subsequent integer in the sequence and not just the first one, because if the first term is positive, $s_1=\Diamond n_1$, then:
$$
s_22_{{}_{V}}s_12_{{}_{V}} = s_22_{{}_{V}}\Diamond n_12_{{}_{V}}
$$
while if $s_1$ is negative, $s_1=n_1\Diamond$:
$$
s_22_{{}_{V}}s_12_{{}_{V}} = s_22_{{}_{V}}n_1\Diamond 2_{{}_{V}} = s_22_{{}_{V}}n_12_{{}_{V}}\Diamond1
$$
which shows that both instances of $2_{{}_{V}}$ in the above expressions differ in sign. When the $\Diamond$ is moved to the right of the rightmost $2_{{}_{V}}$, the negative integer, $s_1$ between the first and the second $2_{{}_{V}}$ changes into a natural number, $n_1$, not into a positive integer, $\Diamond n_1$, and so it affects later reciprocals in the continued fraction.

The second case of equation \ref{fractionfunction}, when $s_1 \geq 0$, recursively evaluates the reciprocal to yield a fraction specified by a numerator and a denominator. The third case, $k=0$, specifies what happens when the number of arguments supplied to $f$ is zero. When that occurs, the recursion terminates: $f()=0$, so $f(x)=\frac{1}{2+x-0}$ when $x$ is positive.

The following procedure constructs the sequence, $s_1, s_2, \cdots, s_k$, starting from a fraction, $\frac{n}{d}$, with $0<n<d$:
\begin{enumerate}
\item 
Let $sign=1$, $i=1$, $n_1=n$ and $d_1=d$.
\item
If $n_i=1$, compute $s_i=sign\times(d_i-2)$ and stop.
\item
Compute $s_i=sign\times({\rm floor}(d_i/n_i)-1)$.
\item
If $s_i=0$, set $sign=-sign$ and $n_i=d_i-n_i$ and go to step 2.
\item
Compute $d_{i+1}=n_i$ and $n_{i+1}=(2+|s_i|)\times n_i - d_i$.
\item 
Increment $i$ and go to step 2.
\end{enumerate}

This procedure is similar to the one that was specified for the case when all the integers were non-negative, but step 4 checks to see if the numerator is larger than half of the denominator, and if it is, it replaces the numerator with the difference between it and the denominator. This difference will be less than half as large as the denominator and will therefore produce a reciprocal that exceeds two. Step 4 also records the occurrence of this replacement in the sign of the resulting sequence elements, so that it can be correctly decoded later.

This gives us a natural encoding of all the rational numbers between $0$ and $1$ within finite sets.

Any rational number that is not itself an integer is equal to an integer, $s_0$, minus a fraction between 0 and 1, so every rational number has a natural encoding in a set given by:
\begin{equation}
\label{rationalnumberencodedinset}
s_k2_{{}_{V}}\cdots s_22_{{}_{V}}s_12_{{}_{V}}s_0
\end{equation}
where $s_1, \cdots, s_k$ are the integers in the continued fraction sequence defined above, and $s_0$ is the integer from which that fraction is subtracted to obtain the full rational number. The fraction is subtracted rather than added because integers such as $s_0$ subtract numbers on their left and a more complicated, less symmetric expression involving a $\Diamond$ would be necessary to add the fraction rather than subtract it.

The sequences for standard continued fractions are usually specified using the notation $[a_0; a_1, a_2, \cdots, a_n]$, where the first entry in the sequence is separated from the rest by a semicolon to indicate its role as the integer part of the rational number. We can use the same notation for sequences specifying natural representations of rational numbers whose first term is an integer when the expression in equation \ref{rationalnumberencodedinset} would be awkward.

As the symmetry in this expression for a general rational number suggests, the integer part of a rational number can be included along with the fractional part in a procedure that handles them both in a consistent way. 

The encoding, $s(f)=[s_0; s_1, \cdots, s_k]$, of an arbitrary rational number, $f$, is given by:
\begin{equation}
s(f)=\begin{cases}
(f),& \ \ \ f-\lfloor f \rfloor=0\\
 (\lceil f \rceil) \ominus s(\frac{1}{f-\lfloor f \rfloor}-2), &\ \ \ f-\lfloor f \rfloor<\frac{1}{2} \\
 (\lceil f \rceil) \oplus s(\frac{1}{1-(f-\lfloor f \rfloor)}-2) , &\ \ \ f-\lfloor f \rfloor\geq\frac{1}{2} \\
\end{cases}
\end{equation}
with $\oplus$ defined as concatenation of tuples:
\begin{equation}
(a)\oplus(b, c, d, \cdots)\ \ \equiv \ \ (a, b, c, d, \cdots)
\end{equation}
and $\ominus$ denoting concatenation with the signs of the second sequence's terms flipped:
\begin{equation}
(a)\ominus(b, c, d, \cdots)\ \ \equiv \ \ (a, -b, -c, -d, \cdots).
\end{equation}

$\lfloor f \rfloor$ refers to the largest integer that doesn't exceed $f$ and $\lceil f \rceil$ is the smallest integer that is not less than $f$. This encoding chooses whether to take the reciprocal of the fractional part of $f$ or one minus that fractional part depending on which of those will yield a number greater than or equal to 2. Subtracting 2 from that number is then guaranteed to give a non-negative result.

Decoding is accomplished by the function:
\begin{equation}
f(s_0, s_1, s_2, \cdots, s_k)=\begin{cases}
s_0& \ \ \ k=0\\
s_0 - \frac{1}{2+f(s_1, s_2, \cdots, s_k)} &\ \ \ s_1 \geq 0 \\
s_0 - 1 + \frac{1}{2+f(-s_1, -s_2, \cdots, -s_k)} &\ \ \ s_1 < 0 \\
\end{cases}
\end{equation}
which uses the sign of each term in the sequence, $[s_0; s_1, \cdots, s_k]$, apart from the first term, to determine whether that term specifies the reciprocal of a fractional part of a number, $f-\lfloor f \rfloor$, or the reciprocal of $1-(f-\lfloor f \rfloor)$.

Two implementations of these functions are provided in appendix \ref{codeappendix}. The Python version is intended to be readable with no concern for efficiency, and works for positive and negative rational numbers, while the C version is intended to be fast and uses unsigned numerators and denominators.

The sequences that encode various rational numbers are shown in table \ref{finalcontinuedfractions}, with the sequences generated by the standard continued fraction shown alongside for comparison.

\begin{table}
\centering
\begin{tabular}{|r|l|l|}
\hline
{} & {\scriptsize Standard Continued Fraction} & {\scriptsize Natural Representation}\\
\hline
1/2 & [0; 2] & [1; 0] \\
1/3 & [0; 3] & [1; -1] \\
1/4 & [0; 4] & [1; -2] \\
1/5 & [0; 5] & [1; -3] \\
2/3 & [0; 1, 2] & [1; 1] \\
3/2 & [1; 2] & [2; 0] \\
2/5 & [0; 2, 2] & [1; -1, 0] \\
5/2 & [2; 2] & [3; 0] \\
3/4 & [0; 1, 3] & [1; 2] \\
4/3 & [1; 3] & [2; -1] \\
3/5 & [0; 1, 1, 2] & [1; 1, 0] \\
5/3 & [1; 1, 2] & [2; 1] \\
4/5 & [0; 1, 4] & [1; 3] \\
21/29 & [0; 1, 2, 1, 1, 1, 2] & [1; 2, 1, 1] \\
89/144 & [0; 1, 1, 1, 1, 1, 1, 1, 1, 1, 2] & [1; 1, 1, 1, 1, 1] \\
\hline
\end{tabular}
\caption{\label{finalcontinuedfractions}Rational numbers expressed as sequences of integers given by their continued fraction expansions and their natural representations.}
\end{table}

The table shows that the sequences generated by the natural encoding of the fractions tend to be shorter and contain numbers which are closer to zero. 

The final entry, $89/144$, in table \ref{finalcontinuedfractions} is a ratio of consecutive terms in the Fibonacci sequence, which provides an approximation to the fractional part of the golden ratio, $\varphi=(1+\sqrt{5})/2$. The golden ratio, $1.618\cdots$, is the irrational number that is most difficult to approximate using rational numbers. Both the natural representations and the standard continued fraction representations of the rational approximations to $\varphi$ contain long strings of consecutive 1's, which provide the smallest possible update to the number with each additional term in the sequence.

For each such approximation, the natural representation encodes the same number in a sequence which is approximately half as long as that of the continued fraction. This compact encoding is accompanied by an improvement in the computational efficiency.

Table \ref{benchmarking} shows the time taken by both the standard and the natural continued fraction algorithms to encode and decode ratios of successive Fibonacci numbers, using the implementations written in C given in the appendix, using \texttt{gcc}'s compiler optimizations and run on a machine with an AMD Ryzen 7 1700x processor.

\begin{table}
\centering
\begin{tabular}{|r|r|r|r|r|}
\hline

{} & \multicolumn{2}{c|}{\scriptsize Continued Fraction} & \multicolumn{2}{c|}{\scriptsize Natural Representation} \\\cline{2-5}
{\scriptsize Ratio} & {\scriptsize Encoding} & {\scriptsize Decoding} & {\scriptsize Encoding} & {\scriptsize Decoding} \\
\hline
$f_5/f_6$       &  156 &  25 &  114 &   22 \\
$f_{10}/f_{11}$ &  341 &  57 &  226 &   56 \\
$f_{20}/f_{21}$ &  791 & 124 &  414 &  111 \\
$f_{30}/f_{31}$ & 1682 & 281 &  605 &  179 \\
$f_{40}/f_{41}$ & 2140 & 353 &  783 &  241 \\
$f_{50}/f_{51}$ & 2667 & 454 &  970 &  305 \\
$f_{60}/f_{61}$ & 3201 & 560 & 1231 &  379 \\
$f_{70}/f_{71}$ & 3712 & 667 & 1352 &  503 \\
$f_{80}/f_{81}$ & 4268 & 772 & 1596 &  581 \\
\hline
\end{tabular}  
\caption{\label{benchmarking}Time Taken to Encode and Decode Ratios of Fibonacci Numbers, $f_1=f_2=1, f_{n+1}=f_n+f_{n-1}$. Each entry in the table shows the time taken in milliseconds for 10,000,000 successive repetitions of the encoding or decoding operation.}
\end{table}

As the table shows, the improvement in performance achieved by using the natural representation instead of the standard continued fraction is quite significant and becomes more pronounced as the computational complexity increases.
Because of this, the natural representation provides a faster algorithm for generating rational approximations to certain irrational numbers. 

For example, the standard continued fraction expansion of $\sqrt{3}$ is $[1;1,2,1,2,1,2,\cdots]$, while the natural version is $[2; 2, 2, 2, \cdots]$. The same approximation to $\sqrt{3}$ can be achieved by evaluating a natural representation's continued fraction with a finite sequence of 2's, or by evaluating a standard continued fraction using a sequence that is twice as long and takes about 50\% more time to compute.

So the project of identifying the natural encoding of rational numbers within sets unexpectedly leads to a more compact representation of those numbers and a faster algorithm for computing them.

\subsection{Evaluation Rules for Division}

A rational number is naturally encoded by a set of the form $s_k2_{{}_{V}}\cdots s_22_{{}_{V}}s_12_{{}_{V}}s_0$ so there must be rules of evaluation that convert a product of an integer and a reciprocal into this form.

For simplicity, we can express the rules for multiplication using ordered pairs in which the second entry is zero:
\begin{equation}
``\cdots\times x" = (x\Diamond, 0).
\end{equation}

This represents the start of a sequence that is initially at zero and has a step size of $x$. The requirement that evaluation must respect associativity ensures that other expressions involving multiplication that can be transformed into an expression containing this one will evaluate to the correct result.

If $x$ is an integer, then:
\begin{equation}
\label{multiplyby1}
[1(x\Diamond, 0)]=[1(x, 0)\Diamond]=[(x\Diamond, 0)\Diamond x]
\end{equation}
which says:
\begin{equation}
``(1+\cdots)\times x" = ``-((1+\cdots)\times -x)" = ``\cdots\times x+x".
\end{equation}

If we consider the product of a natural number greater than 1 with a reciprocal of a natural number:
\begin{equation}
``(\cdots+2+n)\times \frac{1}{2+m}" = n2(m2_{{}_{V}}\Diamond, 0)
\end{equation}
then the integer part of the fraction $(2+n)/(2+m)$ is generated by the evaluation rule:
\begin{equation}
[m2(m2_{{}_{V}}\Diamond, 0)]=[(m2_{{}_{V}}\Diamond, 0)1]
\end{equation}
which is equivalent to:
\begin{equation}
[m2(0, m2_{{}_{V}})]=[(0, m2_{{}_{V}})1].
\end{equation}

After that rule has been applied to any expression that multiplies a natural number by a reciprocal, the resulting expression will contain a natural number, $a$, plus a fractional part with a denominator equal to $2+m$. 

This is equivalent to $a+1$ minus a different fraction with the same denominator. This generates the integer part, $s_0=\Diamond 1a$, of the natural encoding of the rational number $s_k2_{{}_{V}}\cdots s_22_{{}_{V}}s_12_{{}_{V}}s_0$.

If the numerator of the subtracted fraction is zero, the rational number that results from the division has no fractional part and the process of division is finished. If it is 1, then the rational number is of the form $s_1\twov s_0$, where $s_1=\Diamond m$, and the division is finished.

So we can assume that the remaining numerator is at least 2, and we can also assume that it is less than half of the denominator, since the rule $\Diamond\twov=\twov\Diamond1$ can be used to replace $n$ with $m-n$ if it is greater than half of the denominator, with the consequence that the integer $s_1$ that results will be negative.

The resulting fraction can be expressed as:
\begin{equation}
\frac{n_1+2}{m_1+2} = \cfrac{1}{\frac{m_1+2}{n_1+2}} = \cfrac{1}{2+\frac{m_1+2}{n_1+2}-\frac{2\times n_1+4}{n_1+2}}= \cfrac{1}{2+\frac{m_1-2\times n_1-2}{n_1+2}}
\end{equation}
where $m_1-2\times n_1-2$ is non-negative. The set that encodes this expression is:
\begin{equation}
m_1\Diamond n_1n_12\Diamond(n_12_{{}_{V}}\Diamond, 0)2_{{}_{V}}
\end{equation}
and $m_1\Diamond n_1n_12\Diamond$ will evaluate to a non-negative natural number. We can shorten this expression by using the non-positive number $-m_1+2\times n_1+2=m_1\Diamond n_1n_12$:
\begin{equation}
m_1\Diamond n_1n_12(n_12_{{}_{V}}, 0)2_{{}_{V}}
\end{equation}
which is equivalent under the rules of evaluation to:
\begin{equation}
m_1\Diamond n_1(n_12_{{}_{V}}, 0)\Diamond1\Diamond2_{{}_{V}}=m_1\Diamond n_1(n_12_{{}_{V}}\Diamond, 0)1\Diamond2_{{}_{V}}.
\end{equation}

If the numerator in this resulting fraction is greater than $2+n_1$, then equation \ref{multiplyby1} will apply, and will eventually result in an expression of the form:
\begin{equation}
h(n_12_{{}_{V}}\Diamond, 0)s_12_{{}_{V}}
\end{equation}
where $h<2+n_1$ and $s_1$ is an integer.

If $h=0$ then the process has finished. If $h=1$, then it has reached the final reciprocal. If $h>1$, then the expression can be written as:
\begin{equation}
n_22(m_22_{{}_{V}}\Diamond, 0)s_12_{{}_{V}}
\end{equation}
where $m_2=n_1$ and $n_2$ is a natural number satisfying $n_2<m_2$. The above process for $n_1$ and $m_1$ can then be repeated with $n_2$ and $m_2$, after which the expression for the rational number will have $s_22_{{}_{V}}s_12_{{}_{V}}s_0$ on the right of the ordered pair.

The numerators and denominators get smaller in each iteration of this procedure until the process terminates, yielding an expression of the form $(m_k2_{{}_{V}}\Diamond, 0)s_k2_{{}_{V}}\cdots s_22_{{}_{V}}s_12_{{}_{V}}s_0$ where the final multiplicative operator $(m_k2_{{}_{V}}\Diamond, 0)$ has nothing left to multiply and the signs of the integers $s_0$, $s_1$ and so on are chosen to ensure that the numerator is less than half of the denominator before each reciprocal is taken.

The evaluation rule that generates the natural representation of the rational number that results from division is therefore:
\begin{equation}
[n2(m2_{{}_{V}}\Diamond, 0)] = [m\Diamond n(n2_{{}_{V}}\Diamond, 0)1\Diamond2_{{}_{V}}].
\end{equation}

\section{What Does the Natural Encoding of a Rational Number Mean?}

For natural numbers, the set that encodes a number has a very clear relation to that number, and the constituent structure of the set clearly mirrors the structure of the number. With the integers, once we become comfortable with the concept that the structure of the $\Diamond$ set naturally reverses the direction of a sequence, the constituent structure of the set that represents an integer can be understood as an accurate representation of the integer itself.

In this case, the set, $s_k2_{{}_{V}}\cdots2_{{}_{V}}s_1 2_{{}_{V}}s_0$, representing a rational number, $s$, doesn't match our concept of a fraction as a numerator and a denominator. Its constituent structure is just a restatement of its expression: several integers with instances of $2_{{}_{V}}$ between each integer and the next.

We would like to have a clear understanding of what the integers mean, not as an input to an algorithm that computes the rational number, but as a part of the rational number's structure.

The natural numbers and their reciprocals taught us that distinct numbers represented by distinct sets can have identical constituent structures, and we can see how they relate to one another in the constituent structure diagram of a set that contains both of those distinct sets.

So we will consider a set that contains multiple rational numbers and natural numbers as constituents. Figure \ref{rationalstructure} shows the constituent structure of the set
$\{{}^1\!\!/_{\!\!4}, {}^2\!\!/_{\!\!5}, {}^3\!\!/_{\!\!5}, {}^3\!\!/_{\!\!4}, {}^4\!\!/_{\!\!3}, {}^5\!\!/_{\!\!3}, {}^5\!\!/_{\!\!2}, 4\}$.

\begin{figure}
\begin{tikzpicture}[baseline= (a).base]
\node[scale=.6] (a) at (-100,0){
\begin{tikzcd}
{}&{}     & {} & {}            & {} & {}  & {}  & {}  & \ar[-, gray, dashed]{ddlllllll}\ar[-, gray, dashed]{ddlllll}\ar[-, gray, dashed]{ddlll}\ar[-, gray, dashed]{ddl}\bullet\ar[-, gray, dashed]{ddrrrrrrr}\ar[-, gray, dashed]{ddrrrrr}\ar[-, gray, dashed]{ddrrr}\ar[-, gray, dashed]{ddr} & {} & {}  & {}  & {}  & {}  & {}  & {}  & {} & {}  \\
{}           & {}         & {} & {}  & {}  & {}  & {}  & {}  & {}  & {}  & {}  & {}  & {}  & {} & {} & {} \\
{} &[-25pt] 1/4 &[-25pt] {} &[-25pt] 2/5 &[-25pt] {} &[-25pt] 3/5 &[-25pt] {} &[-25pt] 3/4 &[-25pt] {} &[-25pt] 4/3 &[-25pt] {} &[-25pt] 5/3 &[-25pt] {} &[-25pt] 5/2 &[-25pt] {} &[-25pt] 4 \\ 
{} &[-25pt] {} &[-25pt] {} &[-25pt] {} &[-25pt] {} &[-25pt] {} &[-25pt] {} &[-25pt] {} &[-25pt] {} &[-25pt] {} &[-25pt] {} &[-25pt] {} &[-25pt] {} &[-25pt] {} &[-25pt] {} &[-25pt] {} \\ 
{} &[-25pt] {} &[-25pt] {} &[-25pt] {} &[-25pt] {} &[-25pt] {} &[-25pt] {} &[-25pt] {} &[-25pt] {} &[-25pt] {} &[-25pt] {} &[-25pt] {} &[-25pt] {} &[-25pt] {} &[-25pt] {} &[-25pt] {} \\ 
{} &[-25pt] {} &[-25pt] 1/3\ar[-]{uuul}\ar[-, "\Diamond2_{{}_{V}}\Diamond" description]{uuur} &[-25pt] {} &[-25pt] {} &[-25pt] {} &[-25pt] 2/3\ar[-]{uuur}\ar[-, "\Diamond2_{{}_{V}}\Diamond" description]{uuul} &[-25pt] {} &[-25pt] {} &[-25pt] {} &[-25pt] 3/2\ar[-, "\Diamond" near start]{uuur}\ar[-]{uuul} &[-25pt] {} &[-25pt] {} &[-25pt] {} &[-25pt] 3\ar[-]{uuur}\ar[-, "\Diamond2_{{}_{V}}\Diamond" description]{uuul} &[-25pt] {} \\ 
{} &[-25pt] {} &[-25pt] {} &[-25pt] {} &[-25pt] {} &[-25pt] {} &[-25pt] {} &[-25pt] {} &[-25pt] {} &[-25pt] {} &[-25pt] {} &[-25pt] {} &[-25pt] {} &[-25pt] {} &[-25pt] {} &[-25pt] {} \\ 
{} &[-25pt] {} &[-25pt] {} &[-25pt] {} &[-25pt] 1/2\ar[-, "\Diamond" near start]{uurr}\ar[-]{uull} &[-25pt] {} &[-25pt] {} &[-25pt] {} &[-25pt] {} &[-25pt] {} &[-25pt] {} &[-25pt] {} &[-25pt] 2\ar[-]{uurr}\ar[-, "\Diamond2_{{}_{V}}\Diamond" description]{uull} &[-25pt] {} &[-25pt] {} &[-25pt] {} \\ 
{} &[-25pt] {} &[-25pt] {} &[-25pt] {} &[-25pt] {} &[-25pt] {} &[-25pt] {} &[-25pt] {} &[-25pt] 1\ar[-]{urrrr}\ar[-, "\Diamond2_{{}_{V}}\Diamond" description]{ullll} &[-25pt] {} &[-25pt] {} &[-25pt] {} &[-25pt] {} &[-25pt] {} &[-25pt] {} &[-25pt] {} \\ 
0\ar[-]{urrrrrrrr} &[-25pt] {} &[-25pt] {} &[-25pt] {} &[-25pt] {} &[-25pt] {} &[-25pt] {} &[-25pt] {} &[-25pt] {} &[-25pt] {} &[-25pt] {} &[-25pt] {} &[-25pt] {} &[-25pt] {} &[-25pt] {} &[-25pt] {} \\ 
\end{tikzcd}
};
\end{tikzpicture}
\caption{\label{rationalstructure}The constituent structure of a set containing several rational numbers.
The path connecting 0 to a number, $s$, specifies the set, $s_k 2_{{}_{V}}\cdots s_2 2_{{}_{V}} s_1 2_{{}_{V}} s_0$, that encodes that number, and the sequence of integers, $[s_0; s_1, \cdots, s_k]$. Each edge traversed along the path adds symbols to the left of the expression for the constructed set. Unbroken edges add 1; edges containing $\Diamond2_{{}_{V}}\Diamond$ add $\Diamond2_{{}_{V}}\Diamond$; edges marked with $\Diamond$ add the $\Diamond$ to the expression, possibly cancelling an existing $\Diamond$, before the 1 specified by the edge is added. The set for $5/3$, for example, is $1\Diamond\Diamond 2_{{}_{V}}\Diamond2=1 (2_{{}_{V}}\Diamond2)$, indicating the numeric value $2-\frac{1}{3}$ and the natural representation $[2; 1]$.
}
\end{figure}

Each node in the diagram corresponds to a rational number. The structures of the $\Diamond$ and $2_{{}_{V}}$ sets situated on top of numbers have been omitted for clarity. Symbols along the edges connecting one number to another indicate the sets that must be added to the set for one number to reach the set for the other, after evaluating the resulting expression.

The edge connecting 1 to ${}^1\!\!/_{\!\!2}$ contains $\Diamond2_{{}_{V}}\Diamond$, which means that evaluating the expression $\Diamond2_{{}_{V}}\Diamond1$ produces ${}^1\!\!/_{\!\!2}=2_{{}_{V}}$:
$$
[\Diamond2_{{}_{V}}\Diamond1] = [\Diamond[2_{{}_{V}}\Diamond1]] = [\Diamond\Diamond2_{{}_{V}}]=2_{{}_{V}}.
$$

The diagram provides insight into the meaning of the integers, $[s_0; s_1, \cdots, s_k]$, that appear in the continued fraction expression for a rational number encoded using $2_{{}_{V}}$. The integers specify a route within this binary tree structure that connects $0$ to the specific rational number encoded.

A positive integer, $\Diamond n$, or a negative one, $n \Diamond$, specifies the number, $n$, of consecutive steps along the same direction that should be followed before changing direction or reaching the end of the route.

Successive integers in the set representing a rational number, $s_k 2_{{}_{V}}\cdots s_2 2_{{}_{V}} s_1 2_{{}_{V}} s_0$, are separated by an occurrence of $2_{{}_{V}}$, which changes the direction of the route through the binary tree by traversing one of the edges labelled with $\Diamond2_{{}_{V}}\Diamond$ in figure \ref{rationalstructure}. 

This brings the route into a different sequence of numbers by introducing the reciprocal:
$$
m+\cdots\  \longrightarrow \  m - \frac{1}{2+\cdots}
$$

The sign of the integer that follows the $2_{{}_{V}}$ determines which of the two edges that proceed away from the current number will be followed in the next step. Although it may seem counterintuitive, a negative integer, $n\Diamond$, specifies that the edge without the $\Diamond$ symbol should be traversed. The reason for this is that integers subtract numbers that are supplied to them on the left, so in the expression for the set, $s_k 2_{{}_{V}}\cdots s_2 2_{{}_{V}} s_1 2_{{}_{V}} s_0$, each instance of $2_{{}_{V}}$ is subtracted from the integer on its right, leading to a downward sequence that a subsequent negative integer can reverse.

\begin{figure}
\begin{tikzpicture}[baseline= (a).base]
\node[scale=.6] (a) at (-100,0){
\begin{tikzcd}
{} &[-10pt] {}           &[-10pt] {} &[-10pt] {}            &[-10pt] {} &[-10pt] {}  &[-10pt] {}  &[-10pt] {}  &[-10pt] \ar[-, gray, dashed]{ddlllllll}\ar[-, gray, dashed]{ddlllll}\ar[-, gray, dashed]{ddlll}\ar[-, gray, dashed]{ddl}\bullet\ar[-, gray, dashed]{ddrrrrrrr}\ar[-, gray, dashed]{ddrrrrr}\ar[-, gray, dashed]{ddrrr}\ar[-, gray, dashed]{ddr} &[-10pt] {} &[-10pt] {}  &[-10pt] {}  &[-10pt] {}  &[-10pt] {}  &[-10pt] {}  &[-10pt] {}  &[-10pt] {} &[-10pt] {} \\
{} &[-10pt] {} &[-10pt] {} &[-10pt] {} &[-10pt] {} &[-10pt] {} &[-10pt] {} &[-10pt] {} &[-10pt] {} &[-10pt] {} &[-10pt] {} &[-10pt] {} &[-10pt] {} &[-10pt] {} &[-10pt] {} &[-10pt] {} \\ 
{} &[-10pt] \left [1;-2\right ] &[-10pt] {} &[-10pt] \left [1;-1,0\right ] &[-10pt] {} &[-10pt] \left [1;1,0\right ] &[-10pt] {} &[-10pt] \left [1;2\right ] &[-10pt] {} &[-10pt] \left [2;-1\right ] &[-10pt] {} &[-10pt] \left [2;1\right ] &[-10pt] {} &[-10pt] \left [3;0\right ] &[-10pt] {} &[-10pt] \left [4\right ] \\ 
{} &[-10pt] {} &[-10pt] {} &[-10pt] {} &[-10pt] {} &[-10pt] {} &[-10pt] {} &[-10pt] {} &[-10pt] {} &[-10pt] {} &[-10pt] {} &[-10pt] {} &[-10pt] {} &[-10pt] {} &[-10pt] {} &[-10pt] {} \\ 
{} & {} & {} & {} & {} & {} & {} & {} & {} & {} & {} & {} & {} & {} & {} & {} \\ 
{} &[-10pt] {} &[-10pt] \left [1;-1\right ]\ar[-]{uuul}\ar[-, "\Diamond2_{{}_{V}}\Diamond" description]{uuur} &[-10pt] {} &[-10pt] {} &[-10pt] {} &[-10pt] \left [1;1\right ]\ar[-]{uuur}\ar[-, "\Diamond2_{{}_{V}}\Diamond" description]{uuul} &[-10pt] {} &[-10pt] {} &[-10pt] {} &[-10pt] \left [2;0\right ]\ar[-, "\Diamond" near start]{uuur}\ar[-]{uuul} &[-10pt] {} &[-10pt] {} &[-10pt] {} &[-10pt] \left [3\right ]\ar[-]{uuur}\ar[-, "\Diamond2_{{}_{V}}\Diamond" description]{uuul} &[-10pt] {} \\ 
{} &[-10pt] {} &[-10pt] {} &[-10pt] {} &[-10pt] {} &[-10pt] {} &[-10pt] {} &[-10pt] {} &[-10pt] {} &[-10pt] {} &[-10pt] {} &[-10pt] {} &[-10pt] {} &[-10pt] {} &[-10pt] {} &[-10pt] {} \\ 
{} &[-10pt] {} &[-10pt] {} &[-10pt] {} &[-10pt] \left [1;0\right ]\ar[-, "\Diamond" near start]{uurr}\ar[-]{uull} &[-10pt] {} &[-10pt] {} &[-10pt] {} &[-10pt] {} &[-10pt] {} &[-10pt] {} &[-10pt] {} &[-10pt] \left [2\right ]\ar[-]{uurr}\ar[-, "\Diamond2_{{}_{V}}\Diamond" description]{uull} &[-10pt] {} &[-10pt] {} &[-10pt] {} \\ 
{} &[-10pt] {} &[-10pt] {} &[-10pt] {} &[-10pt] {} &[-10pt] {} &[-10pt] {} &[-10pt] {} &[-10pt] \left [1\right ]\ar[-]{urrrr}\ar[-, "\Diamond2_{{}_{V}}\Diamond" description]{ullll} &[-10pt] {} &[-10pt] {} &[-10pt] {} &[-10pt] {} &[-10pt] {} &[-10pt] {} &[-10pt] {} \\ 
\left [0\right ]\ar[-]{urrrrrrrr} &[-10pt] {} &[-10pt] {} &[-10pt] {} &[-10pt] {} &[-10pt] {} &[-10pt] {} &[-10pt] {} &[-10pt] {} &[-10pt] {} &[-10pt] {} &[-10pt] {} &[-10pt] {} &[-10pt] {} &[-10pt] {} &[-10pt] {} \\
\end{tikzcd}
};
\end{tikzpicture}
\caption{\label{naturalrepresentationstructure}The natural representations of the rational numbers specify their locations within the constituent structure of sets that contain many rational numbers.
}
\end{figure}

In figure \ref{rationalstructure}, each step from an integer across an edge labelled with $\Diamond2_{{}_{V}}\Diamond$ effectively subtracts ${}^1\!\!/_{\!\!2}$ from the integer, and further steps in that direction lead to a decreasing sequence of numbers.

The natural representations of the rational numbers are shown in figure \ref{naturalrepresentationstructure} within this structure. The fact that the integers that represent a rational number encode the route to it within the tree is visible from their pattern: All of the nodes that can be reached by traversing the $\Diamond2_{{}_{V}}\Diamond$ edge from 2 have 2 as their first entry, and the same is true of the other integers. 
The route to the number with the sequence $[1;1,0]$ goes through the number with the sequence $[1;1]$, and the route to that number goes through $[1]$.

The representation in figure \ref{naturalrepresentationstructure} makes it evident that which numbers are constituents of which others is indicated in their natural representation, with $2$ being a constituent of $[2; -1]$ and so on. In contrast, figure \ref{rationalstructure} shows that ${}^1\!\!/_{\!\!5}$ is not a constituent of ${}^3\!\!/_{\!\!5}$, and $5$ and $3$ are also not constituents of it, which reveals that the numerator and denominator of a rational number don't specify its constituents.

It is also worth noting that the standard continued fraction representation of a rational number also doesn't reveal its constituents. For example, ${}^1\!\!/_{\!\!3}$ is a constituent of ${}^2\!\!/_{\!\!5}$, but the standard continued fraction representation of ${}^1\!\!/_{\!\!3}$ is $[0;3]$, while that of ${}^2\!\!/_{\!\!5}$ is $[0;2,2]$.

Only the natural representation clearly expresses the constituency relation, with ${}^2\!\!/_{\!\!5}$ having the representation $[1;-1,0]$ which contains the representation of ${}^1\!\!/_{\!\!3}$, namely $[1;-1]$. 

The binary tree structure shown in figures \ref{rationalstructure} and \ref{naturalrepresentationstructure} is already known to mathematicians as the Stern-Brocot tree\cite{sternbrocot}. Many of its properties, such as the occurrence of each positive rational number exactly once in its simplest form, the ordering of numbers from left to right in accordance with their numerical order, and its relation to standard continued fractions, have been known for more than a century and a half\cite{stern,brocot}, as has its usefulness for identifying integers whose ratio best approximates a given number.

What has not previously been made clear is that the tree displays the structure of the rational numbers themselves. Without the concept of constituent structure, the relationship between the numbers in the tree has always seemed to be arithmetical and algorithmic.
Unlike the standard continued fraction representation, each rational number's natural representation is unique, and every finite sequence of integers in which all of them, apart from the first and the last, are required to be non-zero, identifies a unique rational number that can be reached by following that route in the tree.

\begin{landscape}
\begin{figure}
\begin{tikzpicture}[baseline= (a).base]
\node[scale=.20] (a) at (-100,0){\huge
\begin{tikzcd}
-4 &[-45pt] -7/2 &[-45pt] {} &[-45pt] -8/3 &[-45pt] {} &[-45pt] -7/3 &[-45pt] {} &[-45pt] -7/4 &[-45pt] {} &[-45pt] -8/5 &[-45pt] {} &[-45pt] -7/5 &[-45pt] {} &[-45pt] -5/4 &[-45pt] {} &[-45pt] -4/5 &[-45pt] {} &[-45pt] -5/7 &[-45pt] {} &[-45pt] -5/8 &[-45pt] {} &[-45pt] -4/7 &[-45pt] {} &[-45pt] -3/7 &[-45pt] {} &[-45pt] -3/8 &[-45pt] {} &[-45pt] -2/7 &[-45pt] {} &[-45pt] -1/5 &[-45pt] {} &[-45pt] 1/4 &[-45pt] {} &[-45pt] 2/5 &[-45pt] {} &[-45pt] 3/5 &[-45pt] {} &[-45pt] 3/4 &[-45pt] {} &[-45pt] 4/3 &[-45pt] {} &[-45pt] 5/3 &[-45pt] {} &[-45pt] 5/2 &[-45pt] {} &[-45pt] 4 \\ 
{} &[-45pt] {} &[-45pt] {} &[-45pt] {} &[-45pt] {} &[-45pt] {} &[-45pt] {} &[-45pt] {} &[-45pt] {} &[-45pt] {} &[-45pt] {} &[-45pt] {} &[-45pt] {} &[-45pt] {} &[-45pt] {} &[-45pt] {} &[-45pt] {} &[-45pt] {} &[-45pt] {} &[-45pt] {} &[-45pt] {} &[-45pt] {} &[-45pt] {} &[-45pt] {} &[-45pt] {} &[-45pt] {} &[-45pt] {} &[-45pt] {} &[-45pt] {} &[-45pt] {} &[-45pt] {} &[-45pt] {} &[-45pt] {} &[-45pt] {} &[-45pt] {} &[-45pt] {} &[-45pt] {} &[-45pt] {} &[-45pt] {} &[-45pt] {} &[-45pt] {} &[-45pt] {} &[-45pt] {} &[-45pt] {} &[-45pt] {} &[-45pt] {} \\ 
{} &[-45pt] {} &[-45pt] {} &[-45pt] {} &[-45pt] {} &[-45pt] {} &[-45pt] {} &[-45pt] {} &[-45pt] {} &[-45pt] {} &[-45pt] {} &[-45pt] {} &[-45pt] {} &[-45pt] {} &[-45pt] {} &[-45pt] {} &[-45pt] {} &[-45pt] {} &[-45pt] {} &[-45pt] {} &[-45pt] {} &[-45pt] {} &[-45pt] {} &[-45pt] {} &[-45pt] {} &[-45pt] {} &[-45pt] {} &[-45pt] {} &[-45pt] {} &[-45pt] {} &[-45pt] {} &[-45pt] {} &[-45pt] {} &[-45pt] {} &[-45pt] {} &[-45pt] {} &[-45pt] {} &[-45pt] {} &[-45pt] {} &[-45pt] {} &[-45pt] {} &[-45pt] {} &[-45pt] {} &[-45pt] {} &[-45pt] {} &[-45pt] {} \\ 
{} &[-45pt] {} &[-45pt] -3\ar[-]{uuull}\ar[-, "\Diamond2_{{}_{V}}\Diamond" description]{uuul} &[-45pt] {} &[-45pt] -5/2\ar[-, "\Diamond" near start]{uuur}\ar[-]{uuul} &[-45pt] {} &[-45pt] {} &[-45pt] {} &[-45pt] -5/3\ar[-]{uuul}\ar[-, "\Diamond2_{{}_{V}}\Diamond" description]{uuur} &[-45pt] {} &[-45pt] {} &[-45pt] {} &[-45pt] -4/3\ar[-]{uuur}\ar[-, "\Diamond2_{{}_{V}}\Diamond" description]{uuul} &[-45pt] {} &[-45pt] {} &[-45pt] {} &[-45pt] -3/4\ar[-]{uuul}\ar[-, "\Diamond2_{{}_{V}}\Diamond" description]{uuur} &[-45pt] {} &[-45pt] {} &[-45pt] {} &[-45pt] -3/5\ar[-, "\Diamond" near start]{uuul}\ar[-]{uuur} &[-45pt] {} &[-45pt] {} &[-45pt] {} &[-45pt] -2/5\ar[-, "\Diamond" near start]{uuur}\ar[-]{uuul} &[-45pt] {} &[-45pt] {} &[-45pt] {} &[-45pt] -1/4\ar[-]{uuur}\ar[-, "\Diamond2_{{}_{V}}\Diamond" description]{uuul} &[-45pt] {} &[-45pt] {} &[-45pt] {} &[-45pt] 1/3\ar[-]{uuul}\ar[-, "\Diamond2_{{}_{V}}\Diamond" description]{uuur} &[-45pt] {} &[-45pt] {} &[-45pt] {} &[-45pt] 2/3\ar[-]{uuur}\ar[-, "\Diamond2_{{}_{V}}\Diamond" description]{uuul} &[-45pt] {} &[-45pt] {} &[-45pt] {} &[-45pt] 3/2\ar[-, "\Diamond" near start]{uuur}\ar[-]{uuul} &[-45pt] {} &[-45pt] {} &[-45pt] {} &[-45pt] 3\ar[-]{uuur}\ar[-, "\Diamond2_{{}_{V}}\Diamond" description]{uuul} &[-45pt] {} \\ 
{} &[-45pt] {} &[-45pt] {} &[-45pt] {} &[-45pt] {} &[-45pt] {} &[-45pt] {} &[-45pt] {} &[-45pt] {} &[-45pt] {} &[-45pt] {} &[-45pt] {} &[-45pt] {} &[-45pt] {} &[-45pt] {} &[-45pt] {} &[-45pt] {} &[-45pt] {} &[-45pt] {} &[-45pt] {} &[-45pt] {} &[-45pt] {} &[-45pt] {} &[-45pt] {} &[-45pt] {} &[-45pt] {} &[-45pt] {} &[-45pt] {} &[-45pt] {} &[-45pt] {} &[-45pt] {} &[-45pt] {} &[-45pt] {} &[-45pt] {} &[-45pt] {} &[-45pt] {} &[-45pt] {} &[-45pt] {} &[-45pt] {} &[-45pt] {} &[-45pt] {} &[-45pt] {} &[-45pt] {} &[-45pt] {} &[-45pt] {} &[-45pt] {} \\ 
{} &[-45pt] {} &[-45pt] {} &[-45pt] {} &[-45pt] {} &[-45pt] {} &[-45pt] {} &[-45pt] {} &[-45pt] {} &[-45pt] {} &[-45pt] {} &[-45pt] {} &[-45pt] {} &[-45pt] {} &[-45pt] {} &[-45pt] {} &[-45pt] {} &[-45pt] {} &[-45pt] {} &[-45pt] {} &[-45pt] {} &[-45pt] {} &[-45pt] {} &[-45pt] {} &[-45pt] {} &[-45pt] {} &[-45pt] {} &[-45pt] {} &[-45pt] {} &[-45pt] {} &[-45pt] {} &[-45pt] {} &[-45pt] {} &[-45pt] {} &[-45pt] {} &[-45pt] {} &[-45pt] {} &[-45pt] {} &[-45pt] {} &[-45pt] {} &[-45pt] {} &[-45pt] {} &[-45pt] {} &[-45pt] {} &[-45pt] {} &[-45pt] {} \\ 
{} &[-45pt] {} &[-45pt] {} &[-45pt] {} &[-45pt] {} &[-45pt] {} &[-45pt] -2\ar[-]{uuullll}\ar[-, "\Diamond2_{{}_{V}}\Diamond" description]{uuull} &[-45pt] {} &[-45pt] {} &[-45pt] {} &[-45pt] -3/2\ar[-, "\Diamond" near start]{uuurr}\ar[-]{uuull} &[-45pt] {} &[-45pt] {} &[-45pt] {} &[-45pt] {} &[-45pt] {} &[-45pt] {} &[-45pt] {} &[-45pt] -2/3\ar[-]{uuull}\ar[-, "\Diamond2_{{}_{V}}\Diamond" description]{uuurr} &[-45pt] {} &[-45pt] {} &[-45pt] {} &[-45pt] {} &[-45pt] {} &[-45pt] {} &[-45pt] {} &[-45pt] -1/3\ar[-]{uuurr}\ar[-, "\Diamond2_{{}_{V}}\Diamond" description]{uuull} &[-45pt] {} &[-45pt] {} &[-45pt] {} &[-45pt] {} &[-45pt] {} &[-45pt] {} &[-45pt] {} &[-45pt] 1/2\ar[-, "\Diamond" near start]{uuurr}\ar[-]{uuull} &[-45pt] {} &[-45pt] {} &[-45pt] {} &[-45pt] {} &[-45pt] {} &[-45pt] {} &[-45pt] {} &[-45pt] 2\ar[-]{uuurr}\ar[-, "\Diamond2_{{}_{V}}\Diamond" description]{uuull} &[-45pt] {} &[-45pt] {} &[-45pt] {} \\ 
{} &[-45pt] {} &[-45pt] {} &[-45pt] {} &[-45pt] {} &[-45pt] {} &[-45pt] {} &[-45pt] {} &[-45pt] {} &[-45pt] {} &[-45pt] {} &[-45pt] {} &[-45pt] {} &[-45pt] {} &[-45pt] {} &[-45pt] {} &[-45pt] {} &[-45pt] {} &[-45pt] {} &[-45pt] {} &[-45pt] {} &[-45pt] {} &[-45pt] {} &[-45pt] {} &[-45pt] {} &[-45pt] {} &[-45pt] {} &[-45pt] {} &[-45pt] {} &[-45pt] {} &[-45pt] {} &[-45pt] {} &[-45pt] {} &[-45pt] {} &[-45pt] {} &[-45pt] {} &[-45pt] {} &[-45pt] {} &[-45pt] {} &[-45pt] {} &[-45pt] {} &[-45pt] {} &[-45pt] {} &[-45pt] {} &[-45pt] {} &[-45pt] {} \\ 
{} &[-45pt] {} &[-45pt] {} &[-45pt] {} &[-45pt] {} &[-45pt] {} &[-45pt] {} &[-45pt] {} &[-45pt] {} &[-45pt] {} &[-45pt] {} &[-45pt] {} &[-45pt] {} &[-45pt] {} &[-45pt] {} &[-45pt] {} &[-45pt] {} &[-45pt] {} &[-45pt] {} &[-45pt] {} &[-45pt] {} &[-45pt] {} &[-45pt] {} &[-45pt] {} &[-45pt] {} &[-45pt] {} &[-45pt] {} &[-45pt] {} &[-45pt] {} &[-45pt] {} &[-45pt] {} &[-45pt] {} &[-45pt] {} &[-45pt] {} &[-45pt] {} &[-45pt] {} &[-45pt] {} &[-45pt] {} &[-45pt] {} &[-45pt] {} &[-45pt] {} &[-45pt] {} &[-45pt] {} &[-45pt] {} &[-45pt] {} &[-45pt] {} \\ 
{} &[-45pt] {} &[-45pt] {} &[-45pt] {} &[-45pt] {} &[-45pt] {} &[-45pt] {} &[-45pt] {} &[-45pt] {} &[-45pt] {} &[-45pt] {} &[-45pt] {} &[-45pt] {} &[-45pt] {} &[-45pt] -1\ar[-]{uuullllllll}\ar[-, "\Diamond2_{{}_{V}}\Diamond" description]{uuullll} &[-45pt] {} &[-45pt] {} &[-45pt] {} &[-45pt] {} &[-45pt] {} &[-45pt] {} &[-45pt] {} &[-45pt] -1/2\ar[-, "\Diamond" near start]{uuurrrr}\ar[-]{uuullll} &[-45pt] {} &[-45pt] {} &[-45pt] {} &[-45pt] {} &[-45pt] {} &[-45pt] {} &[-45pt] {} &[-45pt] {} &[-45pt] {} &[-45pt] {} &[-45pt] {} &[-45pt] {} &[-45pt] {} &[-45pt] {} &[-45pt] {} &[-45pt] 1\ar[-]{uuurrrr}\ar[-, "\Diamond2_{{}_{V}}\Diamond" description]{uuullll} &[-45pt] {} &[-45pt] {} &[-45pt] {} &[-45pt] {} &[-45pt] {} &[-45pt] {} &[-45pt] {} \\ 
{} &[-45pt] {} &[-45pt] {} &[-45pt] {} &[-45pt] {} &[-45pt] {} &[-45pt] {} &[-45pt] {} &[-45pt] {} &[-45pt] {} &[-45pt] {} &[-45pt] {} &[-45pt] {} &[-45pt] {} &[-45pt] {} &[-45pt] {} &[-45pt] {} &[-45pt] {} &[-45pt] {} &[-45pt] {} &[-45pt] {} &[-45pt] {} &[-45pt] {} &[-45pt] {} &[-45pt] {} &[-45pt] {} &[-45pt] {} &[-45pt] {} &[-45pt] {} &[-45pt] {} &[-45pt] {} &[-45pt] {} &[-45pt] {} &[-45pt] {} &[-45pt] {} &[-45pt] {} &[-45pt] {} &[-45pt] {} &[-45pt] {} &[-45pt] {} &[-45pt] {} &[-45pt] {} &[-45pt] {} &[-45pt] {} &[-45pt] {} &[-45pt] {} \\ 
{} &[-45pt] {} &[-45pt] {} &[-45pt] {} &[-45pt] {} &[-45pt] {} &[-45pt] {} &[-45pt] {} &[-45pt] {} &[-45pt] {} &[-45pt] {} &[-45pt] {} &[-45pt] {} &[-45pt] {} &[-45pt] {} &[-45pt] {} &[-45pt] {} &[-45pt] {} &[-45pt] {} &[-45pt] {} &[-45pt] {} &[-45pt] {} &[-45pt] {} &[-45pt] {} &[-45pt] {} &[-45pt] {} &[-45pt] {} &[-45pt] {} &[-45pt] {} &[-45pt] {} &[-45pt] {} &[-45pt] {} &[-45pt] {} &[-45pt] {} &[-45pt] {} &[-45pt] {} &[-45pt] {} &[-45pt] {} &[-45pt] {} &[-45pt] {} &[-45pt] {} &[-45pt] {} &[-45pt] {} &[-45pt] {} &[-45pt] {} &[-45pt] {} \\ 
{} &[-45pt] {} &[-45pt] {} &[-45pt] {} &[-45pt] {} &[-45pt] {} &[-45pt] {} &[-45pt] {} &[-45pt] {} &[-45pt] {} &[-45pt] {} &[-45pt] {} &[-45pt] {} &[-45pt] {} &[-45pt] {} &[-45pt] {} &[-45pt] {} &[-45pt] {} &[-45pt] {} &[-45pt] {} &[-45pt] {} &[-45pt] {} &[-45pt] {} &[-45pt] {} &[-45pt] {} &[-45pt] {} &[-45pt] {} &[-45pt] {} &[-45pt] {} &[-45pt] {} &[-45pt] 0\ar[-]{uuurrrrrrrr}\ar[-, "\Diamond"']{uuullllllllllllllll}\ar[-, "\Diamond2_{{}_{V}}\Diamond" description]{uuullllllll} &[-45pt] {} &[-45pt] {} &[-45pt] {} &[-45pt] {} &[-45pt] {} &[-45pt] {} &[-45pt] {} &[-45pt] {} &[-45pt] {} &[-45pt] {} &[-45pt] {} &[-45pt] {} &[-45pt] {} &[-45pt] {} &[-45pt] {} \\ 
\end{tikzcd}
};
\end{tikzpicture}
\caption{\label{treewithnegativerationals}The Stern-Brocot tree extended to include all rational numbers, including negative numbers, as specified by their constituency relations and natural representations. Various symmetries are evident: Apart from integers, numbers symmetrically located around $-\nicefrac{1}{2}$ sum to $-1$; numbers symmetric around $1$ or $-1$ are reciprocals; numbers $a$ and $b$ which are symmetric around $0$ or $-2$ satisfy $a=\nicefrac{-b}{(1+b)}$. Symmetries extend outside the binary tree of the number around which the symmetry appears: $\nicefrac{5}{2}$ and $-\nicefrac{7}{2}$ sum to -$1$ despite being outside the binary tree that starts from $-\nicefrac{1}{2}$. The symmetries overlap: The reciprocals symmetric around $-1$ extend from $-\nicefrac{7}{2}$ to $-\nicefrac{2}{7}$, covering almost two thirds of the width of the tree, while the $\nicefrac{-b}{(1+b)}$ symmetry around $0$ extends from $-\nicefrac{4}{5}$ to $4$, covering two thirds of the tree and overlapping with the former symmetry in the middle third.}
\end{figure}
\end{landscape}

\begin{landscape}
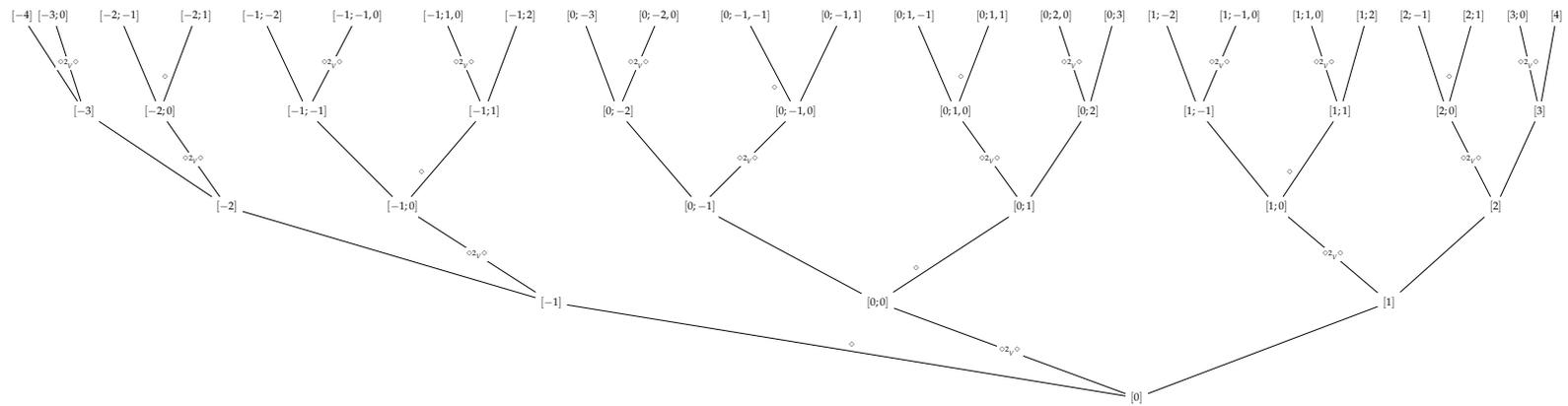
\begin{figure}
\begin{tikzpicture}[baseline= (a).base]
\node[scale=.20] (a) at (-100,0){\huge
\begin{tikzcd}
\left [-4\right ] &[-60pt] \left [-3;0\right ] &[-60pt] {} &[-60pt] \left [-2;-1\right ] &[-60pt] {} &[-60pt] \left [-2;1\right ] &[-60pt] {} &[-60pt] \left [-1;-2\right ] &[-60pt] {} &[-60pt] \left [-1;-1,0\right ] &[-60pt] {} &[-60pt] \left [-1;1,0\right ] &[-60pt] {} &[-60pt] \left [-1;2\right ] &[-60pt] {} &[-60pt] \left [0;-3\right ] &[-60pt] {} &[-60pt] \left [0;-2,0\right ] &[-60pt] {} &[-60pt] \left [0;-1,-1\right ] &[-60pt] {} &[-60pt] \left [0;-1,1\right ] &[-60pt] {} &[-60pt] \left [0;1,-1\right ] &[-60pt] {} &[-60pt] \left [0;1,1\right ] &[-60pt] {} &[-60pt] \left [0;2,0\right ] &[-60pt] {} &[-60pt] \left [0;3\right ] &[-60pt] {} &[-60pt] \left [1;-2\right ] &[-60pt] {} &[-60pt] \left [1;-1,0\right ] &[-60pt] {} &[-60pt] \left [1;1,0\right ] &[-60pt] {} &[-60pt] \left [1;2\right ] &[-60pt] {} &[-60pt] \left [2;-1\right ] &[-60pt] {} &[-60pt] \left [2;1\right ] &[-60pt] {} &[-60pt] \left [3;0\right ] &[-60pt] {} &[-60pt] \left [4\right ] \\ 
{} &[-60pt] {} &[-60pt] {} &[-60pt] {} &[-60pt] {} &[-60pt] {} &[-60pt] {} &[-60pt] {} &[-60pt] {} &[-60pt] {} &[-60pt] {} &[-60pt] {} &[-60pt] {} &[-60pt] {} &[-60pt] {} &[-60pt] {} &[-60pt] {} &[-60pt] {} &[-60pt] {} &[-60pt] {} &[-60pt] {} &[-60pt] {} &[-60pt] {} &[-60pt] {} &[-60pt] {} &[-60pt] {} &[-60pt] {} &[-60pt] {} &[-60pt] {} &[-60pt] {} &[-60pt] {} &[-60pt] {} &[-60pt] {} &[-60pt] {} &[-60pt] {} &[-60pt] {} &[-60pt] {} &[-60pt] {} &[-60pt] {} &[-60pt] {} &[-60pt] {} &[-60pt] {} &[-60pt] {} &[-60pt] {} &[-60pt] {} &[-60pt] {} \\ 
{} &[-60pt] {} &[-60pt] {} &[-60pt] {} &[-60pt] {} &[-60pt] {} &[-60pt] {} &[-60pt] {} &[-60pt] {} &[-60pt] {} &[-60pt] {} &[-60pt] {} &[-60pt] {} &[-60pt] {} &[-60pt] {} &[-60pt] {} &[-60pt] {} &[-60pt] {} &[-60pt] {} &[-60pt] {} &[-60pt] {} &[-60pt] {} &[-60pt] {} &[-60pt] {} &[-60pt] {} &[-60pt] {} &[-60pt] {} &[-60pt] {} &[-60pt] {} &[-60pt] {} &[-60pt] {} &[-60pt] {} &[-60pt] {} &[-60pt] {} &[-60pt] {} &[-60pt] {} &[-60pt] {} &[-60pt] {} &[-60pt] {} &[-60pt] {} &[-60pt] {} &[-60pt] {} &[-60pt] {} &[-60pt] {} &[-60pt] {} &[-60pt] {} \\ 
{} &[-60pt] {} &[-60pt] \left [-3\right ]\ar[-]{uuull}\ar[-, "\Diamond2_{{}_{V}}\Diamond" description]{uuul} &[-60pt] {} &[-60pt] \left [-2;0\right ]\ar[-, "\Diamond" near start]{uuur}\ar[-]{uuul} &[-60pt] {} &[-60pt] {} &[-60pt] {} &[-60pt] \left [-1;-1\right ]\ar[-]{uuul}\ar[-, "\Diamond2_{{}_{V}}\Diamond" description]{uuur} &[-60pt] {} &[-60pt] {} &[-60pt] {} &[-60pt] \left [-1;1\right ]\ar[-]{uuur}\ar[-, "\Diamond2_{{}_{V}}\Diamond" description]{uuul} &[-60pt] {} &[-60pt] {} &[-60pt] {} &[-60pt] \left [0;-2\right ]\ar[-]{uuul}\ar[-, "\Diamond2_{{}_{V}}\Diamond" description]{uuur} &[-60pt] {} &[-60pt] {} &[-60pt] {} &[-60pt] \left [0;-1,0\right ]\ar[-, "\Diamond" near start]{uuul}\ar[-]{uuur} &[-60pt] {} &[-60pt] {} &[-60pt] {} &[-60pt] \left [0;1,0\right ]\ar[-, "\Diamond" near start]{uuur}\ar[-]{uuul} &[-60pt] {} &[-60pt] {} &[-60pt] {} &[-60pt] \left [0;2\right ]\ar[-]{uuur}\ar[-, "\Diamond2_{{}_{V}}\Diamond" description]{uuul} &[-60pt] {} &[-60pt] {} &[-60pt] {} &[-60pt] \left [1;-1\right ]\ar[-]{uuul}\ar[-, "\Diamond2_{{}_{V}}\Diamond" description]{uuur} &[-60pt] {} &[-60pt] {} &[-60pt] {} &[-60pt] \left [1;1\right ]\ar[-]{uuur}\ar[-, "\Diamond2_{{}_{V}}\Diamond" description]{uuul} &[-60pt] {} &[-60pt] {} &[-60pt] {} &[-60pt] \left [2;0\right ]\ar[-, "\Diamond" near start]{uuur}\ar[-]{uuul} &[-60pt] {} &[-60pt] {} &[-60pt] {} &[-60pt] \left [3\right ]\ar[-]{uuur}\ar[-, "\Diamond2_{{}_{V}}\Diamond" description]{uuul} &[-60pt] {} \\ 
{} &[-60pt] {} &[-60pt] {} &[-60pt] {} &[-60pt] {} &[-60pt] {} &[-60pt] {} &[-60pt] {} &[-60pt] {} &[-60pt] {} &[-60pt] {} &[-60pt] {} &[-60pt] {} &[-60pt] {} &[-60pt] {} &[-60pt] {} &[-60pt] {} &[-60pt] {} &[-60pt] {} &[-60pt] {} &[-60pt] {} &[-60pt] {} &[-60pt] {} &[-60pt] {} &[-60pt] {} &[-60pt] {} &[-60pt] {} &[-60pt] {} &[-60pt] {} &[-60pt] {} &[-60pt] {} &[-60pt] {} &[-60pt] {} &[-60pt] {} &[-60pt] {} &[-60pt] {} &[-60pt] {} &[-60pt] {} &[-60pt] {} &[-60pt] {} &[-60pt] {} &[-60pt] {} &[-60pt] {} &[-60pt] {} &[-60pt] {} &[-60pt] {} \\ 
{} &[-60pt] {} &[-60pt] {} &[-60pt] {} &[-60pt] {} &[-60pt] {} &[-60pt] {} &[-60pt] {} &[-60pt] {} &[-60pt] {} &[-60pt] {} &[-60pt] {} &[-60pt] {} &[-60pt] {} &[-60pt] {} &[-60pt] {} &[-60pt] {} &[-60pt] {} &[-60pt] {} &[-60pt] {} &[-60pt] {} &[-60pt] {} &[-60pt] {} &[-60pt] {} &[-60pt] {} &[-60pt] {} &[-60pt] {} &[-60pt] {} &[-60pt] {} &[-60pt] {} &[-60pt] {} &[-60pt] {} &[-60pt] {} &[-60pt] {} &[-60pt] {} &[-60pt] {} &[-60pt] {} &[-60pt] {} &[-60pt] {} &[-60pt] {} &[-60pt] {} &[-60pt] {} &[-60pt] {} &[-60pt] {} &[-60pt] {} &[-60pt] {} \\ 
{} &[-60pt] {} &[-60pt] {} &[-60pt] {} &[-60pt] {} &[-60pt] {} &[-60pt] \left [-2\right ]\ar[-]{uuullll}\ar[-, "\Diamond2_{{}_{V}}\Diamond" description]{uuull} &[-60pt] {} &[-60pt] {} &[-60pt] {} &[-60pt] \left [-1;0\right ]\ar[-, "\Diamond" near start]{uuurr}\ar[-]{uuull} &[-60pt] {} &[-60pt] {} &[-60pt] {} &[-60pt] {} &[-60pt] {} &[-60pt] {} &[-60pt] {} &[-60pt] \left [0;-1\right ]\ar[-]{uuull}\ar[-, "\Diamond2_{{}_{V}}\Diamond" description]{uuurr} &[-60pt] {} &[-60pt] {} &[-60pt] {} &[-60pt] {} &[-60pt] {} &[-60pt] {} &[-60pt] {} &[-60pt] \left [0;1\right ]\ar[-]{uuurr}\ar[-, "\Diamond2_{{}_{V}}\Diamond" description]{uuull} &[-60pt] {} &[-60pt] {} &[-60pt] {} &[-60pt] {} &[-60pt] {} &[-60pt] {} &[-60pt] {} &[-60pt] \left [1;0\right ]\ar[-, "\Diamond" near start]{uuurr}\ar[-]{uuull} &[-60pt] {} &[-60pt] {} &[-60pt] {} &[-60pt] {} &[-60pt] {} &[-60pt] {} &[-60pt] {} &[-60pt] \left [2\right ]\ar[-]{uuurr}\ar[-, "\Diamond2_{{}_{V}}\Diamond" description]{uuull} &[-60pt] {} &[-60pt] {} &[-60pt] {} \\ 
{} &[-60pt] {} &[-60pt] {} &[-60pt] {} &[-60pt] {} &[-60pt] {} &[-60pt] {} &[-60pt] {} &[-60pt] {} &[-60pt] {} &[-60pt] {} &[-60pt] {} &[-60pt] {} &[-60pt] {} &[-60pt] {} &[-60pt] {} &[-60pt] {} &[-60pt] {} &[-60pt] {} &[-60pt] {} &[-60pt] {} &[-60pt] {} &[-60pt] {} &[-60pt] {} &[-60pt] {} &[-60pt] {} &[-60pt] {} &[-60pt] {} &[-60pt] {} &[-60pt] {} &[-60pt] {} &[-60pt] {} &[-60pt] {} &[-60pt] {} &[-60pt] {} &[-60pt] {} &[-60pt] {} &[-60pt] {} &[-60pt] {} &[-60pt] {} &[-60pt] {} &[-60pt] {} &[-60pt] {} &[-60pt] {} &[-60pt] {} &[-60pt] {} \\ 
{} &[-60pt] {} &[-60pt] {} &[-60pt] {} &[-60pt] {} &[-60pt] {} &[-60pt] {} &[-60pt] {} &[-60pt] {} &[-60pt] {} &[-60pt] {} &[-60pt] {} &[-60pt] {} &[-60pt] {} &[-60pt] {} &[-60pt] {} &[-60pt] {} &[-60pt] {} &[-60pt] {} &[-60pt] {} &[-60pt] {} &[-60pt] {} &[-60pt] {} &[-60pt] {} &[-60pt] {} &[-60pt] {} &[-60pt] {} &[-60pt] {} &[-60pt] {} &[-60pt] {} &[-60pt] {} &[-60pt] {} &[-60pt] {} &[-60pt] {} &[-60pt] {} &[-60pt] {} &[-60pt] {} &[-60pt] {} &[-60pt] {} &[-60pt] {} &[-60pt] {} &[-60pt] {} &[-60pt] {} &[-60pt] {} &[-60pt] {} &[-60pt] {} \\ 
{} &[-60pt] {} &[-60pt] {} &[-60pt] {} &[-60pt] {} &[-60pt] {} &[-60pt] {} &[-60pt] {} &[-60pt] {} &[-60pt] {} &[-60pt] {} &[-60pt] {} &[-60pt] {} &[-60pt] {} &[-60pt] \left [-1\right ]\ar[-]{uuullllllll}\ar[-, "\Diamond2_{{}_{V}}\Diamond" description]{uuullll} &[-60pt] {} &[-60pt] {} &[-60pt] {} &[-60pt] {} &[-60pt] {} &[-60pt] {} &[-60pt] {} &[-60pt] \left [0;0\right ]\ar[-, "\Diamond" near start]{uuurrrr}\ar[-]{uuullll} &[-60pt] {} &[-60pt] {} &[-60pt] {} &[-60pt] {} &[-60pt] {} &[-60pt] {} &[-60pt] {} &[-60pt] {} &[-60pt] {} &[-60pt] {} &[-60pt] {} &[-60pt] {} &[-60pt] {} &[-60pt] {} &[-60pt] {} &[-60pt] \left [1\right ]\ar[-]{uuurrrr}\ar[-, "\Diamond2_{{}_{V}}\Diamond" description]{uuullll} &[-60pt] {} &[-60pt] {} &[-60pt] {} &[-60pt] {} &[-60pt] {} &[-60pt] {} &[-60pt] {} \\ 
{} &[-60pt] {} &[-60pt] {} &[-60pt] {} &[-60pt] {} &[-60pt] {} &[-60pt] {} &[-60pt] {} &[-60pt] {} &[-60pt] {} &[-60pt] {} &[-60pt] {} &[-60pt] {} &[-60pt] {} &[-60pt] {} &[-60pt] {} &[-60pt] {} &[-60pt] {} &[-60pt] {} &[-60pt] {} &[-60pt] {} &[-60pt] {} &[-60pt] {} &[-60pt] {} &[-60pt] {} &[-60pt] {} &[-60pt] {} &[-60pt] {} &[-60pt] {} &[-60pt] {} &[-60pt] {} &[-60pt] {} &[-60pt] {} &[-60pt] {} &[-60pt] {} &[-60pt] {} &[-60pt] {} &[-60pt] {} &[-60pt] {} &[-60pt] {} &[-60pt] {} &[-60pt] {} &[-60pt] {} &[-60pt] {} &[-60pt] {} &[-60pt] {} \\ 
{} &[-60pt] {} &[-60pt] {} &[-60pt] {} &[-60pt] {} &[-60pt] {} &[-60pt] {} &[-60pt] {} &[-60pt] {} &[-60pt] {} &[-60pt] {} &[-60pt] {} &[-60pt] {} &[-60pt] {} &[-60pt] {} &[-60pt] {} &[-60pt] {} &[-60pt] {} &[-60pt] {} &[-60pt] {} &[-60pt] {} &[-60pt] {} &[-60pt] {} &[-60pt] {} &[-60pt] {} &[-60pt] {} &[-60pt] {} &[-60pt] {} &[-60pt] {} &[-60pt] {} &[-60pt] {} &[-60pt] {} &[-60pt] {} &[-60pt] {} &[-60pt] {} &[-60pt] {} &[-60pt] {} &[-60pt] {} &[-60pt] {} &[-60pt] {} &[-60pt] {} &[-60pt] {} &[-60pt] {} &[-60pt] {} &[-60pt] {} &[-60pt] {} \\ 
{} &[-60pt] {} &[-60pt] {} &[-60pt] {} &[-60pt] {} &[-60pt] {} &[-60pt] {} &[-60pt] {} &[-60pt] {} &[-60pt] {} &[-60pt] {} &[-60pt] {} &[-60pt] {} &[-60pt] {} &[-60pt] {} &[-60pt] {} &[-60pt] {} &[-60pt] {} &[-60pt] {} &[-60pt] {} &[-60pt] {} &[-60pt] {} &[-60pt] {} &[-60pt] {} &[-60pt] {} &[-60pt] {} &[-60pt] {} &[-60pt] {} &[-60pt] {} &[-60pt] {} &[-60pt] \left [0\right ]\ar[-]{uuurrrrrrrr}\ar[-, "\Diamond"']{uuullllllllllllllll}\ar[-, "\Diamond2_{{}_{V}}\Diamond" description]{uuullllllll} &[-60pt] {} &[-60pt] {} &[-60pt] {} &[-60pt] {} &[-60pt] {} &[-60pt] {} &[-60pt] {} &[-60pt] {} &[-60pt] {} &[-60pt] {} &[-60pt] {} &[-60pt] {} &[-60pt] {} &[-60pt] {} &[-60pt] {} \\ 
\end{tikzcd}
};
\end{tikzpicture}
\caption{\label{naturalrepresentationtreewithnegativerationals}The natural representations of positive and negative rational numbers. Note that, like the rational numbers themselves, the lists of integers appear in increasing order from left to right, when the earliest integers in each list have the highest priority in determining their order. $[-4]$ is to left of $[-3;0]$ because $-4<-3$ and so on. The tree is symmetric around $[0;0]=-\nicefrac{1}{2}$ from that level upwards. Sequences symmetrically located around $[0;0]$ have opposite signs. The height of a sequence within the tree is equal to the sum of the absolute values of its entries plus its length. 
}
\end{figure}
\end{landscape}

\begin{figure}
\centering
\begin{tikzpicture}[baseline= (a).base]
\node[scale=.45] (a) at (-100,0){
\begin{tikzcd}
-4 &[-20pt] {} &[-20pt] -5/2 &[-20pt] {} &[-20pt] -5/3 &[-20pt] {} &[-20pt] -4/3 &[-20pt] {} &[-20pt] -3/4 &[-20pt] {} &[-20pt] -3/5 &[-20pt] {} &[-20pt] -2/5 &[-20pt] {} &[-20pt] -1/4 &[-20pt] {} &[-20pt] 1/4 &[-20pt] {} &[-20pt] 2/5 &[-20pt] {} &[-20pt] 3/5 &[-20pt] {} &[-20pt] 3/4 &[-20pt] {} &[-20pt] 4/3 &[-20pt] {} &[-20pt] 5/3 &[-20pt] {} &[-20pt] 5/2 &[-20pt] {} &[-20pt] 4 \\ 
{} &[-20pt] {} &[-20pt] {} &[-20pt] {} &[-20pt] {} &[-20pt] {} &[-20pt] {} &[-20pt] {} &[-20pt] {} &[-20pt] {} &[-20pt] {} &[-20pt] {} &[-20pt] {} &[-20pt] {} &[-20pt] {} &[-20pt] {} &[-20pt] {} &[-20pt] {} &[-20pt] {} &[-20pt] {} &[-20pt] {} &[-20pt] {} &[-20pt] {} &[-20pt] {} &[-20pt] {} &[-20pt] {} &[-20pt] {} &[-20pt] {} &[-20pt] {} &[-20pt] {} &[-20pt] {} \\ 
{} &[-20pt] {} &[-20pt] {} &[-20pt] {} &[-20pt] {} &[-20pt] {} &[-20pt] {} &[-20pt] {} &[-20pt] {} &[-20pt] {} &[-20pt] {} &[-20pt] {} &[-20pt] {} &[-20pt] {} &[-20pt] {} &[-20pt] {} &[-20pt] {} &[-20pt] {} &[-20pt] {} &[-20pt] {} &[-20pt] {} &[-20pt] {} &[-20pt] {} &[-20pt] {} &[-20pt] {} &[-20pt] {} &[-20pt] {} &[-20pt] {} &[-20pt] {} &[-20pt] {} &[-20pt] {} \\ 
{} &[-20pt] -3\ar[-]{uuul}\ar[-]{uuur} &[-20pt] {} &[-20pt] {} &[-20pt] {} &[-20pt] -3/2\ar[-]{uuul}\ar[-]{uuur} &[-20pt] {} &[-20pt] {} &[-20pt] {} &[-20pt] -2/3\ar[-]{uuul}\ar[-]{uuur} &[-20pt] {} &[-20pt] {} &[-20pt] {} &[-20pt] -1/3\ar[-]{uuur}\ar[-]{uuul} &[-20pt] {} &[-20pt] {} &[-20pt] {} &[-20pt] 1/3\ar[-]{uuul}\ar[-]{uuur} &[-20pt] {} &[-20pt] {} &[-20pt] {} &[-20pt] 2/3\ar[-]{uuur}\ar[-]{uuul} &[-20pt] {} &[-20pt] {} &[-20pt] {} &[-20pt] 3/2\ar[-]{uuur}\ar[-]{uuul} &[-20pt] {} &[-20pt] {} &[-20pt] {} &[-20pt] 3\ar[-]{uuur}\ar[-]{uuul} &[-20pt] {} \\ 
{} &[-20pt] {} &[-20pt] {} &[-20pt] {} &[-20pt] {} &[-20pt] {} &[-20pt] {} &[-20pt] {} &[-20pt] {} &[-20pt] {} &[-20pt] {} &[-20pt] {} &[-20pt] {} &[-20pt] {} &[-20pt] {} &[-20pt] {} &[-20pt] {} &[-20pt] {} &[-20pt] {} &[-20pt] {} &[-20pt] {} &[-20pt] {} &[-20pt] {} &[-20pt] {} &[-20pt] {} &[-20pt] {} &[-20pt] {} &[-20pt] {} &[-20pt] {} &[-20pt] {} &[-20pt] {} \\ 
{} &[-20pt] {} &[-20pt] {} &[-20pt] {} &[-20pt] {} &[-20pt] {} &[-20pt] {} &[-20pt] {} &[-20pt] {} &[-20pt] {} &[-20pt] {} &[-20pt] {} &[-20pt] {} &[-20pt] {} &[-20pt] {} &[-20pt] {} &[-20pt] {} &[-20pt] {} &[-20pt] {} &[-20pt] {} &[-20pt] {} &[-20pt] {} &[-20pt] {} &[-20pt] {} &[-20pt] {} &[-20pt] {} &[-20pt] {} &[-20pt] {} &[-20pt] {} &[-20pt] {} &[-20pt] {} \\ 
{} &[-20pt] {} &[-20pt] {} &[-20pt] -2\ar[-]{uuull}\ar[-]{uuurr} &[-20pt] {} &[-20pt] {} &[-20pt] {} &[-20pt] {} &[-20pt] {} &[-20pt] {} &[-20pt] {} &[-20pt] -1/2\ar[-]{uuull}\ar[-]{uuurr} &[-20pt] {} &[-20pt] {} &[-20pt] {} &[-20pt] {} &[-20pt] {} &[-20pt] {} &[-20pt] {} &[-20pt] 1/2\ar[-]{uuurr}\ar[-]{uuull} &[-20pt] {} &[-20pt] {} &[-20pt] {} &[-20pt] {} &[-20pt] {} &[-20pt] {} &[-20pt] {} &[-20pt] 2\ar[-]{uuurr}\ar[-]{uuull} &[-20pt] {} &[-20pt] {} &[-20pt] {} \\ 
{} &[-20pt] {} &[-20pt] {} &[-20pt] {} &[-20pt] {} &[-20pt] {} &[-20pt] {} &[-20pt] {} &[-20pt] {} &[-20pt] {} &[-20pt] {} &[-20pt] {} &[-20pt] {} &[-20pt] {} &[-20pt] {} &[-20pt] {} &[-20pt] {} &[-20pt] {} &[-20pt] {} &[-20pt] {} &[-20pt] {} &[-20pt] {} &[-20pt] {} &[-20pt] {} &[-20pt] {} &[-20pt] {} &[-20pt] {} &[-20pt] {} &[-20pt] {} &[-20pt] {} &[-20pt] {} \\ 
{} &[-20pt] {} &[-20pt] {} &[-20pt] {} &[-20pt] {} &[-20pt] {} &[-20pt] {} &[-20pt] {} &[-20pt] {} &[-20pt] {} &[-20pt] {} &[-20pt] {} &[-20pt] {} &[-20pt] {} &[-20pt] {} &[-20pt] {} &[-20pt] {} &[-20pt] {} &[-20pt] {} &[-20pt] {} &[-20pt] {} &[-20pt] {} &[-20pt] {} &[-20pt] {} &[-20pt] {} &[-20pt] {} &[-20pt] {} &[-20pt] {} &[-20pt] {} &[-20pt] {} &[-20pt] {} \\ 
{} &[-20pt] {} &[-20pt] {} &[-20pt] {} &[-20pt] {} &[-20pt] {} &[-20pt] {} &[-20pt] -1\ar[-]{uuullll}\ar[-]{uuurrrr} &[-20pt] {} &[-20pt] {} &[-20pt] {} &[-20pt] {} &[-20pt] {} &[-20pt] {} &[-20pt] {} &[-20pt] {} &[-20pt] {} &[-20pt] {} &[-20pt] {} &[-20pt] {} &[-20pt] {} &[-20pt] {} &[-20pt] {} &[-20pt] 1\ar[-]{uuurrrr}\ar[-]{uuullll} &[-20pt] {} &[-20pt] {} &[-20pt] {} &[-20pt] {} &[-20pt] {} &[-20pt] {} &[-20pt] {} \\ 
{} &[-20pt] {} &[-20pt] {} &[-20pt] {} &[-20pt] {} &[-20pt] {} &[-20pt] {} &[-20pt] {} &[-20pt] {} &[-20pt] {} &[-20pt] {} &[-20pt] {} &[-20pt] {} &[-20pt] {} &[-20pt] {} &[-20pt] {} &[-20pt] {} &[-20pt] {} &[-20pt] {} &[-20pt] {} &[-20pt] {} &[-20pt] {} &[-20pt] {} &[-20pt] {} &[-20pt] {} &[-20pt] {} &[-20pt] {} &[-20pt] {} &[-20pt] {} &[-20pt] {} &[-20pt] {} \\ 
{} &[-20pt] {} &[-20pt] {} &[-20pt] {} &[-20pt] {} &[-20pt] {} &[-20pt] {} &[-20pt] {} &[-20pt] {} &[-20pt] {} &[-20pt] {} &[-20pt] {} &[-20pt] {} &[-20pt] {} &[-20pt] {} &[-20pt] {} &[-20pt] {} &[-20pt] {} &[-20pt] {} &[-20pt] {} &[-20pt] {} &[-20pt] {} &[-20pt] {} &[-20pt] {} &[-20pt] {} &[-20pt] {} &[-20pt] {} &[-20pt] {} &[-20pt] {} &[-20pt] {} &[-20pt] {} \\ 
{} &[-20pt] {} &[-20pt] {} &[-20pt] {} &[-20pt] {} &[-20pt] {} &[-20pt] {} &[-20pt] {} &[-20pt] {} &[-20pt] {} &[-20pt] {} &[-20pt] {} &[-20pt] {} &[-20pt] {} &[-20pt] {} &[-20pt] 0\ar[-]{uuurrrrrrrr}\ar[-]{uuullllllll} &[-20pt] {} &[-20pt] {} &[-20pt] {} &[-20pt] {} &[-20pt] {} &[-20pt] {} &[-20pt] {} &[-20pt] {} &[-20pt] {} &[-20pt] {} &[-20pt] {} &[-20pt] {} &[-20pt] {} &[-20pt] {} &[-20pt] {} \\ 
\end{tikzcd}
};
\end{tikzpicture}
\caption{\label{simplertree}
The most intuitively obvious way to extend the Stern-Brocot tree to negative rationals is to take the mirror image of the positive tree and connect the vertices at $1$ and $-1$ to a single root at $0$. The resulting tree contains all rationals in the right order, but is not generated in its entirety by any single procedure, has ad-hoc rules for the vertex at zero, different from the rules for other vertices, and has no overlapping symmetries. Only the numbers within the local binary tree originating from a given number exhibit a symmetry around that number: Numbers symmetrically positioned around $1$ are inverses; and the same is true for $-1$, but there are no numbers that participate in both symmetries. Numbers symmetric around $-\nicefrac{1}{2}$ sum to $-1$, but only if they are within the binary tree originating from $-\nicefrac{1}{2}$.}

\vspace{9mm}

\centering
\begin{tikzpicture}[baseline= (a).base]
\node[scale=.45] (a) at (-100,0){
\begin{tikzcd}
{} &[-25pt] -9/5 &[-25pt] {} &[-25pt] -12/7 &[-25pt] {} &[-25pt] -13/8 &[-25pt] {} &[-25pt] -11/7 &[-25pt] {} &[-25pt] -10/7 &[-25pt] {} &[-25pt] -11/8 &[-25pt] {} &[-25pt] -9/7 &[-25pt] {} &[-25pt] -6/5 &[-25pt] {} &[-25pt] -3/4 &[-25pt] {} &[-25pt] -3/5 &[-25pt] {} &[-25pt] -2/5 &[-25pt] {} &[-25pt] -1/4 &[-25pt] {} &[-25pt] 1/3 &[-25pt] {} &[-25pt] 2/3 &[-25pt] {} &[-25pt] 3/2 &[-25pt] {} &[-25pt] 3 \\ 
{} &[-25pt] {} &[-25pt] {} &[-25pt] {} &[-25pt] {} &[-25pt] {} &[-25pt] {} &[-25pt] {} &[-25pt] {} &[-25pt] {} &[-25pt] {} &[-25pt] {} &[-25pt] {} &[-25pt] {} &[-25pt] {} &[-25pt] {} &[-25pt] {} &[-25pt] {} &[-25pt] {} &[-25pt] {} &[-25pt] {} &[-25pt] {} &[-25pt] {} &[-25pt] {} &[-25pt] {} &[-25pt] {} &[-25pt] {} &[-25pt] {} &[-25pt] {} &[-25pt] {} &[-25pt] {} &[-25pt] {} \\ 
{} &[-25pt] {} &[-25pt] {} &[-25pt] {} &[-25pt] {} &[-25pt] {} &[-25pt] {} &[-25pt] {} &[-25pt] {} &[-25pt] {} &[-25pt] {} &[-25pt] {} &[-25pt] {} &[-25pt] {} &[-25pt] {} &[-25pt] {} &[-25pt] {} &[-25pt] {} &[-25pt] {} &[-25pt] {} &[-25pt] {} &[-25pt] {} &[-25pt] {} &[-25pt] {} &[-25pt] {} &[-25pt] {} &[-25pt] {} &[-25pt] {} &[-25pt] {} &[-25pt] {} &[-25pt] {} &[-25pt] {} \\ 
{} &[-25pt] {} &[-25pt] -7/4\ar[-]{uuul}\ar[-]{uuur} &[-25pt] {} &[-25pt] {} &[-25pt] {} &[-25pt] -8/5\ar[-]{uuur}\ar[-]{uuul} &[-25pt] {} &[-25pt] {} &[-25pt] {} &[-25pt] -7/5\ar[-]{uuur}\ar[-]{uuul} &[-25pt] {} &[-25pt] {} &[-25pt] {} &[-25pt] -5/4\ar[-]{uuur}\ar[-]{uuul} &[-25pt] {} &[-25pt] {} &[-25pt] {} &[-25pt] -2/3\ar[-]{uuul}\ar[-]{uuur} &[-25pt] {} &[-25pt] {} &[-25pt] {} &[-25pt] -1/3\ar[-]{uuur}\ar[-]{uuul} &[-25pt] {} &[-25pt] {} &[-25pt] {} &[-25pt] 1/2\ar[-]{uuur}\ar[-]{uuul} &[-25pt] {} &[-25pt] {} &[-25pt] {} &[-25pt] 2\ar[-]{uuur}\ar[-]{uuul} &[-25pt] {} \\ 
{} &[-25pt] {} &[-25pt] {} &[-25pt] {} &[-25pt] {} &[-25pt] {} &[-25pt] {} &[-25pt] {} &[-25pt] {} &[-25pt] {} &[-25pt] {} &[-25pt] {} &[-25pt] {} &[-25pt] {} &[-25pt] {} &[-25pt] {} &[-25pt] {} &[-25pt] {} &[-25pt] {} &[-25pt] {} &[-25pt] {} &[-25pt] {} &[-25pt] {} &[-25pt] {} &[-25pt] {} &[-25pt] {} &[-25pt] {} &[-25pt] {} &[-25pt] {} &[-25pt] {} &[-25pt] {} &[-25pt] {} \\ 
{} &[-25pt] {} &[-25pt] {} &[-25pt] {} &[-25pt] {} &[-25pt] {} &[-25pt] {} &[-25pt] {} &[-25pt] {} &[-25pt] {} &[-25pt] {} &[-25pt] {} &[-25pt] {} &[-25pt] {} &[-25pt] {} &[-25pt] {} &[-25pt] {} &[-25pt] {} &[-25pt] {} &[-25pt] {} &[-25pt] {} &[-25pt] {} &[-25pt] {} &[-25pt] {} &[-25pt] {} &[-25pt] {} &[-25pt] {} &[-25pt] {} &[-25pt] {} &[-25pt] {} &[-25pt] {} &[-25pt] {} \\ 
{} &[-25pt] {} &[-25pt] {} &[-25pt] {} &[-25pt] -5/3\ar[-]{uuull}\ar[-]{uuurr} &[-25pt] {} &[-25pt] {} &[-25pt] {} &[-25pt] {} &[-25pt] {} &[-25pt] {} &[-25pt] {} &[-25pt] -4/3\ar[-]{uuurr}\ar[-]{uuull} &[-25pt] {} &[-25pt] {} &[-25pt] {} &[-25pt] {} &[-25pt] {} &[-25pt] {} &[-25pt] {} &[-25pt] -1/2\ar[-]{uuurr}\ar[-]{uuull} &[-25pt] {} &[-25pt] {} &[-25pt] {} &[-25pt] {} &[-25pt] {} &[-25pt] {} &[-25pt] {} &[-25pt] 1\ar[-]{uuurr}\ar[-]{uuull} &[-25pt] {} &[-25pt] {} &[-25pt] {} \\ 
{} &[-25pt] {} &[-25pt] {} &[-25pt] {} &[-25pt] {} &[-25pt] {} &[-25pt] {} &[-25pt] {} &[-25pt] {} &[-25pt] {} &[-25pt] {} &[-25pt] {} &[-25pt] {} &[-25pt] {} &[-25pt] {} &[-25pt] {} &[-25pt] {} &[-25pt] {} &[-25pt] {} &[-25pt] {} &[-25pt] {} &[-25pt] {} &[-25pt] {} &[-25pt] {} &[-25pt] {} &[-25pt] {} &[-25pt] {} &[-25pt] {} &[-25pt] {} &[-25pt] {} &[-25pt] {} &[-25pt] {} \\ 
{} &[-25pt] {} &[-25pt] {} &[-25pt] {} &[-25pt] {} &[-25pt] {} &[-25pt] {} &[-25pt] {} &[-25pt] {} &[-25pt] {} &[-25pt] {} &[-25pt] {} &[-25pt] {} &[-25pt] {} &[-25pt] {} &[-25pt] {} &[-25pt] {} &[-25pt] {} &[-25pt] {} &[-25pt] {} &[-25pt] {} &[-25pt] {} &[-25pt] {} &[-25pt] {} &[-25pt] {} &[-25pt] {} &[-25pt] {} &[-25pt] {} &[-25pt] {} &[-25pt] {} &[-25pt] {} &[-25pt] {} \\ 
{} &[-25pt] {} &[-25pt] {} &[-25pt] {} &[-25pt] {} &[-25pt] {} &[-25pt] {} &[-25pt] {} &[-25pt] -3/2\ar[-]{uuurrrr}\ar[-]{uuullll} &[-25pt] {} &[-25pt] {} &[-25pt] {} &[-25pt] {} &[-25pt] {} &[-25pt] {} &[-25pt] {} &[-25pt] {} &[-25pt] {} &[-25pt] {} &[-25pt] {} &[-25pt] {} &[-25pt] {} &[-25pt] {} &[-25pt] {} &[-25pt] 0\ar[-]{uuurrrr}\ar[-]{uuullll} &[-25pt] {} &[-25pt] {} &[-25pt] {} &[-25pt] {} &[-25pt] {} &[-25pt] {} &[-25pt] {} \\ 
{} &[-25pt] {} &[-25pt] {} &[-25pt] {} &[-25pt] {} &[-25pt] {} &[-25pt] {} &[-25pt] {} &[-25pt] {} &[-25pt] {} &[-25pt] {} &[-25pt] {} &[-25pt] {} &[-25pt] {} &[-25pt] {} &[-25pt] {} &[-25pt] {} &[-25pt] {} &[-25pt] {} &[-25pt] {} &[-25pt] {} &[-25pt] {} &[-25pt] {} &[-25pt] {} &[-25pt] {} &[-25pt] {} &[-25pt] {} &[-25pt] {} &[-25pt] {} &[-25pt] {} &[-25pt] {} &[-25pt] {} \\ 
{} &[-25pt] {} &[-25pt] {} &[-25pt] {} &[-25pt] {} &[-25pt] {} &[-25pt] {} &[-25pt] {} &[-25pt] {} &[-25pt] {} &[-25pt] {} &[-25pt] {} &[-25pt] {} &[-25pt] {} &[-25pt] {} &[-25pt] {} &[-25pt] {} &[-25pt] {} &[-25pt] {} &[-25pt] {} &[-25pt] {} &[-25pt] {} &[-25pt] {} &[-25pt] {} &[-25pt] {} &[-25pt] {} &[-25pt] {} &[-25pt] {} &[-25pt] {} &[-25pt] {} &[-25pt] {} &[-25pt] {} \\ 
{} &[-25pt] {} &[-25pt] {} &[-25pt] {} &[-25pt] {} &[-25pt] {} &[-25pt] {} &[-25pt] {} &[-25pt] {} &[-25pt] {} &[-25pt] {} &[-25pt] {} &[-25pt] {} &[-25pt] {} &[-25pt] {} &[-25pt] {} &[-25pt] -1\ar[-]{uuurrrrrrrr}\ar[-]{uuullllllll} &[-25pt] {} &[-25pt] {} &[-25pt] {} &[-25pt] {} &[-25pt] {} &[-25pt] {} &[-25pt] {} &[-25pt] {} &[-25pt] {} &[-25pt] {} &[-25pt] {} &[-25pt] {} &[-25pt] {} &[-25pt] {} &[-25pt] {} \\ 
{} &[-25pt] {} &[-25pt] {} &[-25pt] {} &[-25pt] {} &[-25pt] {} &[-25pt] {} &[-25pt] {} &[-25pt] {} &[-25pt] {} &[-25pt] {} &[-25pt] {} &[-25pt] {} &[-25pt] {} &[-25pt] {} &[-25pt] {} &[-25pt] {} &[-25pt] {} &[-25pt] {} &[-25pt] {} &[-25pt] {} &[-25pt] {} &[-25pt] {} &[-25pt] {} &[-25pt] {} &[-25pt] {} &[-25pt] {} &[-25pt] {} &[-25pt] {} &[-25pt] {} &[-25pt] {} &[-25pt] {} \\ 
{} &[-25pt] {} &[-25pt] {} &[-25pt] {} &[-25pt] {} &[-25pt] {} &[-25pt] {} &[-25pt] {} &[-25pt] {} &[-25pt] {} &[-25pt] {} &[-25pt] {} &[-25pt] {} &[-25pt] {} &[-25pt] {} &[-25pt] {} &[-25pt] {} &[-25pt] {} &[-25pt] {} &[-25pt] {} &[-25pt] {} &[-25pt] {} &[-25pt] {} &[-25pt] {} &[-25pt] {} &[-25pt] {} &[-25pt] {} &[-25pt] {} &[-25pt] {} &[-25pt] {} &[-25pt] {} &[-25pt] {} \\ 
\cdots \ar[-]{uuurrrrrrrrrrrrrrrr} &[-25pt] {} &[-25pt] {} &[-25pt] {} &[-25pt] {} &[-25pt] {} &[-25pt] {} &[-25pt] {} &[-25pt] {} &[-25pt] {} &[-25pt] {} &[-25pt] {} &[-25pt] {} &[-25pt] {} &[-25pt] {} &[-25pt] {} &[-25pt] {} &[-25pt] {} &[-25pt] {} &[-25pt] {} &[-25pt] {} &[-25pt] {} &[-25pt] {} &[-25pt] {} &[-25pt] {} &[-25pt] {} &[-25pt] {} &[-25pt] {} &[-25pt] {} &[-25pt] {} &[-25pt] {} &[-25pt] {} \\ 
\end{tikzcd}
};
\end{tikzpicture}
\caption{\label{implicittree}
Another way to include negative rationals is to repeat the operation that yields the binary tree rooted at $1$ from the tree rooted at $2$ by subtracting $1$ from each number. Subtracting $1$ from every number in the tree of positive rationals generates a tree that includes the negative numbers greater than $-1$. Repeating this indefinitely generates an infinite binary tree with no root that includes every rational number. This tree has fewer symmetries than the tree in figure \ref{treewithnegativerationals}. Numbers symmetric around $-\nicefrac{1}{2}$ sum to $-1$, but only within its own binary tree. Numbers symmetric around $-1$ are not reciprocals. Symmetries are restricted to the binary tree whose root is the number around which the symmetry manifests, so they don't overlap with other symmetries.}
\end{figure}

This includes the negative rational numbers, which are naturally included in the tree as shown in figures \ref{treewithnegativerationals} and \ref{naturalrepresentationtreewithnegativerationals}. The only pattern that is broken by their inclusion is that there are three edges connecting upward from the number $0$ at the root of the tree. The numbers $1$, $-1$ and $-{}^1\!\!/_{\!\!2}$ can all be reached from $0$ by traversing a single edge.

Zero is the only number that can have three edges from it because $[0;0]$, which indicates the number $-{}^1\!\!/_{\!\!2}$, is the only natural representation that contains two consecutive zeros. The sign of an integer appearing in a natural representation encodes the information about whether the $f\leftrightarrow1-f$ operation needs to be applied before taking the corresponding reciprocal. 

There are only two positions in the list for which this information isn't necessary. Specifically, it isn't needed for the first integer in the list, which is the integer part of the rational number and is not involved in a reciprocal, and it isn't needed for the last integer in a list with two or more integers in it, if that integer is zero, indicating $f={}^1\!\!/_{\!\!2}=1-f$. 

So $[0]$ is the only case when a zero can be added after a list of integers that already ends in zero. Every other node in the tree falls into one of two categories, each of which permits only two of the three moves available from $[0]$:

\begin{itemize}
    \item
    The final integer in the list is negative or positive, and the two available moves are ``go forward'' to the next positive or negative integer and ``branch to a new reciprocal'', which adds a $0$ to the end of the list.
    \item
    The final integer is zero, and the two available moves are ``go in the positive direction'' and ``go in the negative direction'', which change the final $0$ to $1$ and $-1$ respectively.
\end{itemize}

There are new, overlapping symmetries present and visible in figure \ref{treewithnegativerationals} in the full tree including the negative numbers, which aren't present in the purely positive part of the tree. These symmetries are naturally present; no choice of ours introduced them. This specific way of connecting the negative numbers is dictated by the natural representations of the rational numbers. The systematic arrangement of the natural representations is visible in figure \ref{naturalrepresentationtreewithnegativerationals}, which makes it clear that the same rules apply to positive and negative rational numbers.

This structure can be contrasted with the tree shown in figure \ref{simplertree}, which shows how one might expect an extension of the Stern-Brocot tree to negative numbers to appear. While it appears at first to be quite symmetric, each symmetry within it is localized to the specific binary tree with the number around which the symmetry manifests at the root. Different symmetries don't overlap, except in the trivial case when one tree is entirely within another, because they don't extend outside of their local binary trees.

Symmetries that don't extend outside of their own binary tree, and don't overlap with other symmetries, are, in a sense, trivial, because the constraints they impose can always be satisfied, regardless of the objects that are arranged within the tree and how the symmetries relate them. 

For example, if we are building a tree with objects at each vertex, whose root is $x$, and we require any objects, $a$ and $b$, positioned symmetrically around $x$ to satisfy $a=f(b)$ for some function, $f$, then we can build the tree on the left side of $x$ first, and then use $f$ to determine what object should be at each position in the tree on the right side of $x$. If we choose $x_l$ to be the object at the vertex reached by moving left from $x$, then the binary tree starting from $x_l$ is not constrained at all by the symmetry around $x$. We can therefore impose any symmetry we want on that tree, choosing a new, arbitrary function, $g$, which determines the objects on the right side given the objects on the left, without any consideration of the symmetries already imposed.

However, if a tree has multiple overlapping symmetries, each of which extends beyond the binary tree whose root is the line of reflection for that symmetry, then the functions that implement the symmetries, $f$, $g$ and so on, cannot be arbitrary. They must obey an algebra that makes them compatible with one another, and if the functions are simple, natural operations on the objects at the vertices of the tree, then the algebra they satisfy reveals information about the structures of those objects.

In our case, with the symmetries $f(x)=-1-x$ and $g(x)=\nicefrac{1}{x}$, where $f$ reflects a number around $-\nicefrac{1}{2}$ and $g$ reflects positive numbers around $1$ and negative numbers around $-1$, the functions need to satisfy $fgf=gfg$, or, equivalently, $fgfgfg=1$, in order for the two numbers symmetric around $-\nicefrac{1}{2}$ that are obtained from two positive reciprocals, $x$ and $\nicefrac{1}{x}$, by reflecting them first around $-\nicefrac{1}{2}$ and then around $-1$, to satisfy the symmetry around $-\nicefrac{1}{2}$.

So the symmetries of the tree specified by the constituent structures and natural representations of rational numbers specify an algebra, $f^2=g^2=(fg)^3=1$, that encodes certain properties of numbers and corresponding operations on them. This is the algebra of the group, $D_3$, of symmetries of an equilateral triangle, with $f$, $g$, and $fgf=gfg$ corresponding to reflections through the lines that bisect the angles of the triangle, and $fg$ corresponding to rotation through 120 degrees.

One might observe that the tree dictated by constituent structure is not a binary tree, but has a consistent feature that binary trees, as commonly defined, lack: Every vertex in it, including the root at zero, has exactly three edges connected to it. If we ignore constituent structure and arrange these three edges at zero in the same way as every other vertex's edges are arranged, we get the tree shown in figure \ref{implicittree}, which is an infinite binary tree with no root.

This tree has a more uniform structure than the tree shown in figure \ref{treewithnegativerationals}, but the arrangement of the numbers within it has fewer symmetries, none of which overlap.

Within the tree shown in figure \ref{naturalrepresentationtreewithnegativerationals}, the sequences of integers can be seen to exhibit a consistent order in their horizontal and vertical arrangement.

Considering the bottom of the tree containing $[0]$ to have a vertical position or height of 1, the vertical position of a sequence of integers is equal to the length of the sequence plus the sum of the absolute values of its entries. So $-{}^4\!\!/_{\!\!7}=[0;-1,1]$ appears at the $5^{\rm th}$ level in the list, because $5=3+|0|+|-1|+|1|$.

The width of the tree at that height, $h=5$, is given by $3\times2^{h-2}=24$, since the total tree consists of three binary trees originating from level 2. The natural representations of each of those rational numbers can be determined from the rule that the length of each natural representation plus the sum of the absolute values of the integers in it is equal to $h$.

As figure \ref{naturalrepresentationtreewithnegativerationals} shows, the $3\times2^{h-2}$ natural representations at a given height in the tree appear in lexicographical order from left to right. So $[-4]<[-3;0]$ because $-4<-3$, and $[0;2,0] < [0;3]$ because $2<3$.

This allows us to specify that, for a given sequence, its horizontal position in the tree is equal to its position within the sorted list of sequences at that height.

Equation \ref{sequenceorder} below gives a formal specification of the order of sequences that define natural representations of rational numbers. Sequences at the same height in the tree will always differ at some position, but the order extends to the entire tree, and to all natural representations, when it is also specified that sequences of two or more integers that end in negative numbers are less than sequences that continue them, and all other sequences are greater than sequences that continue them:
\begin{align}
\label{sequenceorder}
b_1>c_1\  &\Longrightarrow\  &[a_0;a_1, \cdots, a_n, b_{1}, \cdots, b_m] &> [a_0;a_1, \cdots, a_n, c_{1},\cdots, c_k] \nonumber\\
a_n>0\ \text{\sc\small\ or } \ n=0 \ &\Longrightarrow\ &  [a_0;a_1, \cdots, a_n]& > [a_0;a_1, \cdots, a_n, a_{n+1}, \cdots, a_m] \\
a_n<0\ \text{\sc\small\ and } \ n>0 \ &\Longrightarrow\ &  [a_0;a_1, \cdots, a_n] &< [a_0;a_1, \cdots, a_n, a_{n+1}, \cdots, a_m]
\nonumber
\end{align}

So natural representations have the same order, as sequences of integers, as the rational numbers they represent. If $s(f)$ is the sequence of integers for a given fraction, $f$, then 
$$
s(f_1)\leq s(f_2)\Longleftrightarrow f_1\leq f_2
$$
when sequences are ordered according to equation \ref{sequenceorder}. 

This is not true of the standard continued fraction representation of rational numbers, for which ${}^3\!\!/_{\!\!8}$ and ${}^7\!\!/_{\!\!8}$ have the sequences $[0; 2, 1, 2]$ and $[0; 1, 7]$, and these sequences are on the same side of $[0; 1, 1, 1, 2]$, which represents ${}^5\!\!/_{\!\!8}$, in lexicographic order or any simple ordering of sequences based on their entries, although it is possible to define an ordering on standard continued fraction representations that coincides with the order of the rational numbers they represent, by flipping the sign of every second entry and then taking the lexicographic order.

The standard continued fraction representation is also not unique, but can be made unique by requiring the final denominator to be greater than or equal to 2. Like the sign-flipping procedure that makes the ordering of continued fraction representations match the order of the rationals, and which suggests subtraction rather than addition of the next reciprocal, requiring the final denominator to be greater than or equal to 2 brings the standard continued fraction representation closer to the natural representation.

These are hints from the mathematics of continued fractions that a more fundamental representation of a rational number exists, but it was an examination of the constituent structure of the simplest sets that led us to it.

\raggedbottom

\pagebreak

\appendix

\section{Natural Representation and Continued Fraction Implementations}

\label{codeappendix}

\mbox{This code is available on github at:}\linebreak
\href{https://github.com/roflanagan/natural-representation}{\texttt{https://github.com/roflanagan/natural-representation}}

\subsection{Implementation in Python}

\vspace{3mm}

{\scriptsize
\begin{Verbatim}[commandchars=\\\{\}]
\PY{k+kn}{from} \PY{n+nn}{fractions} \PY{k+kn}{import} \PY{n}{Fraction}
\PY{k+kn}{from} \PY{n+nn}{math} \PY{k+kn}{import} \PY{k+kp}{floor}
\PY{k+kn}{from} \PY{n+nn}{numpy} \PY{k+kn}{import} \PY{k+kp}{array}

\PY{k}{def} \PY{n+nf}{natural\PYZus{}representation}\PY{p}{(}\PY{n}{f}\PY{p}{)}\PY{p}{:}
    \PY{l+s+sd}{\PYZdq{}\PYZdq{}\PYZdq{}}
\PY{l+s+sd}{    Compute the sequence of integers in the natural representation of the fraction f}
\PY{l+s+sd}{    \PYZdq{}\PYZdq{}\PYZdq{}}
    \PY{n}{integer\PYZus{}part} \PY{o}{=} \PY{n+nb}{int}\PY{p}{(}\PY{k+kp}{floor}\PY{p}{(}\PY{n}{f}\PY{p}{)}\PY{p}{)}
    \PY{k}{if} \PY{n}{f} \PY{o}{==} \PY{n}{integer\PYZus{}part}\PY{p}{:}
        \PY{k}{return} \PY{p}{[}\PY{n}{integer\PYZus{}part}\PY{p}{]}
    \PY{n}{fractional\PYZus{}part} \PY{o}{=} \PY{n}{f} \PY{o}{\PYZhy{}} \PY{n}{integer\PYZus{}part}
    \PY{k}{if} \PY{n}{fractional\PYZus{}part} \PY{o}{\PYZgt{}}\PY{o}{=} \PY{l+m+mf}{0.5}\PY{p}{:} 
        \PY{k}{return} \PY{p}{[}\PY{n}{integer\PYZus{}part} \PY{o}{+} \PY{l+m+mi}{1}\PY{p}{]} \PY{o}{+} \PY{n}{natural\PYZus{}representation}\PY{p}{(}\PY{l+m+mi}{1} \PY{o}{/} \PY{p}{(}\PY{l+m+mi}{1}\PY{o}{\PYZhy{}}\PY{n}{fractional\PYZus{}part}\PY{p}{)} \PY{o}{\PYZhy{}} \PY{l+m+mi}{2}\PY{p}{)}
    \PY{k}{else}\PY{p}{:}
        \PY{n}{rest} \PY{o}{=} \PY{k+kp}{array}\PY{p}{(}\PY{n}{natural\PYZus{}representation}\PY{p}{(}\PY{l+m+mi}{1} \PY{o}{/} \PY{n}{fractional\PYZus{}part} \PY{o}{\PYZhy{}} \PY{l+m+mi}{2}\PY{p}{)}\PY{p}{)}
        \PY{k}{return} \PY{p}{[}\PY{n}{integer\PYZus{}part} \PY{o}{+} \PY{l+m+mi}{1}\PY{p}{]} \PY{o}{+} \PY{n+nb}{list}\PY{p}{(}\PY{o}{\PYZhy{}}\PY{n}{rest}\PY{p}{)}

\PY{k}{def} \PY{n+nf}{evaluate\PYZus{}natural\PYZus{}representation}\PY{p}{(}\PY{n}{sequence}\PY{p}{)}\PY{p}{:}
    \PY{l+s+sd}{\PYZdq{}\PYZdq{}\PYZdq{}}
\PY{l+s+sd}{    Compute the fraction f from the sequence of integers in its natural representation}
\PY{l+s+sd}{    \PYZdq{}\PYZdq{}\PYZdq{}}
    \PY{n}{rest\PYZus{}of\PYZus{}sequence} \PY{o}{=} \PY{k+kp}{array}\PY{p}{(}\PY{n}{sequence}\PY{p}{[}\PY{l+m+mi}{1}\PY{p}{:}\PY{p}{]}\PY{p}{)}
    \PY{k}{if} \PY{n+nb}{len}\PY{p}{(}\PY{n}{rest\PYZus{}of\PYZus{}sequence}\PY{p}{)} \PY{o}{==} \PY{l+m+mi}{0}\PY{p}{:}
        \PY{k}{return} \PY{n}{sequence}\PY{p}{[}\PY{l+m+mi}{0}\PY{p}{]}
    \PY{k}{if} \PY{n}{rest\PYZus{}of\PYZus{}sequence}\PY{p}{[}\PY{l+m+mi}{0}\PY{p}{]} \PY{o}{\PYZgt{}}\PY{o}{=} \PY{l+m+mi}{0}\PY{p}{:}
        \PY{n}{rest} \PY{o}{=} \PY{n}{evaluate\PYZus{}natural\PYZus{}representation}\PY{p}{(}\PY{n}{rest\PYZus{}of\PYZus{}sequence}\PY{p}{)}
        \PY{k}{return} \PY{n}{sequence}\PY{p}{[}\PY{l+m+mi}{0}\PY{p}{]} \PY{o}{\PYZhy{}} \PY{n}{Fraction}\PY{p}{(}\PY{l+m+mi}{1}\PY{p}{,} \PY{l+m+mi}{2} \PY{o}{+} \PY{n}{rest}\PY{p}{)}
    \PY{k}{else}\PY{p}{:}
        \PY{n}{rest} \PY{o}{=} \PY{n}{evaluate\PYZus{}natural\PYZus{}representation}\PY{p}{(}\PY{o}{\PYZhy{}}\PY{n}{rest\PYZus{}of\PYZus{}sequence}\PY{p}{)}
        \PY{k}{return} \PY{n}{sequence}\PY{p}{[}\PY{l+m+mi}{0}\PY{p}{]} \PY{o}{\PYZhy{}} \PY{l+m+mi}{1} \PY{o}{+} \PY{n}{Fraction}\PY{p}{(}\PY{l+m+mi}{1}\PY{p}{,} \PY{l+m+mi}{2} \PY{o}{+} \PY{n}{rest}\PY{p}{)}

\PY{k}{def} \PY{n+nf}{show\PYZus{}examples}\PY{p}{(}\PY{n}{max\PYZus{}numerator}\PY{p}{,} \PY{n}{max\PYZus{}denominator}\PY{p}{)}\PY{p}{:}
    \PY{l+s+sd}{\PYZdq{}\PYZdq{}\PYZdq{}}
\PY{l+s+sd}{    Encode and decode all positive and negative rational numbers whose numerators}
\PY{l+s+sd}{    and denominators don\PYZsq{}t exceed the specified limits}
\PY{l+s+sd}{    \PYZdq{}\PYZdq{}\PYZdq{}}
    \PY{k+kn}{from} \PY{n+nn}{fractions} \PY{k+kn}{import} \PY{n}{gcd} \PY{k}{as} \PY{n}{greatest\PYZus{}common\PYZus{}divisor}
    \PY{k}{for} \PY{n}{n} \PY{o+ow}{in} \PY{n+nb}{range}\PY{p}{(}\PY{o}{\PYZhy{}}\PY{n}{max\PYZus{}numerator}\PY{p}{,} \PY{n}{max\PYZus{}numerator} \PY{o}{+} \PY{l+m+mi}{1}\PY{p}{)}\PY{p}{:}
        \PY{k}{for} \PY{n}{d} \PY{o+ow}{in} \PY{n+nb}{range}\PY{p}{(}\PY{l+m+mi}{1}\PY{p}{,} \PY{n}{max\PYZus{}denominator} \PY{o}{+} \PY{l+m+mi}{1}\PY{p}{)}\PY{p}{:}
            \PY{k}{if} \PY{n}{greatest\PYZus{}common\PYZus{}divisor}\PY{p}{(}\PY{n}{n}\PY{p}{,} \PY{n}{d}\PY{p}{)} \PY{o}{==} \PY{l+m+mi}{1}\PY{p}{:}
                \PY{n}{fraction} \PY{o}{=} \PY{n}{Fraction}\PY{p}{(}\PY{n}{n}\PY{p}{,} \PY{n}{d}\PY{p}{)}
                \PY{n}{sequence} \PY{o}{=} \PY{n}{natural\PYZus{}representation}\PY{p}{(}\PY{n}{fraction}\PY{p}{)}
                \PY{n}{recovered\PYZus{}fraction} \PY{o}{=} \PY{n}{evaluate\PYZus{}natural\PYZus{}representation}\PY{p}{(}\PY{n}{sequence}\PY{p}{)}
                \PY{k}{print}\PY{p}{(}\PY{n}{fraction}\PY{p}{,} \PY{l+s+s2}{\PYZdq{}}\PY{l+s+se}{\PYZbs{}t}\PY{l+s+s2}{\PYZhy{}\PYZhy{}\PYZgt{}}\PY{l+s+s2}{\PYZdq{}}\PY{p}{,} \PY{n+nb}{str}\PY{p}{(}\PY{n}{sequence}\PY{p}{)}\PY{o}{.}\PY{n}{ljust}\PY{p}{(}\PY{l+m+mi}{12}\PY{p}{)}\PY{p}{,} \PY{l+s+s2}{\PYZdq{}}\PY{l+s+s2}{\PYZhy{}\PYZhy{}\PYZgt{}}\PY{l+s+s2}{\PYZdq{}}\PY{p}{,} \PY{n}{recovered\PYZus{}fraction}\PY{p}{)}


\PY{k}{if} \PY{n+nv+vm}{\PYZus{}\PYZus{}name\PYZus{}\PYZus{}} \PY{o}{==} \PY{l+s+s2}{\PYZdq{}}\PY{l+s+s2}{\PYZus{}\PYZus{}main\PYZus{}\PYZus{}}\PY{l+s+s2}{\PYZdq{}}\PY{p}{:}
    \PY{n}{show\PYZus{}examples}\PY{p}{(}\PY{l+m+mi}{9}\PY{p}{,} \PY{l+m+mi}{9}\PY{p}{)}
\end{Verbatim}

}

\pagebreak

\subsection{Implementation in C}

\vspace{3mm}

{\scriptsize
\begin{Verbatim}[commandchars=\\\{\}]
\PY{c+cp}{\PYZsh{}}\PY{c+cp}{include} \PY{c+cpf}{\PYZlt{}stdio.h\PYZgt{}}
\PY{c+cp}{\PYZsh{}}\PY{c+cp}{include} \PY{c+cpf}{\PYZlt{}math.h\PYZgt{}}
\PY{c+cp}{\PYZsh{}}\PY{c+cp}{include} \PY{c+cpf}{\PYZlt{}string.h\PYZgt{}}
\PY{c+cp}{\PYZsh{}}\PY{c+cp}{include} \PY{c+cpf}{\PYZlt{}sys/time.h\PYZgt{}}
\PY{c+cp}{\PYZsh{}}\PY{c+cp}{include} \PY{c+cpf}{\PYZlt{}stdlib.h\PYZgt{}}

\PY{k+kt}{void} \PY{n+nf}{continued\PYZus{}fraction}\PY{p}{(}\PY{k+kt}{unsigned} \PY{k+kt}{long} \PY{n}{numerator}\PY{p}{,} \PY{k+kt}{unsigned} \PY{k+kt}{long} \PY{n}{denominator}\PY{p}{,} \PY{k+kt}{int} \PY{o}{*}\PY{n}{sequence}\PY{p}{,} \PY{k+kt}{int} \PY{o}{*}\PY{n}{i}\PY{p}{)}
\PY{p}{\PYZob{}}
    \PY{k+kt}{unsigned} \PY{k+kt}{long} \PY{n}{new\PYZus{}numerator}\PY{p}{;}
    \PY{o}{*}\PY{n}{i} \PY{o}{=} \PY{l+m+mi}{0}\PY{p}{;}
    
    \PY{k}{while} \PY{p}{(}\PY{n}{denominator} \PY{o}{\PYZgt{}} \PY{l+m+mi}{1}\PY{p}{)}
    \PY{p}{\PYZob{}}
        \PY{n}{sequence}\PY{p}{[}\PY{o}{*}\PY{n}{i}\PY{p}{]} \PY{o}{=} \PY{n}{numerator} \PY{o}{/} \PY{n}{denominator}\PY{p}{;}
        \PY{n}{new\PYZus{}numerator} \PY{o}{=} \PY{n}{numerator} \PY{o}{\PYZhy{}} \PY{n}{sequence}\PY{p}{[}\PY{o}{*}\PY{n}{i}\PY{p}{]} \PY{o}{*} \PY{n}{denominator}\PY{p}{;}
        \PY{n}{numerator} \PY{o}{=} \PY{n}{denominator}\PY{p}{;}
        \PY{n}{denominator} \PY{o}{=} \PY{n}{new\PYZus{}numerator}\PY{p}{;}
        \PY{o}{*}\PY{n}{i} \PY{o}{=} \PY{o}{*}\PY{n}{i} \PY{o}{+} \PY{l+m+mi}{1}\PY{p}{;}
    \PY{p}{\PYZcb{}}
    \PY{n}{sequence}\PY{p}{[}\PY{o}{*}\PY{n}{i}\PY{p}{]} \PY{o}{=} \PY{n}{numerator}\PY{p}{;}
    \PY{o}{*}\PY{n}{i} \PY{o}{=} \PY{o}{*}\PY{n}{i} \PY{o}{+} \PY{n}{denominator}\PY{p}{;}
\PY{p}{\PYZcb{}}

\PY{k+kt}{void} \PY{n+nf}{evaluate\PYZus{}continued\PYZus{}fraction}\PY{p}{(}\PY{k+kt}{unsigned} \PY{k+kt}{long} \PY{o}{*}\PY{n}{numerator}\PY{p}{,} \PY{k+kt}{unsigned} \PY{k+kt}{long} \PY{o}{*}\PY{n}{denominator}\PY{p}{,} 
                                 \PY{k+kt}{int} \PY{o}{*}\PY{n}{sequence}\PY{p}{,} \PY{k+kt}{int} \PY{o}{*}\PY{n}{sequence\PYZus{}size}\PY{p}{)}
\PY{p}{\PYZob{}}
    \PY{k+kt}{unsigned} \PY{k+kt}{long} \PY{n}{old\PYZus{}numerator}\PY{p}{;}
    \PY{k+kt}{int} \PY{n}{i} \PY{o}{=} \PY{o}{*}\PY{n}{sequence\PYZus{}size} \PY{o}{\PYZhy{}} \PY{l+m+mi}{1}\PY{p}{;}
    \PY{o}{*}\PY{n}{numerator} \PY{o}{=} \PY{n}{sequence}\PY{p}{[}\PY{n}{i}\PY{p}{]}\PY{p}{;}
    \PY{o}{*}\PY{n}{denominator} \PY{o}{=} \PY{l+m+mi}{1}\PY{p}{;}
    
    \PY{k}{while} \PY{p}{(}\PY{n}{i} \PY{o}{\PYZgt{}} \PY{l+m+mi}{0}\PY{p}{)}
    \PY{p}{\PYZob{}}
        \PY{n}{i} \PY{o}{\PYZhy{}}\PY{o}{=} \PY{l+m+mi}{1}\PY{p}{;}
        \PY{n}{old\PYZus{}numerator} \PY{o}{=} \PY{o}{*}\PY{n}{numerator}\PY{p}{;}
        \PY{o}{*}\PY{n}{numerator} \PY{o}{=} \PY{o}{*}\PY{n}{denominator} \PY{o}{+} \PY{n}{sequence}\PY{p}{[}\PY{n}{i}\PY{p}{]} \PY{o}{*} \PY{o}{*}\PY{n}{numerator}\PY{p}{;}
        \PY{o}{*}\PY{n}{denominator} \PY{o}{=} \PY{n}{old\PYZus{}numerator}\PY{p}{;}
    \PY{p}{\PYZcb{}}
\PY{p}{\PYZcb{}}

\PY{k+kt}{void} \PY{n+nf}{natural\PYZus{}representation}\PY{p}{(}\PY{k+kt}{unsigned} \PY{k+kt}{long} \PY{n}{numerator}\PY{p}{,} \PY{k+kt}{unsigned} \PY{k+kt}{long} \PY{n}{denominator}\PY{p}{,} \PY{k+kt}{int} \PY{o}{*}\PY{n}{sequence}\PY{p}{,} \PY{k+kt}{int} \PY{o}{*}\PY{n}{i}\PY{p}{)}
\PY{p}{\PYZob{}}
    \PY{k+kt}{int} \PY{n}{sign} \PY{o}{=} \PY{l+m+mi}{1}\PY{p}{;}
    \PY{k+kt}{unsigned} \PY{k+kt}{long} \PY{n}{integer\PYZus{}part}\PY{p}{,} \PY{n}{fractional\PYZus{}part}\PY{p}{,} \PY{n}{old\PYZus{}denominator}\PY{p}{;}
    \PY{o}{*}\PY{n}{i} \PY{o}{=} \PY{l+m+mi}{0}\PY{p}{;}

    \PY{k}{while} \PY{p}{(}\PY{l+m+mi}{1}\PY{p}{)}
    \PY{p}{\PYZob{}}
        \PY{n}{integer\PYZus{}part} \PY{o}{=} \PY{n}{numerator} \PY{o}{/} \PY{n}{denominator}\PY{p}{;}
        \PY{n}{fractional\PYZus{}part} \PY{o}{=} \PY{n}{numerator} \PY{o}{\PYZhy{}} \PY{n}{integer\PYZus{}part} \PY{o}{*} \PY{n}{denominator}\PY{p}{;}
        \PY{n}{sequence}\PY{p}{[}\PY{o}{*}\PY{n}{i}\PY{p}{]} \PY{o}{=} \PY{n}{sign} \PY{o}{*} \PY{p}{(}\PY{n}{integer\PYZus{}part} \PY{o}{+} \PY{p}{(}\PY{n}{fractional\PYZus{}part} \PY{o}{\PYZgt{}} \PY{l+m+mi}{0}\PY{p}{)}\PY{p}{)}\PY{p}{;}
        \PY{o}{*}\PY{n}{i} \PY{o}{=} \PY{o}{*}\PY{n}{i} \PY{o}{+} \PY{l+m+mi}{1}\PY{p}{;}
        \PY{k}{if} \PY{p}{(}\PY{l+m+mi}{2} \PY{o}{*} \PY{n}{fractional\PYZus{}part} \PY{o}{\PYZgt{}}\PY{o}{=} \PY{n}{denominator}\PY{p}{)}
        \PY{p}{\PYZob{}}
            \PY{n}{denominator} \PY{o}{=} \PY{n}{denominator} \PY{o}{\PYZhy{}} \PY{n}{fractional\PYZus{}part}\PY{p}{;}
            \PY{n}{numerator} \PY{o}{=} \PY{n}{fractional\PYZus{}part} \PY{o}{\PYZhy{}} \PY{n}{denominator}\PY{p}{;}
        \PY{p}{\PYZcb{}}
        \PY{k}{else}
        \PY{p}{\PYZob{}}
            \PY{k}{if} \PY{p}{(}\PY{n}{fractional\PYZus{}part} \PY{o}{=}\PY{o}{=} \PY{l+m+mi}{0}\PY{p}{)}
                \PY{k}{return}\PY{p}{;}
            \PY{n}{old\PYZus{}denominator} \PY{o}{=} \PY{n}{denominator}\PY{p}{;}
            \PY{n}{denominator} \PY{o}{=} \PY{n}{fractional\PYZus{}part}\PY{p}{;}
            \PY{n}{numerator} \PY{o}{=} \PY{n}{old\PYZus{}denominator} \PY{o}{\PYZhy{}} \PY{l+m+mi}{2} \PY{o}{*} \PY{n}{fractional\PYZus{}part}\PY{p}{;}
            \PY{n}{sign} \PY{o}{*}\PY{o}{=} \PY{o}{\PYZhy{}}\PY{l+m+mi}{1}\PY{p}{;}
        \PY{p}{\PYZcb{}}
        
    \PY{p}{\PYZcb{}}
\PY{p}{\PYZcb{}}

\PY{k+kt}{void} \PY{n+nf}{evaluate\PYZus{}natural\PYZus{}representation}\PY{p}{(}\PY{k+kt}{unsigned} \PY{k+kt}{long} \PY{o}{*}\PY{n}{numerator}\PY{p}{,} \PY{k+kt}{unsigned} \PY{k+kt}{long} \PY{o}{*}\PY{n}{denominator}\PY{p}{,} 
                                     \PY{k+kt}{int} \PY{o}{*}\PY{n}{sequence}\PY{p}{,} \PY{k+kt}{int} \PY{o}{*}\PY{n}{sequence\PYZus{}size}\PY{p}{)}
\PY{p}{\PYZob{}}
    \PY{k+kt}{unsigned} \PY{k+kt}{long} \PY{n}{new\PYZus{}denominator}\PY{p}{,} \PY{n}{old\PYZus{}numerator}\PY{p}{;}
    \PY{k+kt}{int} \PY{n}{i} \PY{o}{=} \PY{o}{*}\PY{n}{sequence\PYZus{}size} \PY{o}{\PYZhy{}} \PY{l+m+mi}{1}\PY{p}{,} \PY{n}{previous\PYZus{}sign} \PY{o}{=} \PY{p}{(}\PY{n}{sequence}\PY{p}{[}\PY{n}{i}\PY{p}{]} \PY{o}{\PYZgt{}}\PY{o}{=} \PY{l+m+mi}{0}\PY{p}{)}\PY{p}{;}
    \PY{o}{*}\PY{n}{numerator} \PY{o}{=} \PY{n}{abs}\PY{p}{(}\PY{n}{sequence}\PY{p}{[}\PY{n}{i}\PY{p}{]}\PY{p}{)}\PY{p}{;}
    \PY{o}{*}\PY{n}{denominator} \PY{o}{=} \PY{l+m+mi}{1}\PY{p}{;}
 
    \PY{k}{while} \PY{p}{(}\PY{n}{i} \PY{o}{\PYZgt{}} \PY{l+m+mi}{0}\PY{p}{)}
    \PY{p}{\PYZob{}}
        \PY{n}{i} \PY{o}{\PYZhy{}}\PY{o}{=} \PY{l+m+mi}{1}\PY{p}{;}
        \PY{n}{old\PYZus{}numerator} \PY{o}{=} \PY{o}{*}\PY{n}{numerator}\PY{p}{;}
        \PY{n}{new\PYZus{}denominator} \PY{o}{=} \PY{n}{old\PYZus{}numerator} \PY{o}{+} \PY{l+m+mi}{2} \PY{o}{*} \PY{o}{*}\PY{n}{denominator}\PY{p}{;}
        \PY{o}{*}\PY{n}{numerator} \PY{o}{=} \PY{n}{abs}\PY{p}{(}\PY{n}{sequence}\PY{p}{[}\PY{n}{i}\PY{p}{]}\PY{p}{)} \PY{o}{*} \PY{n}{new\PYZus{}denominator} \PY{o}{\PYZhy{}} \PY{o}{*}\PY{n}{denominator}\PY{p}{;}
        
        \PY{k}{if} \PY{p}{(}\PY{p}{(}\PY{n}{sequence}\PY{p}{[}\PY{n}{i}\PY{p}{]} \PY{o}{\PYZlt{}} \PY{l+m+mi}{0}\PY{p}{)} \PY{o}{=}\PY{o}{=} \PY{n}{previous\PYZus{}sign}\PY{p}{)}
        \PY{p}{\PYZob{}}
            \PY{o}{*}\PY{n}{numerator} \PY{o}{\PYZhy{}}\PY{o}{=} \PY{n}{old\PYZus{}numerator}\PY{p}{;}
            \PY{n}{previous\PYZus{}sign} \PY{o}{=} \PY{l+m+mi}{1} \PY{o}{\PYZhy{}} \PY{n}{previous\PYZus{}sign}\PY{p}{;}
        \PY{p}{\PYZcb{}}
        \PY{o}{*}\PY{n}{denominator} \PY{o}{=} \PY{n}{new\PYZus{}denominator}\PY{p}{;}
    \PY{p}{\PYZcb{}}
\PY{p}{\PYZcb{}}

\PY{k+kt}{unsigned} \PY{k+kt}{long} \PY{n+nf}{current\PYZus{}time}\PY{p}{(}\PY{p}{)}
\PY{p}{\PYZob{}}
    \PY{k}{struct} \PY{n}{timeval} \PY{n}{tv}\PY{p}{;}
    \PY{n}{gettimeofday}\PY{p}{(}\PY{o}{\PYZam{}}\PY{n}{tv}\PY{p}{,}\PY{n+nb}{NULL}\PY{p}{)}\PY{p}{;}
    \PY{k}{return} \PY{l+m+mi}{1000000} \PY{o}{*} \PY{n}{tv}\PY{p}{.}\PY{n}{tv\PYZus{}sec} \PY{o}{+} \PY{n}{tv}\PY{p}{.}\PY{n}{tv\PYZus{}usec}\PY{p}{;}
\PY{p}{\PYZcb{}}

\PY{k+kt}{void} \PY{n+nf}{repeatedly\PYZus{}construct\PYZus{}continued\PYZus{}fraction}\PY{p}{(}\PY{k+kt}{unsigned} \PY{k+kt}{long} \PY{n}{numerator}\PY{p}{,} 
                                             \PY{k+kt}{unsigned} \PY{k+kt}{long} \PY{n}{denominator}\PY{p}{,} 
                                             \PY{k+kt}{int} \PY{o}{*}\PY{n}{sequence}\PY{p}{,} \PY{k+kt}{int} \PY{o}{*}\PY{n}{sequence\PYZus{}size}\PY{p}{,} \PY{k+kt}{int} \PY{n}{iterations}\PY{p}{)}
\PY{p}{\PYZob{}}
    \PY{k+kt}{long} \PY{n}{start}\PY{p}{,} \PY{n}{i}\PY{p}{;}
    \PY{n}{printf}\PY{p}{(}\PY{l+s}{\PYZdq{}}\PY{l+s}{continued fraction:}\PY{l+s+se}{\PYZbs{}n}\PY{l+s}{\PYZdq{}}\PY{p}{)}\PY{p}{;}
    
    \PY{n}{start} \PY{o}{=} \PY{n}{current\PYZus{}time}\PY{p}{(}\PY{p}{)}\PY{p}{;}
    \PY{k}{for} \PY{p}{(}\PY{n}{i} \PY{o}{=} \PY{l+m+mi}{0}\PY{p}{;} \PY{n}{i} \PY{o}{\PYZlt{}} \PY{n}{iterations}\PY{p}{;} \PY{n}{i}\PY{o}{+}\PY{o}{+}\PY{p}{)}
        \PY{n}{continued\PYZus{}fraction}\PY{p}{(}\PY{n}{numerator}\PY{p}{,} \PY{n}{denominator}\PY{p}{,} \PY{n}{sequence}\PY{p}{,} \PY{n}{sequence\PYZus{}size}\PY{p}{)}\PY{p}{;}

    \PY{n}{printf}\PY{p}{(}\PY{l+s}{\PYZdq{}}\PY{l+s}{completed in \PYZpc{}ld microseconds}\PY{l+s+se}{\PYZbs{}n}\PY{l+s}{\PYZdq{}}\PY{p}{,} \PY{n}{current\PYZus{}time}\PY{p}{(}\PY{p}{)}\PY{o}{\PYZhy{}}\PY{n}{start}\PY{p}{)}\PY{p}{;}

    \PY{k}{for} \PY{p}{(}\PY{n}{i}\PY{o}{=}\PY{l+m+mi}{0} \PY{p}{;} \PY{n}{i} \PY{o}{\PYZlt{}} \PY{o}{*}\PY{n}{sequence\PYZus{}size}\PY{p}{;} \PY{n}{i}\PY{o}{+}\PY{o}{+}\PY{p}{)}
        \PY{n}{printf}\PY{p}{(}\PY{l+s}{\PYZdq{}}\PY{l+s}{\PYZpc{}d, }\PY{l+s}{\PYZdq{}}\PY{p}{,} \PY{n}{sequence}\PY{p}{[}\PY{n}{i}\PY{p}{]}\PY{p}{)}\PY{p}{;}
    \PY{n}{printf}\PY{p}{(}\PY{l+s}{\PYZdq{}}\PY{l+s+se}{\PYZbs{}n}\PY{l+s}{\PYZdq{}}\PY{p}{)}\PY{p}{;}
\PY{p}{\PYZcb{}}

\PY{k+kt}{void} \PY{n+nf}{repeatedly\PYZus{}evaluate\PYZus{}continued\PYZus{}fraction}\PY{p}{(}\PY{k+kt}{unsigned} \PY{k+kt}{long} \PY{n}{numerator}\PY{p}{,} 
                                            \PY{k+kt}{unsigned} \PY{k+kt}{long} \PY{n}{denominator}\PY{p}{,} 
                                            \PY{k+kt}{int} \PY{o}{*}\PY{n}{sequence}\PY{p}{,} \PY{k+kt}{int} \PY{o}{*}\PY{n}{sequence\PYZus{}size}\PY{p}{,} \PY{k+kt}{int} \PY{n}{iterations}\PY{p}{)}
\PY{p}{\PYZob{}}
    \PY{k+kt}{long} \PY{n}{start}\PY{p}{,} \PY{n}{i}\PY{p}{;}
    \PY{n}{printf}\PY{p}{(}\PY{l+s}{\PYZdq{}}\PY{l+s+se}{\PYZbs{}n}\PY{l+s+se}{\PYZbs{}n}\PY{l+s}{evaluation of continued fraction:}\PY{l+s+se}{\PYZbs{}n}\PY{l+s}{\PYZdq{}}\PY{p}{)}\PY{p}{;}
    
    \PY{n}{start} \PY{o}{=} \PY{n}{current\PYZus{}time}\PY{p}{(}\PY{p}{)}\PY{p}{;}
    \PY{k}{for} \PY{p}{(}\PY{n}{i} \PY{o}{=} \PY{l+m+mi}{0}\PY{p}{;} \PY{n}{i} \PY{o}{\PYZlt{}} \PY{n}{iterations}\PY{p}{;} \PY{n}{i}\PY{o}{+}\PY{o}{+}\PY{p}{)}
        \PY{n}{evaluate\PYZus{}continued\PYZus{}fraction}\PY{p}{(}\PY{o}{\PYZam{}}\PY{n}{numerator}\PY{p}{,} \PY{o}{\PYZam{}}\PY{n}{denominator}\PY{p}{,} \PY{n}{sequence}\PY{p}{,} \PY{n}{sequence\PYZus{}size}\PY{p}{)}\PY{p}{;}

    \PY{n}{printf}\PY{p}{(}\PY{l+s}{\PYZdq{}}\PY{l+s}{completed in \PYZpc{}ld microseconds}\PY{l+s+se}{\PYZbs{}n}\PY{l+s}{\PYZdq{}}\PY{p}{,} \PY{n}{current\PYZus{}time}\PY{p}{(}\PY{p}{)}\PY{o}{\PYZhy{}}\PY{n}{start}\PY{p}{)}\PY{p}{;}
    \PY{n}{printf}\PY{p}{(}\PY{l+s}{\PYZdq{}}\PY{l+s}{numerator: \PYZpc{}lu   denominator: \PYZpc{}lu}\PY{l+s+se}{\PYZbs{}n}\PY{l+s+se}{\PYZbs{}n}\PY{l+s}{\PYZdq{}}\PY{p}{,} \PY{n}{numerator}\PY{p}{,} \PY{n}{denominator}\PY{p}{)}\PY{p}{;}
\PY{p}{\PYZcb{}}

\PY{k+kt}{void} \PY{n+nf}{repeatedly\PYZus{}construct\PYZus{}natural\PYZus{}representation}\PY{p}{(}\PY{k+kt}{unsigned} \PY{k+kt}{long} \PY{n}{numerator}\PY{p}{,} 
                                                 \PY{k+kt}{unsigned} \PY{k+kt}{long} \PY{n}{denominator}\PY{p}{,} 
                                                 \PY{k+kt}{int} \PY{o}{*}\PY{n}{sequence}\PY{p}{,} \PY{k+kt}{int} \PY{o}{*}\PY{n}{sequence\PYZus{}size}\PY{p}{,} 
                                                 \PY{k+kt}{int} \PY{n}{iterations}\PY{p}{)}
\PY{p}{\PYZob{}}
    \PY{k+kt}{long} \PY{n}{start}\PY{p}{,} \PY{n}{i}\PY{p}{;}

    \PY{n}{printf}\PY{p}{(}\PY{l+s}{\PYZdq{}}\PY{l+s+se}{\PYZbs{}n}\PY{l+s+se}{\PYZbs{}n}\PY{l+s+se}{\PYZbs{}n}\PY{l+s}{natural representation:}\PY{l+s+se}{\PYZbs{}n}\PY{l+s}{\PYZdq{}}\PY{p}{)}\PY{p}{;}
    
    \PY{n}{start} \PY{o}{=} \PY{n}{current\PYZus{}time}\PY{p}{(}\PY{p}{)}\PY{p}{;}
    \PY{k}{for} \PY{p}{(}\PY{n}{i} \PY{o}{=} \PY{l+m+mi}{0}\PY{p}{;} \PY{n}{i} \PY{o}{\PYZlt{}} \PY{n}{iterations}\PY{p}{;} \PY{n}{i}\PY{o}{+}\PY{o}{+}\PY{p}{)}
        \PY{n}{natural\PYZus{}representation}\PY{p}{(}\PY{n}{numerator}\PY{p}{,} \PY{n}{denominator}\PY{p}{,} \PY{n}{sequence}\PY{p}{,} \PY{n}{sequence\PYZus{}size}\PY{p}{)}\PY{p}{;}
    
    \PY{n}{printf}\PY{p}{(}\PY{l+s}{\PYZdq{}}\PY{l+s}{completed in \PYZpc{}ld microseconds}\PY{l+s+se}{\PYZbs{}n}\PY{l+s}{\PYZdq{}}\PY{p}{,} \PY{n}{current\PYZus{}time}\PY{p}{(}\PY{p}{)}\PY{o}{\PYZhy{}}\PY{n}{start}\PY{p}{)}\PY{p}{;}

    \PY{k}{for} \PY{p}{(}\PY{n}{i}\PY{o}{=}\PY{l+m+mi}{0} \PY{p}{;} \PY{n}{i} \PY{o}{\PYZlt{}} \PY{o}{*}\PY{n}{sequence\PYZus{}size}\PY{p}{;} \PY{n}{i}\PY{o}{+}\PY{o}{+}\PY{p}{)}
        \PY{n}{printf}\PY{p}{(}\PY{l+s}{\PYZdq{}}\PY{l+s}{\PYZpc{}d, }\PY{l+s}{\PYZdq{}}\PY{p}{,} \PY{n}{sequence}\PY{p}{[}\PY{n}{i}\PY{p}{]}\PY{p}{)}\PY{p}{;}
    \PY{n}{printf}\PY{p}{(}\PY{l+s}{\PYZdq{}}\PY{l+s+se}{\PYZbs{}n}\PY{l+s+se}{\PYZbs{}n}\PY{l+s}{\PYZdq{}}\PY{p}{)}\PY{p}{;}
\PY{p}{\PYZcb{}}

\PY{k+kt}{void} \PY{n+nf}{repeatedly\PYZus{}evaluate\PYZus{}natural\PYZus{}representation}\PY{p}{(}\PY{k+kt}{unsigned} \PY{k+kt}{long} \PY{n}{numerator}\PY{p}{,} 
                                                \PY{k+kt}{unsigned} \PY{k+kt}{long} \PY{n}{denominator}\PY{p}{,} 
                                                \PY{k+kt}{int} \PY{o}{*}\PY{n}{sequence}\PY{p}{,} \PY{k+kt}{int} \PY{o}{*}\PY{n}{sequence\PYZus{}size}\PY{p}{,}
                                                \PY{k+kt}{int} \PY{n}{iterations}\PY{p}{)}
\PY{p}{\PYZob{}}
    \PY{k+kt}{long} \PY{n}{start}\PY{p}{,} \PY{n}{i}\PY{p}{;}

    \PY{n}{printf}\PY{p}{(}\PY{l+s}{\PYZdq{}}\PY{l+s+se}{\PYZbs{}n}\PY{l+s+se}{\PYZbs{}n}\PY{l+s}{evaluation of natural representation:}\PY{l+s+se}{\PYZbs{}n}\PY{l+s}{\PYZdq{}}\PY{p}{)}\PY{p}{;}
    
    \PY{n}{start} \PY{o}{=} \PY{n}{current\PYZus{}time}\PY{p}{(}\PY{p}{)}\PY{p}{;}
    \PY{k}{for} \PY{p}{(}\PY{n}{i} \PY{o}{=} \PY{l+m+mi}{0}\PY{p}{;} \PY{n}{i} \PY{o}{\PYZlt{}} \PY{n}{iterations}\PY{p}{;} \PY{n}{i}\PY{o}{+}\PY{o}{+}\PY{p}{)}
        \PY{n}{evaluate\PYZus{}natural\PYZus{}representation}\PY{p}{(}\PY{o}{\PYZam{}}\PY{n}{numerator}\PY{p}{,} \PY{o}{\PYZam{}}\PY{n}{denominator}\PY{p}{,} \PY{n}{sequence}\PY{p}{,} \PY{n}{sequence\PYZus{}size}\PY{p}{)}\PY{p}{;}

    \PY{n}{printf}\PY{p}{(}\PY{l+s}{\PYZdq{}}\PY{l+s}{completed in \PYZpc{}ld microseconds}\PY{l+s+se}{\PYZbs{}n}\PY{l+s}{\PYZdq{}}\PY{p}{,} \PY{n}{current\PYZus{}time}\PY{p}{(}\PY{p}{)}\PY{o}{\PYZhy{}}\PY{n}{start}\PY{p}{)}\PY{p}{;}
    \PY{n}{printf}\PY{p}{(}\PY{l+s}{\PYZdq{}}\PY{l+s}{numerator: \PYZpc{}lu   denominator: \PYZpc{}lu}\PY{l+s+se}{\PYZbs{}n}\PY{l+s+se}{\PYZbs{}n}\PY{l+s}{\PYZdq{}}\PY{p}{,} \PY{n}{numerator}\PY{p}{,} \PY{n}{denominator}\PY{p}{)}\PY{p}{;}
\PY{p}{\PYZcb{}}

\PY{k+kt}{int} \PY{n+nf}{main}\PY{p}{(}\PY{k+kt}{int} \PY{n}{argc}\PY{p}{,} \PY{k+kt}{char} \PY{o}{*}\PY{o}{*}\PY{n}{argv}\PY{p}{)}
\PY{p}{\PYZob{}}
    \PY{k+kt}{unsigned} \PY{k+kt}{long} \PY{n}{numerator}\PY{p}{,} \PY{n}{denominator}\PY{p}{;}
    \PY{k+kt}{int} \PY{n}{sequence}\PY{p}{[}\PY{n}{BUFSIZ}\PY{p}{]}\PY{p}{,} \PY{n}{sequence\PYZus{}size}\PY{p}{,} \PY{n}{iterations} \PY{o}{=} \PY{l+m+mi}{10000000}\PY{p}{;}

    \PY{k}{if} \PY{p}{(}\PY{n}{argc} \PY{o}{\PYZlt{}} \PY{l+m+mi}{3}\PY{p}{)}
    \PY{p}{\PYZob{}}
        \PY{n}{printf}\PY{p}{(}\PY{l+s}{\PYZdq{}}\PY{l+s}{Usage: \PYZpc{}s numerator denominator}\PY{l+s+se}{\PYZbs{}n}\PY{l+s}{\PYZdq{}}\PY{p}{,} \PY{n}{argv}\PY{p}{[}\PY{l+m+mi}{0}\PY{p}{]}\PY{p}{)}\PY{p}{;}
        \PY{k}{return} \PY{l+m+mi}{0}\PY{p}{;}
    \PY{p}{\PYZcb{}}
    \PY{n}{sscanf}\PY{p}{(}\PY{n}{argv}\PY{p}{[}\PY{l+m+mi}{1}\PY{p}{]}\PY{p}{,} \PY{l+s}{\PYZdq{}}\PY{l+s}{\PYZpc{}lu}\PY{l+s}{\PYZdq{}}\PY{p}{,} \PY{o}{\PYZam{}}\PY{n}{numerator}\PY{p}{)}\PY{p}{;}
    \PY{n}{sscanf}\PY{p}{(}\PY{n}{argv}\PY{p}{[}\PY{l+m+mi}{2}\PY{p}{]}\PY{p}{,} \PY{l+s}{\PYZdq{}}\PY{l+s}{\PYZpc{}lu}\PY{l+s}{\PYZdq{}}\PY{p}{,} \PY{o}{\PYZam{}}\PY{n}{denominator}\PY{p}{)}\PY{p}{;}

    \PY{n}{printf}\PY{p}{(}\PY{l+s}{\PYZdq{}}\PY{l+s}{numerator: \PYZpc{}lu   denominator: \PYZpc{}lu}\PY{l+s+se}{\PYZbs{}n}\PY{l+s+se}{\PYZbs{}n}\PY{l+s}{\PYZdq{}}\PY{p}{,} \PY{n}{numerator}\PY{p}{,} \PY{n}{denominator}\PY{p}{)}\PY{p}{;}
    
    \PY{n}{repeatedly\PYZus{}construct\PYZus{}continued\PYZus{}fraction}\PY{p}{(}\PY{n}{numerator}\PY{p}{,} \PY{n}{denominator}\PY{p}{,} \PY{n}{sequence}\PY{p}{,} \PY{o}{\PYZam{}}\PY{n}{sequence\PYZus{}size}\PY{p}{,} \PY{n}{iterations}\PY{p}{)}\PY{p}{;}
    \PY{n}{repeatedly\PYZus{}evaluate\PYZus{}continued\PYZus{}fraction}\PY{p}{(}\PY{n}{numerator}\PY{p}{,} \PY{n}{denominator}\PY{p}{,} \PY{n}{sequence}\PY{p}{,} \PY{o}{\PYZam{}}\PY{n}{sequence\PYZus{}size}\PY{p}{,} \PY{n}{iterations}\PY{p}{)}\PY{p}{;}

    \PY{n}{repeatedly\PYZus{}construct\PYZus{}natural\PYZus{}representation}\PY{p}{(}\PY{n}{numerator}\PY{p}{,} \PY{n}{denominator}\PY{p}{,} \PY{n}{sequence}\PY{p}{,} \PY{o}{\PYZam{}}\PY{n}{sequence\PYZus{}size}\PY{p}{,} \PY{n}{iterations}\PY{p}{)}\PY{p}{;}
    \PY{n}{repeatedly\PYZus{}evaluate\PYZus{}natural\PYZus{}representation}\PY{p}{(}\PY{n}{numerator}\PY{p}{,} \PY{n}{denominator}\PY{p}{,} \PY{n}{sequence}\PY{p}{,} \PY{o}{\PYZam{}}\PY{n}{sequence\PYZus{}size}\PY{p}{,} \PY{n}{iterations}\PY{p}{)}\PY{p}{;}
\PY{p}{\PYZcb{}}
\end{Verbatim}

}

\pagebreak

\end{document}